\input ./style/arxiv-general.cfg
\documentclass[MSNbibl,number,citesort,rotating,dvips]{arxbj}
\makeatletter
   \@ifpackageloaded{graphicx}{}{\usepackage{graphicx}}
\makeatother
\usepackage{upgreek}
\usepackage{accents}

\aid{0}
\volume{21}
\issue{3}
\pubyear{2015}
\firstpage{1629}
\lastpage{1669}
\doi{10.3150/14-BEJ616} % kopijuoti is 'New paper accepted'
\makeatletter
\newcommand{\sdbullet}{%
  \hbox{\fontfamily{lmr}\fontsize{4}{0}\selectfont$\circ$}}
\newcommand{\rrVert}{\Vert}
\newcommand{\rrvert}{\vert}
\newcommand{\llVert}{\Vert}
\newcommand{\llvert}{\vert}
\newtheorem{theorem}{Theorem}[section]
\newproclaim{Assumption}{Assumption}

\newtheorem{corollary}[theorem]{Corollary}

\newproclaim{definition}[theorem]{Definition}
\newproclaim{example}[theorem]{Example}
\newproclaim{exercise}[theorem]{Exercise}
\newtheorem{lemma}[theorem]{Lemma}
\newtheorem{Lemma}{Lemma}[section]
\newproclaim{notation}[theorem]{Notation}
\newproclaim{problem}[theorem]{Problem}

\newremark{remark}{Remark}

\makeatother

\begin{document}
\begin{frontmatter}

\title{Robust estimation and inference for heavy tailed GARCH}
\runtitle{Robust estimation and inference for heavy tailed GARCH}

\begin{aug}
%%%% inicialai - be tarpu
\author[1]{\inits{J.B.}\fnms{Jonathan B.}~\snm{Hill}\corref{}\thanksref{1}\ead[label=e1]{jbhill@email.unc.edu}\ead[label=u1,url]{www.unc.edu/\textasciitilde jbhill}}% \and
%\author{\inits{}\fnms{}~\snm{}\thanksref{}\ead[label=e2]{}}
%\author{\inits{}\fnms{}~\snm{}}
%%\runauthor{} %% auto
%\dedicated{}
\address[1]{Department of Economics, University of North
Carolina, Chapel Hill, NC  27599-3305, USA.\\ \printead{e1,u1}}
%\address[]{}
\end{aug}

% HISTORY:
\received{\smonth{12} \syear{2012}}
\revised{\smonth{2} \syear{2014}}

% ABSTRACT
%
\begin{abstract}
We develop two new estimators for a general class of stationary
GARCH models
with possibly heavy tailed asymmetrically distributed errors, covering
processes with symmetric and asymmetric feedback like GARCH,
Asymmetric
GARCH, VGARCH and Quadratic GARCH. The first estimator arises from
negligibly trimming QML criterion equations according to error
extremes. The
second imbeds negligibly transformed errors into QML score equations
for a
Method of Moments estimator. In this case, we exploit a sub-class of
redescending transforms that includes tail-trimming and functions
popular in
the robust estimation literature, and we re-center the transformed
errors to
minimize small sample bias. The negligible transforms allow both
identification of the true parameter and asymptotic normality. We
present a
consistent estimator of the covariance matrix that permits classic inference
without knowledge of the rate of convergence. A simulation study shows both
of our estimators trump existing ones for sharpness and approximate
normality including QML, Log-LAD, and two types of non-Gaussian QML
(Laplace
and Power-Law). Finally, we apply the tail-trimmed QML estimator to
financial data.
\end{abstract}

% KEYWORDS
% visi is mazosios raides ir pagal abecele
%
\begin{keyword}
\kwd{GARCH}
\kwd{heavy tails}
\kwd{QML}
\kwd{robust inference}
\kwd{tail trimming}
\end{keyword}

\end{frontmatter}

%s1 #&#
\section{Introduction}

It is now widely accepted that log-returns of many macroeconomic and
financial time series are heavy tailed, exhibit clustering of large values,
and are asymmetrically distributed. In broader contexts extremes are
encountered in actuarial, meteorological, and telecommunication network
data %
(e.g., Leadbetter \textit{et~al}. \cite{leadbetteretal83}, Embrehts
\textit{et~al}. \cite{Embetal97},
Davis \cite{Davis2010}), while GARCH-type
clustering alone implies higher moments do not exist due to Pareto-like
distribution tails (e.g., Basrak \textit{et~al}. \cite{Basrak2002}, Liu
\cite{Liu2006}).

We develop new methods of robust estimation for a general class of
$\operatorname{GARCH}(1,1)$ models:%
%e1 #&#
\begin{equation}
y_{t}=\sigma_{t}\epsilon_{t}\quad\quad\mbox{with }\sigma
_{t}^{2}=g\bigl(y_{t-1},\sigma
_{t-1}^{2},\theta^{0}\bigr)\geq0\ \mbox{a.s.},
\label{garch}
\end{equation}
where $g(y,\sigma^{2},\theta)$ is a known mapping $g \dvtx  \mathbb{R}
\times [0,\infty) \times\Theta \rightarrow [0,\infty)$ and $%
\Theta$ is a compact subset of $\mathbb{R}^{q}$ for some finite $q
\geq 1$. We assume there exists a unique point $\theta^{0}$ in the
interior of $%
\Theta$ such that $\epsilon_{t} = y_{t}/\sigma_{t}$ is i.i.d. with a
non-degenerate absolutely continuous distribution with support $(-\infty
,\infty)$, $E[\epsilon_{t}] = 0$ and $E[\epsilon_{t}^{2}] = 1$.
Further, $\{y_{t},\sigma_{t}\}$ are stationary and geometrically
$\beta$-mixing. We avoid well known boundary problems by assuming
$\theta^{0}$
lies in the interior of $\Theta$ and $\sigma_{t}^{2}$ has a non-degenerate
distribution, hence (\ref{garch}) is a non-trivial GARCH process. In
Bollerslev's \cite{bollerslev86} classic GARCH model $\sigma_{t}^{2} =
\omega^{0} +
\alpha^{0}y_{t-1}^{2} + \beta^{0}\sigma_{t-1}^{2}$, with $\omega
^{0} > 0$ and $\alpha^{0},\beta^{0} \geq 0$ this requires $%
\alpha^{0} + \beta^{0} > 0$, cf. Andrews \cite{Andrews99} and Francq
and Zako\"{\i}an \cite{FZ04}.

In order to keep technical arguments brief, we assume $\sigma
_{t}^{2}(\theta) := g(y_{t-1},\sigma_{t-1}^{2}(\theta),\theta)$ has
properties similar to a non-trivial classic GARCH model: $\sigma
_{t}^{2}(\theta)$ is twice continuously differentiable, $E[(\sup_{\theta
\in\Theta}|\sigma_{t}^{2}/\sigma_{t}^{2}(\theta)|)^{p}] < \infty$
for any $p > 0$, and $\sup_{\theta\in\mathcal{N}_{0}}\|(\partial
/\partial\theta)^{i}\ln(\sigma_{t}^{2}(\theta))\|$ is $L_{2+\iota
}$-bounded for tiny $\iota > 0$ and some compact $\mathcal{N}_{0}
\subseteq \Theta$ containing $\theta^{0}$, where $\|\cdot\|$ is the
matrix norm (cf. Francq and Zako\"{\i}an \cite{FZ04}). Similarly, we
impose Lipschitz type bounds
on $g$ that ensure an iterated approximation $h_{0}^{2}(\theta) =
\omega$ and $h_{t}(\theta) = g(y_{t-1},h_{t-1}(\theta),\theta)$ for
$%
t = 1,2,\ldots $ satisfies $\sup_{\theta\in\Theta}|h_{t}(\theta) -
\sigma_{t}^{2}(\theta)| \stackrel{p}{\rightarrow} 0$ as $t
\rightarrow \infty$, a key property for feasible estimation %
(see Nelson \cite{Nelson90}, Francq and Zako\"{\i}an \cite{FZ04},
Straumann and Mikosch \cite{StraumannMikosch}). The above properties of
$%
\sigma_{t}^{2}(\theta)$ cover at least Threshold GARCH with a known
threshold, Asymmetric and Nonlinear Asymmetric GARCH, VGARCH, GJR-GARCH,
Smooth Transition GARCH, and Quadratic GARCH. Consult Engle and Ng \cite
{engleng93},
Carrasco and Chen \cite{carrascochen02}, Francq and Zako\"{\i}an \cite
{FZ04,FZ10} and %
Meitz and Saikkonen \cite{meitzsaikkonen08,meitzsaik11}. EGARCH
evidently is not included
here since it is unknown whether $\sup_{\theta\in\Theta
}|h_{t}(\theta)
- \sigma_{t}^{2}(\theta)| \stackrel{p}{\rightarrow} 0$ as $t
\rightarrow \infty$ %
(see Straumann and Mikosch \cite{StraumannMikosch}, Meitz and Saikkonen
\cite{meitzsaikkonen08,meitzsaik11}).

We are interested in heavy tailed errors or innovation outliers, in
particular we allow \mbox{$E[\epsilon_{t}^{4}] = \infty$}, while GARCH
feedback itself may also prompt heavy tails in $y_{t}$ due to a stochastic
recurrence structure (Basrak \textit{et~al}. \cite{Basrak2002}, Liu
\cite{Liu2006}). In this paper, we negligibly
transform QML loss or score equations to obtain asymptotically normal
estimators of $\theta^{0}$ allowing for $E[\epsilon_{t}^{4}] = \infty
$%
.

Define $\epsilon_{t}(\theta) := y_{t}/\sigma_{t}(\theta)$ and $%
\mathfrak{s}_{t}^{2}(\theta) : =(\partial/\partial\theta)\ln\sigma
_{t}^{2}(\theta)$, and let $I(\cdot)$ denote the indicator function. In
Section~\ref{sec:qmttl}, we tackle the fact that $\sigma
_{t}^{2}(\theta)$
is not observed for $t \leq 0$. The first method trims QML criterion
equations $p_{t}(\theta) := \ln(\sigma_{t}^{2}(\theta)) +
\epsilon_{t}^{2}(\theta)$ according to extremes that arise in a first
order expansion and therefore the score $\sum_{t=1}^{n}(\epsilon
_{t}^{2}(\theta) - 1)\mathfrak{s}_{t}^{2}(\theta)$. Since $\mathfrak
{s}%
_{t}^{2}(\theta)$ has an $L_{2}$-bounded envelope near $\theta^{0}$ it
suffices to minimize $\sum_{t=1}^{n}p_{t}(\theta)I(-l \leq \epsilon
_{t}^{2}(\theta) - 1 \leq u)$ for some positive thresholds $%
\{l,u\} $ that increase with the sample size $n$. Identification of
$\theta
^{0}$ coupled with asymptotic normality are assured if $\{l,u\}$ are
replaced with intermediate order statistics of $\epsilon
_{t}^{2}(\theta)
- 1$. The result is the Quasi-Maximum Tail-Trimmed Estimator (QMTTL),
similar to the least tail-trimmed squares estimator for autoregressions in
Hill \cite{Hillltts}.

The second method imbeds negligibly transformed errors in QML score
equations $(\epsilon_{t}^{2}(\theta) - 1)\mathfrak{s}_{t}^{2}(\theta
)$. We then re-center the transformed errors to minimize small sample
bias and
estimate $\theta^{0}$ by the Method of Negligibly Weighted Moments (MNWM).
By re-centering we may simply transform $\epsilon_{t}(\theta)$ itself
symmetrically which requires only one threshold, for example in the simple
trimming case we use $\epsilon_{t}^{2}(\theta)I(|\epsilon_{t}(\theta)|
\leq c)$ for some $c > 0$. In order to simplify proofs we focus on
simple trimming, and related bounded but smooth weighted redescending
transforms $\epsilon_{t}^{2}(\theta)\varpi(\epsilon_{t}^{2}(\theta
),c)I(|\epsilon_{t}(\theta)| \leq c)$ where $\varpi(\cdot,c)$ is
continuously differentiable in $c$, and $\varpi(\epsilon
_{t}^{2}(\theta
),c) \rightarrow 1$ a.s. as $c \rightarrow \infty$. Weights
related to simple indicators include Hampel's three-part function, and
smooth transforms include Tukey's bisquare and an exponential version %
(cf. Andrews \textit{et~al}. \cite{Andrewsetal72}, Hampel \textit
{et~al}. \cite{Hampleetal86}). See Sections~\ref{sec:qmttl}
and \ref{sec:mm}.

We show how trimming and distribution tail parameters impact efficiency,
while the negligible amount of trimming never affects the asymptotic
covariance matrix when $E[\epsilon_{t}^{4}] < \infty$. Fixed quantile
trimming or truncation always impact efficiency irrespective of higher
moments, and cause bias due to $\epsilon_{t}^{2} - 1$ having an
asymmetric distribution in general (Sakata and White \cite
{SakataWhite98}, Mancini \textit{et~al}. \cite{Mancinietal05}).
Mancini \textit{et~al}. \cite{Mancinietal05} use simulation based
methods to solve the bias, but
this requires knowledge of the error distribution
(see also Cantoni and Ronchetti \cite{CantoniRonchetti}, Ronchetti and
Trojani \cite{ronchettitrojani01}).

The convergence rate of our estimators is $\mathrm{o}(\sqrt{n})$ when $E[\epsilon
_{t}^{4}] = \infty$, but can be assured to be $\sqrt{n}/g_{n}$ for any
sequence of positive numbers $\{g_{n}\}$ that satisfies $g_{n}
\rightarrow
\infty$ as slowly as we choose by following simple rules of thumb for
choosing the threshold $c$. Thus when $E[\epsilon_{t}^{4}] = \infty$
our estimators converge faster than QML (cf. Hall and Yao \cite
{hallyao03}) but
slower than $\sqrt{n}$-convergent estimators in Peng and Yao \cite
{pengyao03}, Berkes and Horvath \cite{BerkesHorvath04} and Zhu and Ling
\cite{ZhuLing}, although the latter two are not for
standard GARCH models in which $E[\epsilon_{t}^{2}] = 1$ identifies the
volatility process. See below for literature details. We do not tackle
\textit{optimal} threshold selection in order to conserve space. We do,
however, show explicitly how threshold selection impacts the convergence
rate which suggests simple rules for trimming. We also discuss practical
considerations for trimming in terms of small sample bias control. See
Sections~\ref{ver_2} and \ref{sec:rate}.

In Section~\ref{sec:infer}, we show classic inference applies as long as self-normalization
is used, a nice convenience since tail thickness and the precise rate of
convergence need never be known. We complete the paper with simulation and
empirical studies in Sections~\ref{sim} and \ref{s6}. In particular, we give evidently the
first comparison of various heavy tail robust estimators for GARCH models,
and show our estimators obtain in general lower bias and are closer to
normally distributed in small samples and therefore lead to better inference.

A complete theory of QML for a variety of strong-GARCH models is presented
in Lee and Hansen \cite{leehansen94}, Berkes \textit{et~al}. \cite
{Berkesal03}, Francq and Zako\"{\i}an \cite{FZ04}, Straumann and
Mikosch \cite{StraumannMikosch} and Meitz and Saikkonen \cite
{meitzsaik11} amongst others, while at least a
finite fourth moment $E[\epsilon_{t}^{4}] < \infty$ is standard. The
allowance of heavier tails $E[\epsilon_{t}^{4}] = \infty$, with
Gaussian asymptotics, evidently only exists for the classic GARCH
model, and
in most cases requires a non-Gaussian QML criterion and non-standard moment
conditions to ensure Fischer consistency (i.e., consistency for the true
parameter $\theta^{0}$). Peng and Yao \cite{pengyao03} propose $\sqrt {n}$-convergent
Log-LAD, requiring $\ln\epsilon_{t}^{2}$ to have a zero median in
order to
identify $\theta^{0}$. Berkes and Horvath \cite{BerkesHorvath04}
characterize a general QML
criterion class that potentially allows for Fischer consistency, $\sqrt {n}$%
-convergence and asymptotic normality even when $E[\epsilon_{t}^{4}] =
\infty$. They treat Gaussian QML, and various non-Gaussian QML like Laplace
QML which requires $E|\epsilon_{t}| = 1$ and $E[\epsilon_{t}^{2}] <
\infty$, and Power-Law QML (PQML) with index $\vartheta > 1$
requiring that $\epsilon_{t}$ have an infinitessimal moment and $%
E[|\epsilon_{t}|/(1+|\epsilon_{t}|)] = 1/\vartheta$. Student's $t$%
-QML is Fischer consistent when $\epsilon_{t}$ is $t$-distributed, and
otherwise may only be consistent for some $\tilde{\theta} \neq \theta
^{0}$ (cf. Newey and Steigerwald \cite{NeweySteigerwald97}, Sakata and
White \cite{SakataWhite98}, Fan \textit{et~al}. \cite{Fanetal12}).

Zhu and Ling \cite{ZhuLing} combine Berkes and Horvath \cite
{BerkesHorvath04} Laplace class with Ling's \cite{ling07} weighting
method for Weighted Laplace QML (WLQML) under the
assumptions $\epsilon_{t}$ has a zero median, $E|\epsilon_{t}| = 1$
and $E[\epsilon_{t}^{2}] < \infty$. The estimator is $\sqrt{n}$%
-convergent and asymptotically normal when $E[\epsilon_{t}^{4}] =
\infty$, but the suggested weights at time $t$ are based on the infinite
past $y_{t-1},y_{t-2},\ldots $\,. Although the authors use a central order
statistic for a threshold and fix $y_{t} = 0$ for $t \leq 0$ in the
weights for the sake of simulations, they do not prove either is valid.
Indeed, for a $\operatorname{GARCH}(1,1)$ the restriction $y_{t} = 0$ for $t \leq 0$
in their weight (2.4) does not support asymptotic normality
(see Zhu and Ling \cite{ZhuLing}, Assumption~2.4 and the discussion on
weight (2.4)). Thus, the estimator is not
evidently feasible.

Assumptions like $E|\epsilon_{t}| = 1$ or $E[|\epsilon
_{t}|/(1+|\epsilon_{t}|)] = 1/\vartheta$ replace the usual $%
E[\epsilon_{t}^{2}] = 1$ to identify $\theta^{0}$. Of course, if $%
E[\epsilon_{t}^{2}] \neq1$ then model (\ref{garch}) is not a standard
GARCH model since $E[y_{t}^{2}|y_{t-1},y_{t-2},\ldots] \neq \sigma
_{t}^{2} $ with positive probability is possible, and Gaussian QML
leads to
asymptotic bias. Thus, asymptotic normality \textit{and} Fischer consistency
are assured precisely by changing the criterion \textit{and} model
assumptions and therefore the model by imposing a non-standard moment
condition. In practice, this may be untenable as many analysts in economics
and finance first impose a version of (\ref{garch}) with $E[\epsilon
_{t}^{2}] = 1$ and then seek a robust estimator. In order to sidestep
such unpleasant moment conditions, Fan \textit{et~al}. \cite{Fanetal12}
introduce a three-step
non-Gaussian QML method. In the first stage, Gaussian QML residuals are
generated. In a second stage, a scale parameter is estimated to ensure
identification in the third non-Gaussian QML stage without imposing
non-standard moment conditions. See also Newey and Steigerwald \cite
{NeweySteigerwald97}. Our
QMTTL and MNWM estimators are computed in one-step and are asymptotically
normal and Fischer consistent by imposing negligible weighting on extremes
couched in a Gaussian QML criterion.

Evidently simulation experiments demonstrating the robustness
properties of
Peng and Yao's \cite{pengyao03} Log-LAD, Berkes and Horvath's \cite
{BerkesHorvath04} non-Gaussian QML and
Zhu and Ling's \cite{ZhuLing} WLQML does not exist, while Fan \textit
{et~al}. \cite{Fanetal12} only inspect
the root-mean squared error of their estimator which masks possible
bias. In
general, the empirical bias and approximate normality properties of these
estimators, as well as their ability to gain accurate inference in small
samples (e.g., Wald tests), are unknown.

In a simulation experiment, we show QMTTL and MNWM trump QML, Log-LAD, WLQML,
and PQML in all cases in terms of bias, approximate normality and $t$-test
performance, and has lower mean-squared-error than every estimator except
PQML (PQML has higher bias and lower dispersion). Overall QMTTL performs
best. The dominant performance of QMTTL and MNWM follows since only they
directly counter the influence of large errors in small \textit{and} large
samples by trimming observations with an error extreme. We show this matters
even when $\epsilon_{t}$ is Gaussian: negligible trimming always improves
QML performance, while untrimmed QML, Log-LAD, WLQML and PQML are
comparatively more sensitive to large errors. Moreover, even PQML,
which we
design as in Berkes and Horvath \cite{BerkesHorvath04} to ensure
identification for Paretian
errors with an infinite fourth moment, has greater bias and is farther from
normality in small samples than QMTTL and MNWM. Thus, the advantages of
non-Gaussian QML for GARCH processes with heavy tailed errors are not clear,
at least as seen by our controlled experiments. We emphasize this last point
by tail-trimming PQML in a way that removes adverse sample extremes and
leaves the estimator asymptotically unbiased. We show in most cases tail
trimming helps PQML in terms of bias, approximate normality and inference,
yet overall QMTTL is still better. Indeed, PQML is infeasible unless the
tail index of $\epsilon_{t}$ is known or estimated using some filtration
for $\epsilon_{t}$ (e.g., QML residuals), and is not Fischer
consistent if $%
\epsilon_{t}$ has any other distribution.

In the literature on additive outlier robust estimation, negligible trimming
is an example of a \textit{redescending} transformation $\psi \dvtx  \mathbb
{%
R} \rightarrow\mathbb{R}$ where in general $\psi(u) \rightarrow
0 $ as $|u| \rightarrow \infty$, and typically $\psi(u) = 0$ when
$|u| > c$ for some $c$ as we use here. See Huber \cite{Huber1964} and
Hampel \textit{et~al}.~\cite{Hampleetal86}. Evidently a complete theory
of redescending M-estimators
exists only for estimates of location for i.i.d. data %
(Shevlyakov and Shurygin \cite{Shevlyakovetal08}). In this paper, our
QML estimator has a score
equation that effectively uses $\psi(\epsilon_{t}) = (\epsilon
_{t}^{2} - 1)I(-l \leq \epsilon_{t}^{2} - 1 \leq u)$
where $l,u \rightarrow \infty$ as \mbox{$n \rightarrow \infty$}. Our
Method of Moments estimator is more generic since it uses either re-centered
$\psi(\epsilon_{t}) = \epsilon_{t}^{2}I(|\epsilon_{t}| \leq c)$
with $c \rightarrow \infty$ as $n \rightarrow \infty$, or
related variants like Hampel's three-part weight, as well as smooth weights
like Tukey's bisquare. In all cases, the increasing thresholds ensure
bias is
eradicated asymptotically.

We ignore additive or isolated outliers, and so-called one-off events
in $%
\{y_{t}\}$ for the sake of brevity. In this case, we would observe
$y_{t} =
y_{t}^{\ast} + x_{t}$ where $y_{t}^{\ast}$ is generated by (\ref%
{garch}) and, for example, $x_{t} = 0$ in most periods $t$. The
challenge here is controlling the propagation of an aberrant
observation due
to $x_{t} \neq 0$ through the volatility mechanism. See, for example,
Charles and Darn\'{e} \cite{CharlesDarne}, Muler and Yohai \cite
{MulerYohai08}, and Boudt \textit{et~al}. \cite{Boudtetal2011}, and
see Mendes \cite{Mendes00} for anecdotal evidence of QML estimator bias.
Incorporating additive outliers in (\ref{garch}) with innovation outliers
would require additional robustness techniques like those employed in these
and related papers (e.g., Muler \textit{et~al}. \cite{Mulleretal09}).
Some methods, however,
are proposed to detect outliers in a GARCH process under the assumption of
thin tailed errors: a few large values are simply assumed to be due to a
non-heavy tailed outlier.\footnote{Charles and Darn\'{e} \cite
{CharlesDarne} extend ideas
developed in \citeauthor{ChenLiu93}  to test for, and
control, additive and
innovation outliers in a GARCH process with Gaussian errors. These
papers do
not provide asymptotic theory, hence the Gaussian assumption can likely be
relaxed. The trimming methods used in the present paper can be extended to
their test statistics which involve a residual variance estimator %
(cf. Hill \cite{Hillltts}, Hill and Aguilar \cite{HillAguilar13}), but
a rigorous theory would need to
be developed.} Other estimators, contrary to claims, do not identify
$\theta
^{0}$ and/or are not robust to heavy tailed errors.\footnote{Muler and
Yohai \cite{MulerYohai08} present a robust M-estimator
$\accentset{\sdbullet}{\theta}_{n} =
\arg\inf_{\theta\in\Theta}\{\sum_{t=1}^{n}\rho(\ln
(y_{t}^{2}/h_{t}^{\ast2}(\theta))\}$ where $h_{t}^{\ast2}(\theta)$
is a
filtered version of $\sigma_{t}^{2}(\theta)$ that restricts the
propagation of outliers. They assume $\rho$ is thrice continuously
differentiable with bounded derivatives. Although claimed to be heavy tail
robust and identify the true $\theta^{0}$ (see their Theorem~3), they do
not prove any such $\rho$ exists. In their simulations, for example, they
use truncated QML with $\rho(u) = \psi_{c}(\exp\{u\} - u)$ where $%
\psi_{c}$ truncates at a fixed threshold $c \dvt \psi_{c}(x) = K$ for
all $x > c$. Thus $\rho(u)$ is non-differentiable at $\exp\{u\} -
u = c$, and at all other points no derivative is bounded which implies
non-robustness to heavy tails. The problem is the QML score is not bounded
when $\rho(u)$ is truncated according to its large values. Our approach,
however, negligibly trims according to properties of the QML score and
therefore ensures heavy tail robustness \textit{and} identification of $
\theta^{0}$.} Further, all such robust estimators are proposed for the
classic GARCH model, hence existing theory does not necessarily extend to
the broader model class (\ref{garch}).

Finally, our methods can be easily extended to higher order GARCH models,
GARCH-in-Mean, and models of the conditional mean and variance like
nonlinear ARMA--GARCH, as well as other estimators like non-Gaussian QML
(Berkes and Horvath \cite{BerkesHorvath04}, Zhu and Ling \cite
{ZhuLing}, Fan \textit{et~al}. \cite{Fanetal12}), LAD
(Peng and Yao \cite{pengyao03}), etc.
We show trimming matters for PQML in our simulation study, and we expect
negligible trimming to improve upon non-Gaussian QML estimators in general,
provided they are Fischer consistent in the first place.

We use the following notation conventions. The indicator function
$I(\cdot)$
is $I(a) = 1$ if $a$ is true, and otherwise $I(a) = 0$. The spectral
norm of matrix $A$ is $\|A\| = \lambda_{\max}(A^{\prime}A)^{1/2}$
with $\lambda_{\max}(\cdot)$ the maximum eigenvalue. If $z$ is a scalar,
we write $(z)_{+} := \max\{0,z\}$. $K$ denotes a positive finite
constant whose value may change from line to line; $\iota> 0$ is an
arbitrarily tiny constant. $\stackrel{p}{\rightarrow}$ and $\stackrel
{d}{%
\rightarrow}$ denote probability and distribution convergence. $x_{n}
\sim a_{n}$ implies $x_{n}/a_{n} \rightarrow 1$. $L(n)$ is a slowly
varying function that may change with the context.

%s2 #&#
\section{Quasi-maximum tail-trimmed likelihood}\label{sec:qmttl}

The observed sample is $\{y_{t}\}_{t=0}^{n}$ with sample size $n + 1
\geq 1$. We start at $t = 0$ to simplify notation since we condition
on the first observation $y_{0}$ and a volatility constant defined below.
Estimation requires a volatility function on $\Theta$,
\[
\sigma_{t}^{2}(\theta)=g \bigl( y_{t-1},
\sigma_{t-1}^{2}(\theta ),\theta \bigr),
\]
hence $\sigma_{t}^{2} = \sigma_{t}^{2}(\theta^{0})$. It is convenient
to assume $\Theta$ is a compact subset of points $\theta$ on which
$\sigma
_{t}^{2}(\theta)$ is stationary:%
%e2 #&#
\begin{equation}
\Theta\subseteq \bigl\{ \theta\in\mathbb{R}^{q}\dvtx  \bigl\{ \sigma
_{t}^{2}(\theta) \bigr\} \mbox{ has a stationary solution}
\bigr\} . \label{theta}
\end{equation}

In practice $\sigma_{t}^{2}(\theta)$ for $t \leq 0$ is not observed,
so define an iterated volatility approximation%
%e3 #&#
\begin{equation}
h_{0}(\theta)=\tilde{\omega}>0\quad\mbox{and}\quad h_{t}(\theta)=g
\bigl( y_{t-1},h_{t-1}(\theta),\theta \bigr) \quad\quad\mbox{for
}t=1,2,\ldots, \label{h_iter}
\end{equation}
where $\tilde{\omega}$ is not necessarily an element of $\theta$. We
initially develop an infeasible robust estimator based on the QML
equations $%
\ln\sigma_{t}^{2}(\theta) + y_{t}^{2}/\sigma_{t}^{2}(\theta)$. We
then show a feasible version based on $\ln h_{t}(\theta) +
y_{t}^{2}/h_{t}(\theta)$ has the same limit distribution.

%s2.1 #&#
\subsection{Tail-trimming}

In order to understand when and where trimming should be applied,
define the
GARCH error function, and a scaled volatility function and its
derivative
\begin{eqnarray*}
 \epsilon_{t}(\theta)&:=&\frac{y_{t}}{\sigma_{t}(\theta)}=\frac
{y_{t}}{%
g(y_{t-1},\sigma_{t-1}^{2}(\theta),\theta)},
\\
 \mathfrak{s}_{t}(\theta)&=& \bigl[ \mathfrak{s}_{i,t}(
\theta) \bigr] _{i=1}^{q}:=\frac{1}{\sigma_{t}^{2}(\theta)}\frac{\partial}{\partial
\theta}
\sigma_{t}^{2} ( \theta )\quad \mbox{and}\quad\mathfrak{d}%
_{t}(\theta)= \bigl[ \mathfrak{d}_{i,j,t}(\theta) \bigr]
_{i,j=1}^{q}:= %%
\frac{\partial}{\partial\theta}
\mathfrak{s}_{t}(\theta).
\end{eqnarray*}
Throughout, we drop $\theta^{0}$ and write $\epsilon_{t} = \epsilon
_{t}(\theta^{0})$, $\mathfrak{s}_{t} = \mathfrak{s}_{t}(\theta^{0})$,
$%
\mathfrak{d}_{t} = \mathfrak{d}_{t}(\theta^{0})$ and so on. Gaussian
asymptotics for QML are grounded on the score equations $m_{t}(\theta
)$ and
their Jacobian $G_{t}(\theta)$:
%e4 #&#
\begin{eqnarray}\label{mG}
m_{t}(\theta)&:=& \bigl( \epsilon_{t}^{2}(
\theta)-1 \bigr) \mathfrak{s}%
_{t}(\theta)\quad\mbox{and}\nonumber
\\[-8pt]\\[-8pt]
G_{t}(\theta)&:=&\frac{\partial}{\partial
\theta}m_{t}(\theta)= \bigl(
\epsilon_{t}^{2}(\theta)-1 \bigr) \mathfrak{d%
}_{t}(\theta)-\epsilon_{t}^{2}(\theta)
\mathfrak{s}_{t}(\theta )\mathfrak{s%
}_{t}(
\theta)^{\prime}. \nonumber
\end{eqnarray}

We assume $\mathfrak{s}_{t}(\theta)$ and $\mathfrak{d}_{t}(\theta)$
have $%
L_{2+\iota}$-bounded envelopes near $\theta^{0}$ for tiny $\iota > 0$, thus asymptotic normality hinges entirely on $\epsilon_{t}^{2} - 1$.
See below for all assumptions. It therefore suffices to trim $\ln\sigma
_{t}^{2}(\theta) + \epsilon_{t}^{2}(\theta)$ negligibly when $%
\epsilon_{t}^{2}(\theta) - 1$ surpasses a large negative or positive
threshold. As long as those thresholds represent intermediate order
statistics, we can identify $\theta^{0}$ and have an asymptotically normal
estimator. Write
\[
\mathcal{E}_{t}(\theta):=\epsilon_{t}^{2}(
\theta)-1,
\]
and denote left and right tail observations and their order statistics
for $%
\mathcal{E}_{t}(\theta)$:
\begin{eqnarray*}
 \mathcal{E}_{t}^{(-)}(\theta)&:=&\mathcal{E}_{t}(
\theta)I \bigl( \mathcal{E}%
_{t}(\theta)<0 \bigr) \quad\mbox{and}\quad\mathcal{E}%
_{(1)}^{(-)}(\theta)\leq
\cdots\leq\mathcal{E}_{(n)}^{(-)}(\theta )\leq0,
\\
 \mathcal{E}_{t}^{(+)}(\theta)&:=&\mathcal{E}_{t}(
\theta)I \bigl( \mathcal{E}%
_{t}(\theta)\geq0 \bigr)\quad \mbox{and}\quad\mathcal{E}%
_{(1)}^{(+)}(\theta)\geq
\cdots\geq\mathcal{E}_{(n)}^{(+)}(\theta )\geq0.
\end{eqnarray*}
The determination of the number of trimmed large $\mathcal
{E}_{t}(\theta)$
in a sample of size $n$ is made by intermediate order sequences $%
\{k_{1,n},k_{2,n}\}$, hence (e.g., Leadbetter \textit{et~al}. \cite
{leadbetteretal83})
\[
k_{i,n}\in \{ 1,\ldots,n-1 \} ,\quad\quad k_{i,n}\rightarrow\infty
\quad\mbox{and}\quad k_{i,n}/n\rightarrow0.
\]
Define an indicator selection function for trimming
\[
\hat{I}_{n,t}^{(\mathcal{E})}(\theta):=I \bigl( \mathcal{E}%
_{(k_{1,n})}^{(-)}(\theta)\leq\mathcal{E}_{t}(\theta)\leq
\mathcal{E} %%
_{(k_{2,n})}^{(+)}(\theta) \bigr) .
\]
The QMTTL estimator therefore solves
\[
\hat{\theta}_{n}=\operatorname{arg\,min}\limits
_{\theta\in\Theta
} \Biggl\{
\frac{1}{%
n}\sum_{t=1}^{n} \bigl( \ln
\sigma_{t}^{2}(\theta)+\epsilon _{t}^{2}(
\theta ) \bigr) \times\hat{I}_{n,t}^{(\mathcal{E})}(\theta) \Biggr\} =
\operatorname{arg\,min}\limits
_{\theta\in\Theta} \bigl\{ \hat {Q}_{n} ( \theta ) \bigr
\} .
\]
Each $k_{i,n}$ represents the number of trimmed $\ln\sigma
_{t}^{2}(\theta
) + \epsilon_{t}^{2}(\theta)$ due to large negative or positive $%
\mathcal{E}_{t}(\theta) = \epsilon_{t}^{2}(\theta) - 1$. We
require $k_{i,n} \rightarrow \infty$ for asymptotic normality, while
negligibility $k_{i,n}/n \rightarrow 0$ ensures identification of $%
\theta^{0}$ asymptotically. Since $\mathcal{E}_{t}(\theta)$ in
general has
an asymmetric distribution, identification of $\theta^{0}$ is assured
asymptotically if we negligibly trim asymmetrically by $\mathcal{E}%
_{t}(\theta)$. In a method of moments framework, however, we can
re-centered trimmed errors allowing for \textit{symmetric} trimming where
negative and positive thresholds are the same: see Section~\ref{sec:mm}.

In practical terms, $\hat{\theta}_{n}$ can be easily computed using standard
iterative optimization routines. In fact, under distribution continuity
arguments developed in Cizek \cite{cizek08}, Lemma~2.1, page 29, apply
for almost sure twice differentiability of the otherwise
non-differentiable $%
\hat{Q}_{n}(\theta)$. In particular, we have \textit{almost surely} $%
(\partial/\partial\theta)\hat{Q}_{n}(\theta) = 1/n%
\sum_{t=1}^{n}m_{t}(\theta)\hat{I}_{n,t}^{(\mathcal{E})}(\theta)$ and
$%
(\partial/\partial\theta)^{2}\hat{Q}_{n}(\theta) =
1/n\sum_{t=1}^{n}G_{t}(\theta)\hat{I}_{n,t}^{(\mathcal{E})}(\theta)$. This
implies standard estimation algorithms that exploit the gradient and Hessian
apply.

In order to characterize the limit distribution of $\hat{\theta}_{n}$, we
require non-random quantiles which the order statistics $\mathcal{E}%
_{(k_{1,n})}^{(-)}(\theta)$ and $\mathcal{E}_{(k_{2,n})}^{(+)}(\theta)$
approximate. Define sequences $\{\mathcal{L}_{n}(\theta),\mathcal{U}%
_{n}(\theta)\}$ denoting the lower $k_{1,n}/n$ and upper $k_{2,n}/n$
quantiles of $\mathcal{E}_{t}(\theta)$:%
%e5 #&#
\begin{equation}\label{ck}
P \bigl( \mathcal{E}_{t}(\theta)\leq-\mathcal{L}_{n}(
\theta) \bigr) =\frac{%
k_{1,n}}{n}\quad\mbox{and}\quad P \bigl( \mathcal{E}_{t}(
\theta)\geq\mathcal{U}%
_{n}(\theta) \bigr) =
\frac{k_{2,n}}{n}.
\end{equation}
The selection indicator is then
\[
I_{n,t}^{(\mathcal{E})}(\theta):=I \bigl( -\mathcal{L}_{n}(
\theta)\leq \mathcal{E}_{t}(\theta)\leq\mathcal{U}_{n}(
\theta) \bigr) .
\]
Notice $\mathcal{E}_{t}(\theta) \in [-1,\infty)$ and $k_{i,n}/n
\rightarrow 0$ imply $\mathcal{L}_{n}(\theta) \rightarrow 1$ and $%
\mathcal{U}_{n}(\theta) \rightarrow \infty$. The quantiles $\{%
\mathcal{L}_{n}(\theta)$, $\mathcal{U}_{n}(\theta)\}$ exist for each $
\theta$ and any choice of fractiles $\{k_{1,n},k_{2,n}\}$ since
$\epsilon
_{t}$ has a smooth distribution. By construction the order statistics
$\{%
\mathcal{E}_{(k_{1,n})}^{(-)}(\theta)$, $\mathcal{E}_{(k_{2,n})}^{(+)}(
\theta)\}$ estimate $\{\mathcal{L}_{n}(\theta),\mathcal{U}_{n}(\theta
)\}$%
, and are uniformly consistent in view of the $\beta$-mixing condition
detailed in Assumption~\ref{ass1} below, for example $\sup_{\theta\in\Theta}|\mathcal
{E}%
_{(k_{2,n})}^{(+)}(\theta)/\mathcal{U}_{n}(\theta) - 1| =
\mathrm{O}_{p}(1/k_{2,n}^{1/2})$. See Appendix~\ref{app:lemmas} for supporting limit
theory.

Finally, define equation variances $\Sigma_{n}$ and $\mathcal{S}_{n}$, and
a scale $\mathcal{V}_{n}$ for standardizing $\hat{\theta}_{n}$:
\begin{eqnarray*}
 \Sigma_{n}&:=&E \bigl[ \mathcal{E}_{t}^{2}I_{n,t}^{(\mathcal{E})}
\bigr] \times E \bigl[ \mathfrak{s}_{t}\mathfrak{s}_{t}^{\prime}
\bigr] \quad\mbox{and}\quad\mathcal{S}_{n}:=E \Biggl[ \Biggl(
\frac{1}{n^{1/2}}%
\sum_{t=1}^{n}m_{t}I_{n,t}^{(\mathcal{E})}
\Biggr) \Biggl( \frac
{1}{n^{1/2}}%
\sum_{t=1}^{n}m_{t}I_{n,t}^{(\mathcal{E})}
\Biggr) ^{\prime} \Biggr],
\\
 \mathcal{V}_{n}&=& [ \mathcal{V}_{i,j,n} ]
_{i,j=1}^{q}:=nE \bigl[ \mathfrak{s}_{t}
\mathfrak{s}_{t}^{\prime} \bigr] \mathcal {S}_{n}^{-1}E
\bigl[ \mathfrak{s}_{t}\mathfrak{s}_{t}^{\prime}
\bigr] \sim\frac{n}{E [
\mathcal{E}_{t}^{2}I_{n,t}^{(\mathcal{E})} ] }E \bigl[ \mathfrak {s}_{t}%
\mathfrak{s}_{t}^{\prime} \bigr] .
\end{eqnarray*}
The scale form $\mathcal{V}_{n} = nE[\mathfrak{s}_{t}\mathfrak{s}%
_{t}^{\prime}]\mathcal{S}_{n}^{-1}E[\mathfrak{s}_{t}\mathfrak{s}%
_{t}^{\prime}]$ is standard for M-estimators. In view of identification
Assumption~\ref{ass2} and equation (\ref{sigma}), below, and independence it is
easily verified that the long-run variance satisfies $\mathcal{S}_{n} =
\Sigma_{n}(1 + \mathrm{o}(1))$. Thus $\mathcal{V}_{n} \sim n(E[\epsilon
_{t}^{4}I_{n,t}^{(\mathcal{E})}] - 1)^{-1}E[\mathfrak{s}_{t}\mathfrak{s}
_{t}^{\prime}]$, which is positive definite for our data generating process.

%s2.2 #&#
\subsection{Main results}

We require two assumptions concerning the error distribution,
properties of
the volatility response $g$, and parameter identification. Let $\kappa$
denote the moment supremum of $\epsilon_{t}$:
\[
\kappa:=\arg\sup \bigl\{ \xi>0\dvtx E\llvert \epsilon_{t}\rrvert
^{\xi
}<\infty \bigr\} >2.
\]

\begin{Assumption}[(Data generating process)]\label{ass1}
\begin{enumerate}[(d)]
 \item[(a)] There exists a unique point $\theta^{0} = [\omega^{0},\alpha
^{0},\beta^{0}]^{\prime}$
in the interior
of a compact subset $\Theta$ of $\mathbb{R}^{q}$ such
that $\epsilon_{t} = y_{t}/\sigma_{t}$
is i.i.d., $E[\epsilon_{t}] = 0$ and $%
E[\epsilon_{t}^{2}] = 1$.

 \item[(b)] $\epsilon_{t}$ has an absolutely continuous, non-degenerate,
and uniformly bounded distribution on $(-\infty,\infty) :
\sup_{a\in\mathbb{R}}\{(\partial/\partial a)P(\epsilon_{t} \leq a)\}
< \infty$.
If $E[\epsilon_{t}^{4}] = \infty$
then $P(|\epsilon_{t}| > a) = da^{-\kappa}(1 + \mathrm{o}(1))$, where $d > 0$
and $\kappa \in (2,4]$.

 \item[(c)] $g(\cdot,\cdot,\theta)$ is twice continuously
differentiable in $\theta$; $(\partial/\partial\theta)^{i}g(\cdot
,\cdot,\theta)$ is for each $\theta \in \Theta$
and $i = 0,1,2$ Borel measurable; $E[\sup_{\theta\in\Theta
}|\sigma_{t}^{2}/\sigma_{t}^{2}(\theta)|^{p}] < \infty$ for
any $p > 0$; $E[\sup_{\theta\in\mathcal{N}_{0}}\|(\partial
/\partial\theta)^{i}\ln(\sigma_{t}^{2}(\theta))\|^{2+\iota}] <
\infty$ for $i = 1,2$, tiny $\iota > 0$, and some compact $\mathcal
{N}_{0} \subseteq \Theta$
containing $\theta^{0}$ and having positive Lebesgue measure.

 \item[(d)] $\{y_{t}\}$ and $\{\sigma_{t}^{2}(\theta)\}$ for $%
\theta \in \Theta$ are stationary and geometrically $\beta$-mixing.
\end{enumerate}
\end{Assumption}

%re1 #&#
%
\begin{remark}
The tail index $\kappa$ in (b) is identically the moment supremum %
(see Resnick \cite{Resnick87}). The volatility moment bounds in (c)
imply only
the tails of $\epsilon_{t}$ matter for Gaussian asymptotics, and can be
relaxed at the expense of added notation for trimming also according to
$%
\mathfrak{s}_{t}$. Verification of (c) for the classic GARCH model is in
Francq and Zako\"{\i}an \cite{FZ04}, and related proofs for asymmetric
models are in Francq and Zako\"{\i}an \cite{FZ10}.
\end{remark}

%re2 #&#
%
\begin{remark}
Geometric $\beta$-mixing (d) implies mixing in the ergodic sense,
hence ergodicity (see Petersen \cite{Petersen83}). Lipschitz type
conditions on the
volatility response $g$ combined with a smooth bounded distribution for
$%
\epsilon_{t}$ suffice, covering a large variety of models %
(Carrasco and Chen~\cite{carrascochen02}, Straumann and Mikosch \cite
{StraumannMikosch},
Meitz and Saikkonen \cite{meitzsaikkonen08}, Meitz and Saikkonen
\cite{meitzsaik11}).
See Theorem~\ref{th:feas} below for one such set of conditions. In the
classic GARCH model $y_{t} = \sigma_{t}\epsilon_{t}$ and $\sigma
_{t}^{2}(\theta) = \omega + \alpha y_{t-1}^{2} + \beta\sigma
_{t-1}^{2}(\theta)$, for example, where $\omega > 0$, $\alpha,\beta \geq
0$\vadjust{\goodbreak} and $E[\ln(\alpha^{0}\epsilon_{t}^{2} + \beta^{0})] < 0$
ensure stationarity and ergodicity, and combined with $E[\epsilon
_{t}^{2}] = 1$ this allows for IGARCH and mildly explosive cases
$\alpha^{0} +
\beta^{0} \geq 1$ (Nelson \cite{Nelson90}). If additionally $\epsilon_{t}$
has a continuous distribution that is positive on $(-\infty,\infty)$
then $%
\{y_{t},\sigma_{t}^{2}(\theta)\}$ are geometrically $\beta$-mixing %
(Carrasco and Chen \cite{carrascochen02}).
\end{remark}

In the \hyperref[app]{Appendices}, we show $\hat{\theta}_{n}$ obtains the expansion
$\mathcal{%
V}_{n}^{1/2}(\hat{\theta}_{n} - \theta^{0}) = n^{-1/2}\Sigma
_{n}^{-1/2}\*\sum_{t=1}^{n}m_{t}I_{n,t}^{(\mathcal{E})}(1+\mathrm{o}_{p}(1))$,
hence $n^{1/2}\Sigma_{n}^{-1/2}E[m_{t}I_{n,t}^{(\mathcal{E})}]
\rightarrow0$ must hold for asymptotic unbiasedness of $\hat{\theta
}_{n}$%
. This reduces to assuming $n^{1/2}(E[\mathcal
{E}_{t}^{2}I_{n,t}^{(\mathcal{E%
})}])^{-1/2}E[\mathcal{E}_{t}I_{n,t}^{(\mathcal{E})}]\rightarrow0$
since by independence $E[m_{t}I_{n,t}^{(\mathcal{E})}]=E[\mathcal{E}%
_{t}I_{n,t}^{(\mathcal{E})}]\times E[\mathfrak{s}_{t}]$, while $\Sigma
_{n}=E[\mathcal{E}_{t}^{2}I_{n,t}^{(\mathcal{E})}]\times E[%
\mathfrak{s}_{t}\mathfrak{s}_{t}^{\prime}]$ and $\|E[\mathfrak{s}_{t}%
\mathfrak{s}_{t}^{\prime}]\| \in (0,\infty)$.

\begin{Assumption}[(Identification)]\label{ass2}
The fractile sequences $%
\{k_{1,n},k_{2,n}\}$ satisfy
$n^{1/2}\times\break (E[\mathcal{E}%
_{t}^{2}I_{n,t}^{(\mathcal{E})}])^{-1/2} E[\mathcal
{E}_{t}I_{n,t}^{(\mathcal{E%
})}]\rightarrow0$ where $\mathcal{E}_{t}:=\epsilon
_{t}^{2}-1$.
\end{Assumption}

%re3 #&#
%
\begin{remark}
We do not require $E[\mathcal{E}_{t}I_{n,t}^{(\mathcal{E})}] = 0$
for finite $n$ since our results are asymptotic, while $E[\mathcal{E}%
_{t}I_{n,t}^{(\mathcal{E})}] \rightarrow E[\epsilon_{t}^{2}-1] = 0$
automatically holds by dominated convergence and negligibility $k_{i,n}/n
= \mathrm{o}(1)$. Since $n^{1/2}/(E[\mathcal{E}_{t}^{2}I_{n,t}^{(\mathcal{E}%
)}])^{1/2} \rightarrow\infty$ as verified in Section~\ref{sec:rate}
below, we require $E[\mathcal{E}_{t}I_{n,t}^{(\mathcal{E})}]
\rightarrow 0$ fast enough, else there is asymptotic bias.
\end{remark}

%re4 #&#
%
\begin{remark}
There always exists a sequence $\{k_{1,n},k_{2,n}\}$ such that $E[%
\mathcal{E}_{t}I_{n,t}^{(\mathcal{E})}]$ is closer to zero than
$(E[\mathcal{%
E}_{t}^{2}I_{n,t}^{(\mathcal{E})}])^{1/2}/n^{1/2}$ as $n$ increases. In
general $\mathcal{E}_{t} \in [-1,\infty)$ is skewed right hence,
counterintuitively, asymptotic unbiasedness requires $k_{1,n} > k_{2,n}$: a few trimmed large positive values promotes asymptotic normality, but
forces us to trim many negative values to ensure identification. See
Section~\ref{ver_2} for discussion and examples. In a method of moments framework,
however, identification is assured by re-centering the trimmed errors, hence
Assumption~\ref{ass2} is not required. See Section~\ref{sec:mm}.
\end{remark}

%re5 #&#
%
\begin{remark}
Define $m_{n,t} := m_{t}I_{n,t}^{(\mathcal{E})}$. Assumption~\ref{ass2}
ensures $E[\{m_{n,s} - E[m_{n,s}]\}\{m_{n,t} - E[m_{n,t}]\}^{\prime
}] = E[m_{n,s}m_{n,t}^{\prime}] + \mathrm{o}(\|\Sigma_{n}\|/n)$ for all $%
s,t $, and $\|\sum_{i=1}^{n-1}E[m_{n,1}m_{n,i+1}^{\prime}]\| \leq n
 \mathrm{o}\times\break (\|\Sigma_{n}\|/n) = \mathrm{o}(\|\Sigma_{n}\|)$ by Minkowski and
Cauchy--Schwarz inequalities. Hence, $\Sigma_{n}$ is asymptotically
equal to
the long-run covariance matrix $\mathcal{S}_{n}$ of $n^{-1/2}\sum_{t=1}^{n}%
\{m_{n,t} - E[m_{n,t}]\}$ since
%e6 #&#
\begin{eqnarray}
 &&E \Biggl[ \Biggl( \frac{1}{n^{1/2}}\sum_{t=1}^{n}
\bigl\{m_{n,t}- E[m_{n,t}]\bigr\} \Biggr) \Biggl(
\frac{1}{n^{1/2}}\sum_{t=1}^{n}\bigl
\{m_{n,t}-%
E[m_{n,t}]\bigr\} \Biggr)
^{\prime} \Biggr] \label{sigma}
\nonumber
\\[-8pt]
\\[-8pt]
&&\quad =\Sigma_{n}\times \bigl( 1+\mathrm{o} ( 1 ) \bigr) +2\sum
_{i=1}^{n-1} \biggl( 1-\frac{i}{n} \biggr) E
\bigl[ m_{n,1}m_{n,i+1}^{\prime} \bigr] =
\Sigma_{n}\times \bigl( 1+\mathrm{o} ( 1 ) \bigr) .
\nonumber
\end{eqnarray}
\end{remark}

We are now ready to state the main results of this section. The
expansion $%
\mathcal{V}_{n}^{1/2}(\hat{\theta}_{n} - \theta^{0}) =
n^{-1/2}\Sigma_{n}^{-1/2}\sum_{t=1}^{n}\hspace*{-0.5pt}m_{t}I_{n,t}^{(\mathcal{E})}(1 +
\mathrm{o}_{p}(1))$ requires Jacobian consistency 
$1/n\sum_{t=1}^{n}\hspace*{-0.5pt}G_{t}(\hat
{\theta}%
_{n})\times\break \hat{I}_{n,t}^{(\mathcal{E})}(\hat{\theta}_{n}) \stackrel{p}{%
\rightarrow} -E [ \mathfrak{s}_{t}\mathfrak{s}_{t}^{\prime}
] $
and therefore consistency $\hat{\theta}_{n} \stackrel{p}{%
\rightarrow} \theta^{0}$ from first principles. Proofs of main
results are contained in Appendices \ref{app:th21_22} and \ref{app:remain}.

%th2.1 #&#
%
\begin{theorem}[(QMTTL consistency)]
\label{th:consist} Under Assumptions \ref{ass1} and \ref{ass2} $\hat{\theta}_{n} \stackrel
{p}{\rightarrow} \theta^{0}$.
\end{theorem}

%th2.2 #&#
%
\begin{theorem}[(QMTTL normality)]
\label{th:norm} Under Assumptions \ref{ass1} and \ref{ass2} $\mathcal{V}_{n}^{1/2}(\hat{%
\theta}_{n} - \theta^{0})\stackrel{d}{\rightarrow}N(0,I_{q})$
where $\mathcal{V}_{n} = nE[\mathfrak{s}_{t}\mathfrak{s}_{t}^{\prime}]%
\mathcal{S}_{n}^{-1}E[\mathfrak{s}_{t}\mathfrak{s}_{t}^{\prime}] \sim
n(E[\mathcal{E}_{t}^{2}I_{n,t}^{(\mathcal{E})}])^{-1}E[\mathfrak{s}_{t}%
\mathfrak{s}_{t}^{\prime}]$ and each $\mathcal{V}_{i,i,n}
\rightarrow \infty$.
\end{theorem}

Now consider feasible QMTTL. Define $\tilde{\epsilon}_{t}(\theta) :=
y_{t}^{2}/h_{t}(\theta)$ based on the iterated process $\{h_{t}(\theta
)\}$
in (\ref{h_iter}), and $\mathcal{\tilde{E}}_{t}(\theta) := \tilde{%
\epsilon}_{t}^{2}(\theta) - 1$. The feasible estimator is
\[
\tilde{\theta}_{n}=\operatorname{arg\,min}\limits
_{\theta\in\Theta
} \Biggl\{
\frac{1}{%
n}\sum_{t=1}^{n} \bigl( \ln
h_{t}(\theta)+\tilde{\epsilon }_{t}^{2}(\theta )
\bigr) \times I \bigl( \mathcal{\tilde{E}}_{(k_{1,n})}^{(-)}(\theta )
\leq \mathcal{\tilde{E}}_{t}(\theta)\leq\mathcal{\tilde
{E}}_{(k_{2,n})}^{(+)}(%
\theta) \bigr) \Biggr\} .
\]
Under the following Lipschitz bounds for the response $g$ and its
derivatives we show $\tilde{\theta}_{n}$ has the same limit
distribution as
the infeasible $\hat{\theta}_{n}$, cf. Meitz and Saikkonen \cite
{meitzsaikkonen08}. Related
ideas are contained in Straumann and Mikosch \cite{StraumannMikosch}.

Drop arguments: $g = g(y,s,\theta)$, and let $g_{a}$ and $g_{a,b}$
denote first and second derivatives for $a,b \in \{y,s,\theta\}$. We
say a matrix function $\xi(y,s,\theta)$ is Lipschitz in $s$ if $\|\xi
(y,s_{1},\theta) - \xi(y,s_{2},\theta)\| \leq K|s_{1} -
s_{2}|$ $\forall s_{1},s_{2} \in [0,\infty)$ and $y,\theta \in
\mathbb{R} \times \Theta$.

\begin{Assumption}[(Response bounds)]\label{ass3}
\begin{enumerate}[(b)]
 \item[(a)]  $g \leq \rho s + K(1 + y^{2})$ for some $\rho\in
(0,1)$ and $\inf_{y\in\mathbb{R},s\in\mathbb{R}_{+},\theta\in
\Theta}\{|g|\} > 0$;

 \item[(b)] $\|g_{a}\|$ and $\|g_{a,b}\|$ are bounded by $K(1 + y^{2} + s)$ for
each $a,b \in \{y,\theta\}$;

 \item[(c)] $g$, $g_{a}$ and $g_{a,b}$ are Lipschitz in $s$, for each $a,b
\in \{y,s,\theta\}$.
\end{enumerate}
\end{Assumption}

Assumption~\ref{ass3} ensures $h_{t}(\theta), h_{t}^{\theta}(\theta) :=
(\partial/\partial\theta)h_{t}(\theta)$ and $h_{t}^{\theta,\theta
}(\theta) := (\partial/\partial\theta)h_{t}^{\theta}(\theta)$ have
stationary ergodic solutions $\{h_{t}^{\ast}(\theta),h_{t}^{\theta
\ast
}(\theta),h_{t}^{\theta,\theta\ast}(\theta)\}$ with the geometric
property\vspace*{2pt}\break  $E[(\sup_{\theta\in\Theta}|a_{t}^{\ast}(\theta) -
a_{t}(\theta)|)^{\iota}] = \mathrm{o}(\rho
^{t}) $ for each\vspace*{2pt} $a_{t}(\theta) \in \{h_{t}(\theta),h_{i,t}^{\theta
}(\theta),h_{i,j,t}^{\theta,\theta}(\theta)\} \mbox{ and }a_{t}^
{\ast}(\theta)
\in\{h_{t}^{\ast}(\theta),h_{i,t}^{\theta\ast}(\theta
), h_{i,j,t}^{\theta,\theta\ast}(\theta)\}$ and some $\rho
\in(0,1)$. See Lemma~\ref{lm:sol} in Appendix~\ref{app:remain}%
. This leads to the next result.

%th2.3 #&#
%
\begin{theorem}[(Feasible QMTTL)]\label{th:feas} Under Assumptions \ref{ass1}--\ref{ass3} $\mathcal{V}_{n}^{1/2}(\tilde{%
\theta}_{n}-\hat{\theta}_{n}) \stackrel{p}{\rightarrow} 0$.
\end{theorem}

%re6 #&#
%
\begin{remark}
In the remainder of the paper, we focus on the infeasible $\hat{\theta}
_{n}$ for notational economy.
\end{remark}

As stated above, we need only trim by error extremes since first order
asymptotics rests solely on whether $\epsilon_{t}$ has a fourth moment or
not. However, in small samples a large $y_{t-1}$ may cause $\mathfrak{s}_{t}$
or $\mathfrak{d}_{t}$ to spike and therefore the score equation to
exhibit a
sample extreme value. Consider, for example, that in the linear volatility
model $\sigma_{t}^{2}(\theta) = \omega + \alpha y_{t-1}^{2} +
\beta\sigma_{t-1}^{2}(\theta)$ the score weight at the origin
$\mathfrak{s%
}_{t}(\theta)|_{\alpha,\beta=0} = \omega^{-1} \times
[1,y_{t-1}^{2},\omega]^{\prime}$ obtains an extreme value \textit{if and
only if} $|y_{t-1}|$ does. In general $\mathfrak{s}_{t}$ exhibits spikes
when $|y_{t-1}|$ does for $\alpha^{0}$ and $\beta^{0}$ near zero. This
same properly applies to a large variety of GARCH models. Thus,
although $%
\hat{\theta}_{n}$ is consistent and asymptotically normal, for improved
small sample performance trimming by large values of $y_{t-1}$ appears
to be
highly useful in practice. This is not surprising since true additive
outliers render QML biased (see Mendes \cite{Mendes00}, Muler and Yohai
\cite{MulerYohai08}, cf.
Cavaliere and Georgiev \cite{CavaliereGeorgiev09}, Muler \textit
{et~al}. \cite{Mulleretal09}).\vadjust{\goodbreak}

Let $\{\tilde{k}_{n}\}$ be an intermediate order sequence and define
$\hat{I}%
_{n,t}^{(y)} := I(|y_{t}| \leq y_{(\tilde{k}_{n})}^{(a)})$ where $%
y_{(i)}^{(a)}$ are order statistics of $y_{t}^{(a)} := |y_{t}|$. The
estimator in this case is
\[
\hat{\theta}_{n}^{(y)}=\operatorname{arg\,min}\limits
_{\theta\in
\Theta}
\Biggl\{ \frac{1}{n}\sum_{t=1}^{n}
\bigl( \ln\sigma_{t}^{2}(\theta)+\epsilon
_{t}^{2}(\theta) \bigr) \times\hat{I}_{n,t}^{(\mathcal{E})}(
\theta )\hat{I}%
_{n,t-1}^{(y)} \Biggr\} .
\]
Since $\hat{I}_{n,t-1}^{(y)} \stackrel{p}{\rightarrow} 1$, the score
equations $\mathfrak{s}_{t}$ are square integrable, and $\epsilon_{t}$ is
i.i.d., asymptotic normality does not depend on whether $y_{t}$ is heavy
tailed. Indeed, it is easy to show $\hat{\theta}_{n}^{(y)}$ is
asymptotically equivalent to $\hat{\theta}_{n}$. The same property extends
to feasible QMTTL with trimming by $y_{t-1}$, denoted~$\tilde{\theta}%
_{n}^{(y)}$. We therefore omit the proof of the next result.

%co2.4 #&#
%
\begin{corollary}
Under Assumptions \ref{ass1} and \ref{ass2}, trimming by $y_{t-1}$ does not
impact the limit distributions of infeasible and feasible QMTTL
estimators:
$\mathcal{V}_{n}^{1/2}(\hat{\theta}_{n}^{(y)} - \hat{\theta
}_{n})\stackrel%
{p}{\rightarrow}0$ and $\mathcal{V}_{n}^{1/2}(\tilde{\theta}_{n}^{(y)} -
\tilde{\theta}_{n})\stackrel{p}{\rightarrow}0$. Moreover, infeasible and
feasible estimators are asymptotically equivalent: $\mathcal
{V}_{n}^{1/2}(%
\tilde{\theta}_{n}^{(y)}-\hat{\theta}_{n}^{(y)})\stackrel
{p}{\rightarrow}0$.
\end{corollary}

%s2.3 #&#
\subsection{Verification of identification Assumption \texorpdfstring{\protect\ref{ass2}}{2}}\label{ver_2}

We require an explicit model of $P(|\epsilon_{t}| > c)$ in order to
verify Assumption~\ref{ass2}. In our simulation study, we use distributions with
either power law or exponential tail decay.

%s2.3.1 #&#
\subsubsection{Paretian tails}

In the simulation experiment we use
%e7 #&#
\begin{equation}
P \bigl( \llvert \epsilon_{t}\rrvert >c \bigr) = ( 1+c )
^{-\kappa}\quad\quad\mbox{with }\kappa\in(2,4), \label{pow_e_twotail}
\end{equation}
hence $\mathcal{E}_{t}$ has left and right tails:%
%e8 #&#
\begin{eqnarray}\label{pow_E_asym}
 P ( \mathcal{E}_{t}<-c ) &=&P \bigl( \epsilon_{t}^{2}<1-c
\bigr) =0%
\quad\quad\mbox{if }c\geq1 \nonumber,
\\
  &=&1-P \bigl( \epsilon _{t}^{2}>1-c \bigr) =1- ( 2-c )
^{-\kappa}\quad\quad\mbox{if }c\in [ 0,1 ],
\\
 P ( \mathcal{E}_{t}>c ) &=&P \bigl( \epsilon_{t}^{2}>1+c
\bigr) = ( 2+c ) ^{-\kappa/2}.
\nonumber
\end{eqnarray}
We show below identification $n^{1/2}(E[\mathcal{E}_{t}^{2}I_{n,t}^{(%
\mathcal{E})}])^{-1/2}E[\mathcal{E}_{t}I_{n,t}^{(\mathcal{E})}]
\rightarrow 0$ holds if $k_{1,n} \rightarrow \infty$, $k_{1,n}/n
\rightarrow0$ and:
%e9 #&#
\begin{eqnarray}\label{k1k2}
\biggl( \frac{k_{2,n}}{n} \biggr) ^{1-2/\kappa}&=&\frac{\kappa-2}{2} \biggl(
-1+ \biggl( \frac{1}{1-k_{1,n}/n} \biggr) ^{2/\kappa}+\frac{2}{\kappa
-2}
\frac{%
k_{1,n}}{n} \biggr)\nonumber
\\[-8pt]\\[-8pt]
&&{} +\mathrm{o} \biggl( \biggl( \frac{n}{k_{1,n}} \biggr)
^{2/\kappa
-1/2}%
\frac{1}{n^{1/2}} \biggr) . \nonumber
\end{eqnarray}
In practice, (\ref{k1k2}) is greatly simplified asymptotically by
noting $%
(n/k_{1,n})^{2/\kappa-1/2}n^{-1/2} = \mathrm{o}(1)$ and $(1 -
k_{1,n}/n)^{-2/\kappa} - 1 \sim (2/\kappa)(k_{1,n}/n)$, hence
identification applies if $(k_{2,n}/n)^{1-2/\kappa} \sim ((2\kappa
- 2)/\kappa)(k_{1,n}/n)$ or%
%e10 #&#
\begin{equation}
\frac{k_{2,n}}{k_{1,n}^{\kappa/(\kappa-2)}}\sim2 \biggl( 1-\frac
{1}{\kappa}%
\biggr)
^{\kappa/(\kappa-2)}\frac{1}{n^{\kappa/(\kappa-2)-1}}. \label{kn2}
\end{equation}
A similar condition applies in the second order power law case
$P(|\epsilon
_{t}| > c) = dc^{-\kappa}(1+ec^{-\xi})$ with $d,e > 0$, $\xi > 0$
and $\kappa \in (2,4)$, while a less sharp result arises under
$P(|\epsilon_{t}| > c) = dc^{-\kappa}(1+\mathrm{o}(1))$.

In order to show (\ref{k1k2}), we must characterize the moments
$E[\mathcal{E}%
_{t}I_{n,t}^{(\mathcal{E})}] = E[\mathcal{E}_{t}I(-\mathcal{L}_{n}
\leq \mathcal{E}_{t} \leq \mathcal{U}_{n})]$ and $E[\mathcal{E}%
_{t}^{2}I_{n,t}^{(\mathcal{E})}]$. Use (\ref{pow_E_asym}) to deduce $%
\mathcal{U}_{n} = (n/k_{2,n})^{2/\kappa} - 2 \rightarrow
\infty$ and $\mathcal{L}_{n} = 2 - (n/(n-k_{1,n}))^{2/\kappa}
\in [0,1]$ as $n \rightarrow \infty$. Therefore,
%e11 #&#
\begin{eqnarray}\label{EE}
&&E \bigl[ \mathcal{E}_{t}I ( -\mathcal{L}_{n}\leq
\mathcal{E}_{t}\leq \mathcal{U}_{n} ) \bigr]\nonumber
\\
 &&\quad=- \bigl\{ E
\bigl[ \mathcal{E}_{t}I ( \mathcal{E}_{t}>
\mathcal{U}_{n} ) \bigr] +E \bigl[ \mathcal{E}%
_{t}I ( \mathcal{E}_{t}<-\mathcal{L}_{n} ) \bigr]
\bigr\} \nonumber
\\[-8pt]\\[-8pt]
&&\quad=- \biggl\{ \int_{\mathcal{U}_{n}}^{\infty} ( 2+u )
^{-\kappa
/2}\,\mathrm{d}u-\int_{\mathcal{L}_{n}}^{1} \bigl( 1- ( 2-u )
^{-\kappa
/2} \bigr) \,\mathrm{d}u \biggr\}
\nonumber
\\
&&\quad=- \biggl\{ \frac{2}{\kappa-2} \biggl( \frac{k_{2,n}}{n} \biggr)
^{1-2/\kappa
}+1- \biggl( \frac{n}{n-k_{1,n}} \biggr) ^{2/\kappa}-
\frac{2}{\kappa
-2}\frac{%
k_{1,n}}{n} \biggr\} .
\nonumber
\end{eqnarray}
Next $E[\mathcal{E}_{t}^{2}I_{n,t}^{(\mathcal{E})}] \sim K (
n/k_{2,n} ) ^{4/\kappa-1}$ follows from (\ref{mom_e4}) below. Combined
with (\ref{EE}) and by rearranging terms, Assumption~\ref{ass2} holds when $E[%
\mathcal{E}_{t}I_{n,t}^{(\mathcal{E})}] = \mathrm{o}((E[\mathcal{E}%
_{t}^{2}I_{n,t}^{(\mathcal{E})}])^{1/2}/n^{1/2})$, hence when%
%e12 #&#
\begin{eqnarray}\label{k1k2*}
\biggl( \frac{k_{2,n}}{n} \biggr) ^{1-2/\kappa}&=&\frac{\kappa-2}{2} \biggl(
-1+ \biggl( \frac{1}{1-k_{1,n}/n} \biggr) ^{2/\kappa}+\frac{2}{\kappa
-2}
\frac{%
k_{1,n}}{n} \biggr) \nonumber
\\[-8pt]\\[-8pt]
&&{}+\mathrm{o} \biggl( \biggl( \frac{n}{k_{2,n}} \biggr)
^{2/\kappa
-1/2}%
\frac{1}{n^{1/2}} \biggr) . \nonumber
\end{eqnarray}
Notice $k_{2,n}$ appears on both sides of the equality. In order to
achieve (%
\ref{k1k2}), note $k_{2,n}/k_{1,n} \rightarrow 0$. This follows since $
(n/k_{2,n})^{2/\kappa-1/2}n^{-1/2} = \mathrm{o}(1)$ and by the
mean-value-theorem $(1 - k_{1,n}/n)^{-2/\kappa} - 1 \sim
(2/\kappa)k_{1,n}/n$ hence $(k_{2,n}/n)^{1-2/\kappa} \sim Kk_{1,n}/n$, therefore
\[
\biggl( \frac{k_{2,n}}{k_{1,n}} \biggr) ^{1-2/\kappa}= \biggl( \frac
{k_{2,n}/n}{%
k_{1,n}/n}
\biggr) ^{1-2/\kappa}\sim K\frac{ ( k_{1,n}/n )
}{ (
k_{1,n}/n ) ^{1-2/\kappa}}=K ( k_{1,n}/n )
^{2/\kappa
}\rightarrow0.
\]
Now combine $k_{2,n}/k_{1,n} \rightarrow 0$ and (\ref{k1k2*}) to deduce
(\ref{k1k2}).

There are several things to note from (\ref{kn2}). First, there are
arbitrarily many valid $\{k_{1,n},k_{2,n}\}$. Second, $\{
k_{1,n},k_{2,n}\}$
requires knowledge of $\kappa$, which can be consistently estimated for
many processes defined by (\ref{garch}) (see Hill \cite{Hill10}).
However, the
method of moments estimator in Section~\ref{sec:mm} only requires one
two-tailed fractile \textit{without} knowledge of $\kappa$.

Third, $k_{2,n}/k_{1,n} \rightarrow 0$ since $k_{1,n}/n \rightarrow
0$ and $\kappa > 2$. This logically follows since $\mathcal{E}_{t}$
has support $[-1,\infty)$. The right tail is heavier, hence trimming a
positive extreme must be off-set by trimming more negative observations in
order to get $E[\mathcal{E}_{t}I_{n,t}^{(\mathcal{E})}] \approx 0$.

Fourth, $k_{1,n} \sim n/g_{1,n}$ for slowly varying $g_{2,n}
\rightarrow \infty$ implies $k_{2,n} \sim n/g_{2,n}$ for slowly
varying $g_{2,n} , g_{2,n}/g_{1,n} \rightarrow \infty$. Similarly, $%
k_{1,n} \sim \lambda_{1}n^{\delta_{1}}$ for $\lambda_{1} \in
(0,1)$ and $\delta_{1} \in (2/\kappa,1)$ implies $k_{2,n} \sim
\lambda_{2}n^{\delta_{2}}$ for $\lambda_{2} \in (0,1)$ and $\delta
_{2} \in (0,\delta_{1})$. Further, slowly varying $k_{1,n}
\rightarrow \infty$ is not valid since \mbox{$k_{2,n} \rightarrow 0$} is
then required which leads to asymptotic non-normality when $E[\epsilon
_{t}^{4}] = \infty$.

Fifth, we need monotonically larger $k_{1,n}$ as $\kappa \searrow2$,
but always $\lim\sup_{n\rightarrow\infty}(k_{2,n}/k_{1,n}) < 1$.
Exponential tails treated in Section~\ref{sec:expon} reveals an extreme
case: there are no limitations on how we set $\{k_{1,n},k_{2,n}\}$ outside
of an upper bound, although $k_{1,n} > k_{2,n}$ always reduces small
sample bias.

Finally, as a numerical example suppose $\kappa = 2.5$ and $n = 100$.
If $k_{2,n} = 1$ then $k_{1,n} = 33$ renders (\ref{kn2}) a near
equality, although any $k_{1,n} \in \{29,\ldots,35\}$ aligns with
$k_{2,n} = 1$ by rounding. This is striking: we need to trim roughly 33
times as
many negative $\mathcal{E}_{t}(\theta)$ as positive $\mathcal
{E}_{t}(\theta
)$ to approach unbiasedness at $n = 100$. If $n = 800$ then, for
example, $k_{2,n} = 2$ aligns with roughly $k_{1,n} = 200$.

%s2.3.2 #&#
\subsubsection{Exponential tails}\label{sec:expon}

Now suppose $\epsilon_{t}$ has a Laplace distribution:
\[
P ( \epsilon_{t}\leq-c ) =\tfrac{1}{2}\exp \{ -\sqrt{2}c \} \quad\quad\mbox{for
}c>0\quad\mbox{and}\quad P ( \epsilon_{t}>c ) =\tfrac
{1}{2}\exp \{ -\sqrt{2}c \}\quad\quad
\mbox{for }c\geq0.
\]
We use a normal distribution in our simulation study, but the exposition
here is greatly simplified under Laplace, while the conclusions are the same.

We have $P(\mathcal{E}_{t} \leq -c) = 1 - \exp\{-\sqrt{2}(1 -
c)^{1/2}\}$ and $P(\mathcal{E}_{t} \geq c) = \exp\{-\sqrt{2}(1 +
c)^{1/2}\}$. The following are then straightforward to verify: $\mathcal
{L}%
_{n} = 1 - (\ln(n/(n-k_{1,n})))^{2}$ and $\mathcal{U}_{n} =
(\ln(0.5n/k_{2,n}))^{2} - 1$, hence
\begin{eqnarray*}
&&E \bigl[ \mathcal{E}_{t}I ( -\mathcal{L}_{n}\leq
\mathcal{E}_{t}\leq \mathcal{U}_{n} ) \bigr]
\\
&&\quad=2 \biggl(
\frac{k_{1,n}}{n}-\frac
{k_{2,n}}{n}%
\biggr) -\ln \biggl(
\frac{n}{n-k_{1,n}} \biggr) \biggl\{ -\ln \biggl( \frac
{n}{%
n-k_{1,n}} \biggr) +2
\biggl( \frac{n-k_{1,n}}{n} \biggr) \biggr\} .
\end{eqnarray*}
Observe $E[\mathcal{E}_{t}I(-\mathcal{L}_{n} \leq \mathcal{E}_{t}
\leq \mathcal{U}_{n})] \approx 0$ when $k_{1,n}>k_{2,n}$, hence if $%
k_{1,n}/n \rightarrow 0$ then $k_{2,n}/n \rightarrow 0$ must hold.

Since $E[\epsilon_{t}^{4}] < \infty$ we need $E[\mathcal{E}_{t}I(-%
\mathcal{L}_{n} \leq \mathcal{E}_{t} \leq \mathcal{U}_{n})] =
\mathrm{o}(1/n^{1/2})$. Notice $\ln(n/(n - k_{1,n})) \sim k_{1,n}/n$. Hence
if simply each $k_{i,n} = \mathrm{o}(n^{1/2})$, then we achieve $E[\mathcal{E}%
_{t}I(-\mathcal{L}_{n} \leq \mathcal{E}_{t} \leq \mathcal{U}_{n})] =
\mathrm{o}(1/n^{1/2})$. This implies that technically we do not even need
asymmetric trimming $k_{1,n} > k_{2,n}$ as long as we set $k_{1,n} =
k_{2,n} = \mathrm{o}(n^{1/2})$. This follows since tails are so thin that in
general extremes on $[0,\infty)$ are not much larger than extremes on $
[-1,0)$ in small samples. Similarly, we can use any form of asymmetric
trimming that satisfies $k_{i,n} = \mathrm{o}(n^{1/2})$. We show by simulation
that as $n$ gets large, bias evaporates irrespective of $k_{i,n}$, but $
k_{1,n} > k_{2,n}$ always leads to lower small sample bias.

%s2.3.3 #&#
\subsubsection{Remarks}\label{sec:bias_remarks}

We demonstrate by simulation in Section~\ref{sim} that using $k_{1,n} =
10k_{2,n}$ or $k_{1,n} = 35k_{2,n}$ for either $n \in \{100,800\}$
and either Paretian or Gaussian $\epsilon_{t}$ leads to a superb QMTTL
estimator. Indeed, simply using symmetric trimming $k_{1,n} = k_{2,n}$
still leads to a better estimator than Log-LAD and Weighted Laplace QML in
terms of small sample bias and approximate normality, although Power-Law
QML tends to have lower bias and be closer to normal. In general using bias
minimizing fractiles, like $k_{1,n} = 100k_{2,n}$ for Paretian $%
\epsilon_{t}$ when $n = 800$, is not evidently required for obtaining
low bias in finite samples, as long as $k_{1,n}$ is comparatively large
relative to $k_{2,n}$ in which case QMTTL trumps Log-LAD, WLQML and PQML.

We also find that our method of moments estimator in Section~\ref{sec:mm}
dominates Log-LAD, WLQML and PQML, although QMTTL with $k_{1,n} =
35k_{2,n}$ leads to smaller bias and is closer to normally distributed in
nearly every case. Nevertheless, the method of moments estimator is always
asymptotically unbiased and easier to implement because trimming is
symmetric. Which estimator is chosen in practice depends on the analyst's
preferences: method of moments is guaranteed to be asymptotically unbiased,
but QMTTL has superior small sample properties even if $\{
k_{1,n},k_{2,n}\}$
are not chosen to ensure asymptotic unbiasedness.

%s2.4 #&#
\subsection{QMTTL scale and rate of convergence}\label{sec:rate}

The scale $\mathcal{V}_{n}$ and rate of convergence depend on the error tail
index $\kappa > 2$. If $E[\epsilon_{t}^{4}] < \infty$ then by
dominated convergence $E[\mathcal{E}_{t}^{2}I_{n,t}^{(\mathcal{E})}] =
E[(\epsilon_{t}^{2} - 1)^{2}I_{n,t}^{(\mathcal{E})}] \rightarrow %
E[(\epsilon_{t}^{2} - 1)^{2}] = E[\epsilon_{t}^{4}] - 1$, thus
$\mathcal{V}_{n} \sim n(E[\epsilon_{t}^{4}] - 1])^{-1}E[\mathfrak{s%
}_{t}\mathfrak{s}_{t}^{\prime}]$, the classic QML asymptotic covariance
matrix. This implies trimming does not affect efficiency asymptotically.
Hence, we now assume $E[\epsilon_{t}^{4}] = \infty$.

Let the intermediate order sequences $\{k_{n}\}$ and positive
thresholds $\{%
\mathcal{C}_{n}(\theta)\}$ satisfy
\[
P \bigl( \bigl\llvert \mathcal{E}_{t}(\theta)\bigr\rrvert \geq
\mathcal{C}%
_{n}(\theta) \bigr) =\frac{k_{n}}{n}.
\]
The rate $E[\mathcal{E}_{t}^{2}(\theta)I_{n,t}^{(\mathcal{E})}(\theta
)] %
\rightarrow \infty$ is logically governed by the right tail of
$\mathcal{%
E}_{t}(\theta) = \epsilon_{t}^{2}(\theta) - 1 \in [-1,\infty
)$ since by dominated convergence:
\begin{eqnarray*}
E \bigl[ \mathcal{E}_{t}^{2}(\theta)I_{n,t}^{(\mathcal{E})}(
\theta ) \bigr] &=&E \bigl[ \mathcal{E}_{t}^{2}(\theta)I \bigl(
-\mathcal{L}_{n}(\theta )\leq \mathcal{E}_{t}(\theta)\leq
\mathcal{U}_{n}(\theta) \bigr) \bigr]
\\
 &\sim& E%
\bigl[
\mathcal{E}_{t}^{2}(\theta)I \bigl( \bigl\llvert \mathcal
{E}_{t}(\theta )\bigr\rrvert \leq\mathcal{C}_{n}(\theta)
\bigr) \bigr]
\end{eqnarray*}
as though $\mathcal{E}_{t}(\theta)$ were symmetrically trimmed with
thresholds and fractile%
%e13 #&#
\begin{equation}\label{CU}
\mathcal{C}_{n}(\theta)=\mathcal{U}_{n}(\theta)\quad\mbox{and}\quad%
k_{n}=k_{2,n}.
\end{equation}
Note $E[\mathcal{E}_{t}^{2}(\theta)I_{n,t}^{(\mathcal{E})}(\theta)]
\sim
E[\mathcal{E}_{t}^{2}(\theta)I(|\mathcal{E}_{t}(\theta)| \leq %
\mathcal{C}_{n}(\theta))]$ is useful for characterizing the convergence
rate, but identification Assumption~\ref{ass2} in general requires $k_{1,n} > %
k_{2,n}$ hence $\mathcal{L}_{n}(\theta) < \mathcal{U}_{n}(\theta)$.

As long as $E[\epsilon_{t}^{4}] = \infty$, then the rate of convergence
is $\mathcal{V}_{n}^{1/2} = \mathrm{o}(n^{1/2})$: heavy tailed errors can only
adversely affect the convergence rate. The exact rate can be deduced by
observing that from $P(|\epsilon_{t}| > a) = da^{-\kappa}(1 + %
\mathrm{o}(1))$ the variable $\mathcal{E}_{t} = \epsilon_{t}^{2} - 1$ has a
tail sum dominated by the right tail:%
%e14 #&#
\begin{eqnarray}\label{pow_E}
P \bigl( \llvert \mathcal{E}_{t}\rrvert >a \bigr) &=&P \bigl(
\epsilon _{t}^{2}>1+a \bigr) +P \bigl( \epsilon_{t}^{2}<1-a
\bigr)\nonumber
\\[-8pt]\\[-8pt]
&=&d ( 1+a ) ^{-\kappa/2} \bigl( 1+\mathrm{o}(1) \bigr) =da^{-\kappa
/2} \bigl(
1+\mathrm{o}(1) \bigr) \quad\quad\mbox{as }a\rightarrow\infty.
\nonumber
\end{eqnarray}
Hence, the thresholds $\mathcal{C}_{n}$ can always be chosen as
$\mathcal{C}%
_{n} = d^{2/\kappa}(n/k_{n})^{2/\kappa}$. Now use an implication of
Karamata's theorem to obtain as $n \rightarrow \infty$ (e.g., Resnick
\cite{Resnick87}, Theorem~0.6):\footnote{Note if $\kappa = 4$ then for finite $a > 0$ there exists $K > %
0 $ such that $E[\mathcal{E}_{t}^{2}I_{n,t}^{(\mathcal{E})}] \sim %
\int_{0}^{\mathcal{C}_{n}^{2}}P(\mathcal{E}_{t} > u^{1/2})\,\mathrm{d}u = K + \int_{a}^{\mathcal{C}_{n}^{2}}P(\mathcal{E}_{t} > u^{1/2})\,\mathrm{d}u \sim %
K + d\int_{a+1}^{\mathcal{C}_{n}^{2}}u^{-\kappa/4}\,\mathrm{d}u \sim K + %
d\ln(\mathcal{C}_{n}^{4}) \sim d\ln(n)$ since $\mathcal{C}_{n} = %
K(n/k_{n})^{1/4}$.}%
%e15 #&#
\begin{eqnarray}\label{mom_e4}
 \kappa&=&4 : E \bigl[ \mathcal{E}_{t}^{2}I_{n,t}^{(\mathcal{E}%
)}
\bigr] \sim d\ln(n) ,\nonumber
\\[-8pt]\\[-8pt]
 \kappa&\in& ( 2,4 ) : E \bigl[ \mathcal{E}%
_{t}^{2}I_{n,t}^{(\mathcal{E})}
\bigr] \sim \biggl( \frac{\kappa}{%
4-\kappa} \biggr) \mathcal{C}_{n}^{2}P
\bigl( \llvert \mathcal{E}%
_{t}\rrvert >
\mathcal{C}_{n} \bigr) = \biggl( \frac{\kappa}{4-\kappa
}%
\biggr)
d^{4/\kappa} \biggl( \frac{n}{k_{n}} \biggr) ^{4/\kappa-1}=\mathrm{o}(n).
\nonumber
\end{eqnarray}

The following claim summarizes the above details.

%th2.5 #&#
%
\begin{theorem}[(Convergence rate)]\label{th:converg}
Under Assumptions \ref{ass1} and \ref{ass2} if $\kappa > 4$
then $\mathcal{V}_{n} \sim n(E[\epsilon_{t}^{4}] - %
1])^{-1}E[\mathfrak{s}_{t}\mathfrak{s}_{t}^{\prime}]$. If $\kappa
\leq 4$ then for $i = 1,\ldots,q$
%e16 #&#
\begin{eqnarray}\label{V}
\kappa&=&4 :\mathcal{V}_{i,i,n}^{1/2}\sim \biggl(
\frac{n}{\ln (
n ) } \biggr) ^{1/2}d^{-1/2} \bigl( E \bigl[
\mathfrak {s}_{i,t}^{2} \bigr] \bigr) ^{1/2},\nonumber
\\[-8pt]\\[-8pt]
\kappa&\in& ( 2,4 ) :\mathcal{V}_{i,i,n}^{1/2}\sim
n^{1/2} \biggl( \frac{k_{n}}{n} \biggr) ^{2/\kappa-1/2}d^{-2/\kappa}
\biggl( \frac{4-\kappa}{\kappa} \biggr) ^{1/2} \bigl( E \bigl[ \mathfrak
{s}_{i,t}^{2}%
\bigr] \bigr) ^{1/2}.
\nonumber
\end{eqnarray}
\end{theorem}

There are several key observations. First, as long as $\kappa \in %
(2,4)$ then elevating $k_{n}$ arbitrarily close to a fixed percent of $n$,
that is $k_{n} \approx \lambda n$ for $\lambda \in (1,0)$, will
optimize the convergence rate. This is logical  since large errors
adversely affect efficiency. In general this implies
%e17 #&#
\begin{equation}
k_{n}\sim n/g_{n}\quad\quad\mbox{for }g_{n}\rightarrow
\infty\mbox{ at a slow rate,} \label{kng}
\end{equation}
ensures $\mathcal{V}_{i,i,n}^{1/2} \sim n^{1/2}/g_{n}^{2/\kappa-1/2}$
for any $\kappa \in (2,4]$. Hence, $\mathcal{V}_{i,i,n}^{1/2} %
\rightarrow \infty$ can be driven as close to rate $n^{1/2}$ as we
choose by setting $g_{n} \rightarrow \infty$ very slowly (e.g.,
$g_{n} = \ln(\ln(n))$). Further, the rate monotonically
$n^{1/2}/g_{n}^{2/%
\kappa-1/2} \nearrow n^{1/2}$ as $\kappa \nearrow 4$. Hall and
Yao \cite{hallyao03} show the QML rate is $n^{1-2/\kappa}/L(n)$ for some slowly
varying $L(n) \rightarrow \infty$ and any $\kappa \in (2,4]$, hence
QMTTL can be assured to be faster for every $%
\kappa \in (2,4)$. Conversely, Peng and Yao's \cite{pengyao03}
Log-LAD and non-Gaussian QML are $n^{1/2}$-convergent %
(cf. Berkes and Horvath \cite{BerkesHorvath04}, Zhu and Ling
\cite{ZhuLing}), but the higher rate is not without
costs: (i) these estimators are not robust to error extremes in small
samples: see Section~\ref{sim}; (ii) Log-LAD requires $\ln\epsilon
_{t}^{2}$ to have a zero median; and (iii) non-Gaussian QML requires
additional moment conditions for Fischer consistency, for example,
WLQML requires $%
E|\epsilon_{t}| = 1$: see Section~1 for discussion.

Second, if $\kappa < 4$ and we use a fractile form $k_{n}\sim\lambda
n/g_{n}$ for slow $g_{n} \rightarrow \infty$ and $\lambda \in %
(0,1]$, then%
%e18 #&#
\begin{eqnarray}\label{rate_ng}
\frac{n^{1/2}}{g_{n}^{2/\kappa-1/2}} \bigl( \hat{\theta}_{n}-\theta ^{0}
\bigr) &\stackrel{d} {\rightarrow}&N \biggl( 0,\lambda^{- (
2/\kappa
-1/2 ) } \biggl(
\frac{\kappa}{4-\kappa} \biggr) d^{4/\kappa} \bigl( E%
\bigl[
\mathfrak{s}_{t}\mathfrak{s}_{t}^{\prime} \bigr]
\bigr) ^{-1} \biggr) \nonumber
\\[-8pt]\\[-8pt]
&=&N \bigl( 0,\mathcal{V} ( \lambda,\kappa,d ) \bigr)
. \nonumber
\end{eqnarray}
For example, in our simulation study we use $k_{n} \sim \lambda n/\ln
(n)$, hence $\hat{\theta}_{n}$ is $n^{1/2}/\break (\ln(n))^{2/\kappa-1/2}$%
-convergent with asymptotic variance $\mathcal{V}(\lambda,\kappa,d)$. The
asymptotic variance $\mathcal{V}(\lambda,\kappa,d)$ can always by
decreased by increasing $\lambda$ and therefore removing more extremes per
sample.

Third, in view of $k_{n} = k_{2,n}$ by (\ref{CU}), trimming rule (\ref%
{kng}) only concerns the amount of trimmed positive observations of $%
\mathcal{E}_{t} = \epsilon_{t}^{2} - 1$: the left tail of $%
\mathcal{E}_{t}$ is bounded, hence only the rate of right tail trimming
of $%
\mathcal{E}_{t}$ matters for the \textit{convergence rate}. In terms of
\textit{identification}, however, as discussed in Section~\ref{ver_2} the
number of trimmed left and right tail observations $k_{1,n}$ and $k_{2,n}$
must be balanced when $\epsilon_{t}$ is governed by a heavy tailed
distribution. For example, if $P(|\epsilon_{t}| > c) = (1 + %
c)^{-\kappa}$ with $\kappa \in (2,4)$, and $k_{1,n} \sim %
\lambda n/\ln(n)$, both as in our simulation study, then Assumption~\ref{ass2} holds
when $k_{2,n} \sim Kk_{1,n}^{\kappa/(\kappa-2)}/n^{\kappa/(\kappa
-2)-1} \sim Kn/(\ln(n))^{\kappa/(\kappa-2)}$, hence from (\ref%
{rate_ng}) the rate of convergence is $n^{1/2}/((\ln(n))^{\kappa
/(\kappa
-2)})^{2/\kappa-1/2} = n^{1/2}/(\ln(n))^{(4-\kappa)/[2(\kappa-2)]}$.

As a practical matter, naturally too much trimming in any given sample can
lead to small sample bias in $\hat{\theta}_{n}$. In Section~\ref{sim},
we use
$k_{n} \sim \lambda n/\ln(n)$ with $\lambda = 0.025$ for both very
thin and thick tailed error distributions: values much larger than $0.025$
(e.g., $\lambda=0.10$) leads to substantial bias, and values much smaller
(e.g., $\lambda=0.01$) are not effective for rendering $\hat{\theta}_{n}$
approximately normal in small samples. In general any value $\lambda
\in [0.02,0.05]$ leads to roughly the same results. Similar trimming
schemes are
found to be highly successful in other robust estimation and inference
contexts: see Hill \cite{Hilltest,Hillltts} and Hill and Aguilar \cite
{HillAguilar13}.

Last, there are several proposed methods in the robust statistics
literature for selecting trimming parameters like $\lambda$, but in this
literature the seeming universal approach for data transformations
involve a
fixed quantile threshold hence $k_{n} \sim \lambda n$ %
(cf. Huber \cite{Huber1964}, Hampel \textit{et~al}. \cite
{Hampleetal86}, Jureckova and Sen \cite{JurSen96}). Such methods include
covariance determinant or asymptotic variance minimization where a unique
internal solution for $\lambda$ exists. These methods are ill posed here
since they lead to corner solutions: consider that minimizing $\mathcal
{V}%
(\lambda,\kappa,d)$ above on $\lambda \in [\underline{\lambda},\bar
{%
\lambda}]$ leads to $\lambda = \bar{\lambda}$. See Hill and Aguilar
\cite{HillAguilar13}
for references and simulation evidence. In terms of inference more choices
exist, including test statistic functionals over $\lambda$ like the
supremum, and empirical process techniques for $p$-value computation %
(see Hill \cite{Hilltest}).

%s3 #&#
\section{Method of moments with re-centering}\label{sec:mm}

Our second estimator uses the method of moments based on negligibly weighted
errors imbedded in a QML score equation. This gives us the advantage of
re-centering to ensure identification. It therefore allows us to use a
greater variety of error transforms, as well as symmetric transforms
even if
the errors have an asymmetric distribution. Define $\Im_{t} := \sigma
(y_{\tau} \dvtx  \tau \leq t)$.

The class of transformations we consider have the general form%
%e19 #&#
\begin{equation}
\psi ( u,c ) :=u\times\varpi(u,c)\times I \bigl( \llvert u\rrvert
 \leq c \bigr)
, \label{psi}
\end{equation}
where $\varpi(\cdot,c)$ is for each $c$ a Borel function, and%
%e20 #&#
\begin{equation}
\lim_{c\rightarrow\infty}\varpi(u,c)\times I \bigl( \llvert u\rrvert \leq c
\bigr) =1. \label{w_lim}
\end{equation}
Thus, $\psi(u,c)$ is a \textit{redescending} function (see Andrews
\textit{et~al}.
\cite{Andrewsetal72} and Hampel \textit{et~al}. \cite{Hampleetal86}).
In the literature typically $c$ is fixed, but the only way we can
identify $%
\theta^{0}$ and obtain Fischer consistency without an additional simulation
step is to enforce $c \rightarrow \infty$ as $n \rightarrow %
\infty$.\footnote{%
See, for example, Sakata and White \cite{SakataWhite98}, Cantoni and Ronchetti
\cite{CantoniRonchetti} and Mancini \textit{et~al}. \cite
{Mancinietal05}.} Notice as $c \rightarrow \infty$ the transform
satisfies $\psi(u,c) \rightarrow u$ hence it applies a negligible
weight to $u$. Further, it operates similar to tail-trimming since by
(\ref%
{w_lim})%
%e21 #&#
\begin{equation}
\psi ( u,c ) =uI \bigl( \llvert u\rrvert \leq c \bigr) \times \bigl( 1+\mathrm{o} ( 1 )
\bigr) \quad\quad\mbox{as }c\rightarrow\infty. \label{psi_I}
\end{equation}

We focus on two types of weights $\varpi$. First, the simple trimming
case $%
\psi(u,c) = uI(|u| \leq c)$, hence
\[
\varpi(u,c)=1.
\]
The theory developed below easily extends to related redescending
functions $%
\psi(u,c)$, like Hampel's three-part trimming function with thresholds
$0 %
< a < b < c$ (see Andrews \textit{et~al}. \cite{Andrewsetal72}):
\[
\cases{ u, &\quad $0\leq\llvert u\rrvert \leq a$,
\cr
a\times \operatorname{sign}(u), &\quad $a<\llvert u
\rrvert \leq b$,
\cr
\displaystyle\frac{a\times ( c-\llvert  u\rrvert  ) }{c-b}\times \operatorname{sign}(u), &\quad $b<\llvert u\rrvert \leq c$,\cr
0 ,&\quad$ c<\llvert u\rrvert $.
}
\]
This can be identically written as (\ref{psi}) with%
%e22 #&#
\begin{equation}
\varpi(u,c)=I \bigl( \llvert u\rrvert \leq a \bigr) +\frac{a}{%
\llvert  u\rrvert }\times I
\bigl( a<\llvert u\rrvert \leq b \bigr) +\dfrac{a ( c-\llvert  u\rrvert  ) }{\llvert
u\rrvert  ( c-b ) }\times I \bigl( b<
\llvert u\rrvert \leq c \bigr) . \label{Hamp}
\end{equation}
Of course, we abuse notation since there are three thresholds $\{a,b,c\}$.
By construction $\varpi(u,c) \in [0,1]$, while negligibility requires
the smallest threshold $a \rightarrow \infty$, hence $\varpi(u,c) %
\rightarrow 1$ as $a \rightarrow \infty$.

Second, we use smooth weights $\varpi(u,c)$ that are continuously
differentiable in $c$, with%
%e23 #&#
\begin{equation}
\biggl\llvert \frac{\partial}{\partial c}\varpi(u,c)\biggr\rrvert \times I \bigl( \llvert
u\rrvert \leq c \bigr) \leq K\frac{1}{c}. \label{dh}
\end{equation}
Notice the simple trimming case $\varpi(u,c) = 1$ trivially satisfies
 (\ref{dh}). Thus, as $c \rightarrow \infty$ the transform derivative $%
(\partial/\partial c)\psi ( u,c ) \rightarrow 0$ at rate $%
\mathrm{O}(1/c)$ for all $|u| \neq c$. An example is Tukey's
bisquare $\varpi
(u,c) = (1 - (u/c)^{2})^{2}$ with $(\partial/\partial c)\varpi
(u,c) = 2(1 - (u/c)^{2})u^{2}/c^{3}$ hence (\ref{dh}) holds. A
second example is the exponential $\varpi(u,c) = \exp\{-|u|/c\}$ with
$%
(\partial/\partial c)\varpi(u,c) = \exp\{-|u|/c\}|u|/c^{2}$.

\begin{Assumption}[(Redescending transforms)]\label{ass4}
Let $\psi(u,c)$
satisfy (\ref{psi}), (\ref{w_lim}) and (\ref{dh}).
\end{Assumption}

Now define two-tailed observations $\epsilon_{t}^{(a)}(\theta) := %
|\epsilon_{t}(\theta)|$ and their order statistics $\epsilon
_{(1)}^{(a)}(\theta) \geq \epsilon_{(2)}^{(a)}(\theta) \geq %
\cdots$\,, and let $\{k_{n}\}$ be an intermediate order sequence. Write
\begin{eqnarray*}
 \hat{I}_{n,t}^{(\epsilon)}(\theta)&:=&I \bigl( \bigl\llvert \epsilon
_{t}(\theta)\bigr\rrvert \leq\epsilon_{(k_{n})}^{(a)}(
\theta) \bigr),
\\
 \hat{\psi}_{n,t}(\theta)&:=&\psi \bigl( \epsilon_{t}(
\theta),\epsilon _{(k_{n})}^{(a)}(\theta) \bigr) =
\epsilon_{t}(\theta)\times\varpi \bigl( \epsilon_{t}(
\theta),\epsilon_{(k_{n})}^{(a)}(\theta) \bigr) \hat{I} %%
_{n,t}^{(\epsilon)}(\theta),
\end{eqnarray*}
and define re-centered equations and a Method of Negligibly-Weighted
Moments
(MNWM) estimator
\begin{eqnarray*}
\widehat{\check{m}}_{n,t}(\theta)&:=& \Biggl( \hat{\psi}_{n,t}^{2}(
\theta )-%
\frac{1}{n}\sum_{t=1}^{n}
\hat{\psi}_{n,t}^{2}(\theta) \Biggr) \times
\mathfrak{s}_{t}(\theta)\quad\quad\mbox{and}
\\
\hat{\theta }_{n}^{(m)}&:=&
\operatorname{arg\,min}\limits
_{\theta\in\Theta} \Biggl( \sum_{t=1}^{n}
\widehat{\check{m}}%
_{n,t}(\theta) \Biggr) ^{\prime}
\Biggl( \sum_{t=1}^{n}\widehat{\check
{m}}%
_{n,t}(\theta) \Biggr) .
\end{eqnarray*}
Any positive definite symmetric weight matrix $W\in\mathbb{R}^{q\times q}$
leads to the same solution $\operatorname{arg\,min}_{\theta\in\Theta
}\sum_{t=1}^{n}\widehat{\check{m}}_{n,t}(\theta)^{\prime} \times W %
\times \sum_{t=1}^{n}\widehat{\check{m}}_{n,t}(\theta)$. Similarly,
any $%
\Im_{t-1}$-measurable uniformly $L_{2+\iota}$-bounded vector
$z_{t}(\theta
) \in \mathbb{R}^{r}$, $r \geq q$, can be used instead of $%
\mathfrak{s}_{t}(\theta)$ for a GMM estimator (Hansen \cite
{hansen82}). The scaled
volatility derivative $\mathfrak{s}_{t}(\theta)$, however, provides an
analogue to QML. Finally, as discussed in Section~\ref{sec:qmttl} small sample performance
appears to be improved if we also trim by $y_{t-1}$, while asymptotics are
unchanged if trimming is negligible. The estimator in this case uses
the transformed error $\epsilon_{t}(\theta)\varpi(\epsilon_{t}(\theta),\epsilon
_{(k_{n})}^{(a)}(\theta))\hat{I}_{n,t}^{(\epsilon)}(\theta)\hat{I}%
_{n,t-1}^{(y)}$.

Next, for asymptotics let $\{\mathcal{C}_{n}(\theta)\}$ satisfy
\[
P \bigl( \bigl\llvert \epsilon_{t}(\theta)\bigr\rrvert \geq
\mathcal{C}%
_{n}(\theta) \bigr) =\frac{k_{n}}{n},
\]
write compactly
\begin{eqnarray*}
I_{n,t}^{(\epsilon)}(\theta)&:=&I \bigl( \bigl\llvert
\epsilon_{t}(\theta )\bigr\rrvert \leq\epsilon_{(k_{n})}^{(a)}(
\theta) \bigr) ,
\\
 \psi _{n,t}(\theta)&:=&\psi \bigl( \epsilon_{t}(
\theta),\mathcal {C}_{n}(\theta ) \bigr) \quad\mbox{and}\quad
\epsilon_{n,t}(\theta):=\epsilon_{t}(\theta
)I_{n,t}^{(\epsilon)}(\theta),
\end{eqnarray*}
and define equations with non-random thresholds
\[
\check{m}_{n,t}(\theta):= \bigl( \psi_{n,t}^{2}(
\theta)-E \bigl[ \psi _{n,t}^{2}(\theta) \bigr] \bigr) \times
\bigl( \mathfrak {s}_{t}(\theta)-E%
\bigl[
\mathfrak{s}_{t}(\theta) \bigr] \bigr) .
\]
In view of re-centering in $\widehat{\check{m}}_{n,t}(\theta)$ it can be
shown that, asymptotically, $\check{m}_{n,t}$ and $\widehat{\check{m}}_{n,t}$
are interchangeable. See the \hyperref[app]{Appendix}.

Since $\epsilon_{t}$ is i.i.d. and has a smooth distribution, the transform
is negligible in that $\psi_{n,t}(\theta) \stackrel{\mathrm{a.s.}}{\rightarrow
} %
\epsilon_{t}(\theta)$, and $\mathfrak{s}_{t}$ is $\Im_{t-1}$-measurable,
it follows for all $n \geq N$ and some large $N \in \mathbb{N}$
\[
E \bigl[ \check{m}_{n,t}(\theta)|\Im_{t-1} \bigr] =0\quad\mbox{if and only if}\quad\theta=\theta^{0},
\]
hence an identification condition like Assumption~\ref{ass2} automatically holds.
Similarly, by negligibility $\psi(u,c) = uI(|u| \leq c) \times (1 +
\mathrm{o}(1))$ as $c \rightarrow \infty$ and $E [ \epsilon
_{n,t}^{2} ] \rightarrow 1$,\vadjust{\goodbreak} hence by independence of the errors
\begin{eqnarray*}
E \bigl[ \check{m}_{n,t}\check{m}_{n,t}^{\prime}
\bigr] &=&E \bigl( \bigl( \psi ( \epsilon_{t},\mathcal{C}_{n}
) ^{2}-E \bigl[ \psi ( \epsilon_{t},\mathcal{C}_{n}
) ^{2} \bigr] \bigr) ^{2} \bigr) \times E%
\bigl[
\bigl( \mathfrak{s}_{t}-E [ \mathfrak{s}_{t} ] \bigr) \bigl(
\mathfrak{s}_{t}-E [ \mathfrak{s}_{t} ] \bigr)
^{\prime}%
\bigr]
\\
&=&E \bigl( \bigl( \epsilon_{n,t}^{2}-E \bigl[
\epsilon_{n,t}^{2} \bigr] \bigr) ^{2} \bigr) \times
E \bigl[ \bigl( \mathfrak{s}_{t}-E [ \mathfrak{%
s}_{t} ] \bigr) \bigl( \mathfrak{s}_{t}-E [ \mathfrak
{s}_{t} ] \bigr) ^{\prime} \bigr] \times \bigl( 1+\mathrm{o}(1) \bigr)
\\
&=& \bigl( E \bigl[ \epsilon_{n,t}^{4} \bigr] -1 \bigr)
\times E \bigl[ \bigl( \mathfrak{s}_{t}-E [ \mathfrak{s}_{t}
] \bigr) \bigl( \mathfrak{s}%
_{t}-E [ \mathfrak{s}_{t}
] \bigr) ^{\prime} \bigr] \times \bigl( 1+\mathrm{o}(1) \bigr) .
\end{eqnarray*}
The MNWM scale is therefore%
%e24 #&#
\begin{equation}\label{V_MTTM}
\mathcal{\accentset{\sdbullet}{V}}_{n}=\frac{n}{E [ \epsilon_{n,t}^{4} ]
-1}%
\times E
\bigl[ \bigl( \mathfrak{s}_{t}-E [ \mathfrak{s}_{t} ] \bigr)
\bigl( \mathfrak{s}_{t}-E [ \mathfrak{s}_{t} ] \bigr)
^{\prime} \bigr] ,
\end{equation}
which is positive definite under Assumption~\ref{ass1}.

%th3.1 #&#
%
\begin{theorem}[(MNWM)]\label{th:mnwm}
Under Assumptions \ref{ass1} and \ref{ass4} $\mathcal{\accentset{\sdbullet}{V}}%
_{n}^{1/2}(\hat{\theta}_{n}^{(m)} - \theta^{0}) \stackrel{d}{%
\rightarrow} N(0,I_{q})$. Further each $\mathcal{\accentset{\sdbullet}{V}}%
_{i,i,n} \rightarrow \infty$ and $\mathcal{\accentset{\sdbullet}{V}}_{i,i,n}/%
\mathcal{V}_{i,i,n} \rightarrow (0,1)$.
\end{theorem}

%re7 #&#
%
\begin{remark}
In general a direct comparison of QMTTL and MNWM scales $\mathcal{V}%
_{n}$ and $\mathcal{\accentset{\sdbullet}{V}}_{n}$ is difficult for a particular
$n$ due
to the different trimming strategies. Notice, however, that $E[\epsilon
_{n,t}^{4}] - 1 =E[\mathcal{E}_{t}^{2}I_{n,t}^{(\mathcal{E})}] %
\times (1 + \mathrm{o}(1))$ if $\mathcal{C}_{n} = (\mathcal{U}_{n} + %
1)^{1/2}$. This follows by noting $E[\epsilon_{t}^{2}] = 1$, $\mathcal
{E%
}_{t}^{2} \in [-1,\infty)$, negligibility and dominated convergence
imply
\begin{eqnarray*}
E \bigl[ \mathcal{E}_{t}^{2}I_{n,t}^{(\mathcal{E})}
\bigr] &=&E \bigl[ \bigl( \epsilon_{t}^{4}-2
\epsilon_{t}^{2}+1 \bigr) I \bigl( -\mathcal
{L}_{n}\leq \epsilon_{t}^{2}-1\leq
\mathcal{U}_{n} \bigr) \bigr]
\\
&=&E \bigl[ \epsilon_{t}^{4}I \bigl( ( 1-
\mathcal{L}_{n} ) ^{1/2}\leq\llvert \epsilon_{t}
\rrvert \leq ( \mathcal{U}%
_{n}+1 ) ^{1/2} \bigr)
\bigr] \times \bigl( 1+\mathrm{o}(1) \bigr)
\\
&=&E \bigl[ \epsilon_{t}^{4}I \bigl( \llvert
\epsilon_{t}\rrvert \leq ( \mathcal{U}_{n}+1 )
^{1/2} \bigr) \bigr] \times \bigl( 1+\mathrm{o}(1) \bigr) .
\end{eqnarray*}
Thus, $\mathcal{\accentset{\sdbullet}{V}}_{n} \times \mathcal{V}_{n}^{-1} = E[(%
\mathfrak{s}_{t} - E[\mathfrak{s}_{t}])(\mathfrak{s}_{t} - E[%
\mathfrak{s}_{t}]^{\prime}] \times E[\mathfrak{s}_{t}\mathfrak{s}%
_{t}^{\prime}]$ as $n \rightarrow \infty$. Therefore $\mathcal{%
\accentset{\sdbullet}{V}}_{n}$ is smaller than $\mathcal{V}_{n}$ due to the centered
term $E[(\mathfrak{s}_{t} - E[\mathfrak{s}_{t}])(\mathfrak{s}_{t} - %
E[\mathfrak{s}_{t}])^{\prime}]$, hence identification is assured at a cost
of efficiency.
\end{remark}

%re8 #&#
%
\begin{remark}
Since $ \mathcal{\accentset{\sdbullet}{V}}_{n} \sim \mathcal{K}_{n}\mathcal{V}%
_{n}$ for some sequence of positive definite matrices $\{\mathcal
{K}_{n}\}$,
the Section~\ref{sec:rate} discourse on the QMTTL rate of convergence
carries over here.
\end{remark}

%s4 #&#
\section{Inference}\label{sec:infer}

In view of $\mathcal{V}_{n} \sim n(E[\mathcal{E}_{t}^{2}I_{n,t}^{(%
\mathcal{E})}])^{-1}E[\mathfrak{s}_{t}\mathfrak{s}_{t}^{\prime}]$, a
natural estimator of the QMTTL scale $\mathcal{V}_{n}$ is%
%e25 #&#
\begin{equation}\label{V_hat}
\mathcal{\hat{V}}_{n}=\mathcal{\hat{V}}_{n}(\hat{
\theta}_{n})=n\times \frac{1%
}{1/n\sum_{t=1}^{n}1/n\sum_{t=1}^{n}\mathcal{E}_{t}^{2}(\hat{\theta
}_{n})%
\hat{I}_{n,t}^{(\mathcal{E})}(\hat{\theta}_{n})}\times\frac{1}{n}%
\sum_{t=1}^{n}\mathfrak{s}_{t}(
\hat{\theta}_{n})\mathfrak {s}_{t}^{\prime}(%
\hat{\theta}_{n}).
\end{equation}

%th4.1 #&#
%
\begin{theorem}
\label{th:scale} Under Assumptions \ref{ass1} and \ref{ass2} $\mathcal{\hat{V}}_{n} = %
\mathcal{V}_{n}(1 + \mathrm{o}_{p}(1))$.\vadjust{\goodbreak}
\end{theorem}

%re9 #&#
%
\begin{remark}
Notice $\mathcal{\hat{V}}_{n} = \mathcal{V}_{n}(1 + \mathrm{o}_{p}(1))$
only reduces to $\mathcal{\hat{V}}_{n} = \mathcal{V}_{n} + \mathrm{o}_{p}(1)$
when $E[\epsilon_{t}^{4}] < \infty$. In general classic inference is
available without knowing the true rate of convergence, nor even if trimming
is required.
\end{remark}

%re10 #&#
%
\begin{remark}
A consistent estimator of the MNWM scale $\mathcal{\accentset{\sdbullet}{V}}_{n}$
can similarly be constructed.
\end{remark}

A Wald statistic naturally follows for a test of (non)linear parameter
restrictions $R(\theta^{0}) = 0$ where $R \dvtx  \mathbb{R}^{q} %
\rightarrow \mathbb{R}^{J}$ and $J \geq 1$. Assume $R$ is
differentiable with a gradient $\mathcal{D}(\theta) = (\partial
/\partial\theta)R(\theta)$ that is continuous, differentiable and has
full column rank. The test statistic with the QMTTL estimator as a plug-in
is
\[
\mathcal{W}_{n}=R(\hat{\theta}_{n})^{\prime} \bigl(
\mathcal{D}(\hat {\theta}%
_{n})\mathcal{\hat{V}}_{n}^{-1}(
\hat{\theta}_{n})\mathcal{D}(\hat{\theta }%
_{n})^{\prime} \bigr) ^{-1}R(\hat{\theta}_{n}).
\]
Use Theorems \ref{th:norm} and \ref{th:scale} to deduce $\mathcal{W}_{n} \stackrel{d}{%
\rightarrow} \chi^{2}(J)$ under the null, and if $R(\theta^{0})
\neq 0$ then \mbox{$\mathcal{W}_{n} \stackrel{p}{\rightarrow} \infty$}.

Similarly, the proof of Theorem~\ref{th:consist} shows the QMTTL first order
condition is $1/n\*\sum_{t=1}^{n}m_{t}(\hat{\theta}_{n})I_{n,t}^{(\mathcal
{E}%
)}(\hat{\theta}_{n}) = 0$ a.s. This naturally suggests the possibility
of a score or Lagrange Multiplier test since a QMTTL estimator under the
constraint $R(\theta^{0}) = 0$, denoted $\hat{\theta}_{n}^{(c)}$, also
satisfies $1/n\sum_{t=1}^{n}m_{t}(\hat{\theta
}_{n}^{(c)})I_{n,t}^{(\mathcal{E%
})}(\hat{\theta}_{n}^{(c)}) \stackrel{p}{\rightarrow} 0$ if the
constraint is true. A heavy tail robust test of $R(\theta^{0}) = 0$ can
therefore be coached as a tail-trimmed moment condition test as in Hill
and Aguilar \cite{HillAguilar13}.

%s5 #&#
\section{Simulation}\label{sim}

We now compare our robust QML and Method of Moments estimators with various
estimators in the literature. In order to draw the best comparisons between
QMTTL and MNWM, we initially focus on simple trimming for MNWM. We compare
our estimators to QML as a benchmark, as well as Log-LAD, Weighted Laplace
QML (WLQML) and Power-Law QML (PQML) due to their heavy tail robustness
properties. Finally, we investigate other redescending transforms as
alternatives for MNWM, and whether tail-trimming can improve the small
sample properties of PQML.

%s5.1 #&#
\subsection{Data generation and estimators}

Let $P_{\kappa}$ denote a symmetric Pareto distribution: if $\epsilon_{t}$
is distributed $P_{\kappa}$ then $P(\epsilon_{t} \leq -a) = %
P(\epsilon_{t} \geq a) = 0.5(1 + a)^{-\kappa}$ for $a > 0$%
. We draw $20n$ observations for $n \in \{100,800\}$ from the GARCH
process $y_{t} = \sigma_{t}\epsilon_{t}$ and $\sigma_{t}^{2} = %
0.05 + 0.05y_{t-1}^{2} + 0.90\sigma_{t-1}^{2}$ with a starting value $%
\sigma_{1}^{2} = 0.05$, and retain the last $n$ observations for the
sample. This is repeated to produce 10,000 samples $\{y_{t}\}_{t=1}^{n}$.
Our choice of parameter values are indicative of values we obtain in the
empirical study below, and frequently encountered in macroeconomic and
financial data. The error $\epsilon_{t}$ is i.i.d. $N(0,1)$, or $P_{2.5}$
standardized such that $E[\epsilon_{t}^{2}] = 1$.

We compute the feasible QMTTL and MNWM estimators conditional on the first
observation, with parameter space is $\Theta = [\iota,2] \times %
[\iota,1-\iota] \times [\iota,1-\iota]$ where $\iota = %
10^{-10} $. The iterated volatility variable is $h_{1}(\theta) = \omega
$ and $h_{t}(\theta) = \omega + \alpha y_{t-1}^{2} + \beta
h_{t-1}(\theta)$ where we initialize $h_{1}(\theta) = \omega$ for
QMTTL and $h_{1}^{\theta}(\theta) = [1,0,0]^{\prime}$ for MNWM.

As a benchmark for QMTTL we use \textit{strong} \textit{asymmetric} trimming
with error fractiles $k_{2,n} = \max\{1,[0.025n/\ln(n)]\}$ and
$k_{1,n} = 35k_{2,n}$. This equates to $\{k_{1,n},k_{2,n}\} = \{1,35\}$
and $%
\{3,105\}$ for $n = 100$ and $800$. The fractile for trimming by $%
y_{t-1} $ is $\tilde{k}_{n} = \max\{1,[0.1\ln(n)]\}$: asymptotics do
not require such trimming, while removing a very few criterion
equations due
to large $y_{t-1}$ appears to improve the estimator's performance. The
benchmark for MNWM is simple trimming $\psi(u,c) = uI(|u| \leq c)$%
. The error fractile is as above $k_{n} = \max\{1,[0.025n/\ln(n)]\}$
and the fractile for trimming by $y_{t-1}$ is again $\tilde{k}_{n}$.

In addition to the benchmark estimates, we compute MNWM with Tukey's
bisquare and exponential transforms. We also compute QMTTL with \textit{weak}
\textit{asymmetric} ($k_{1,n} = 10k_{2,n}$) and \textit{symmetric} ($%
k_{1,n} = k_{2,n}$) trimming. Recall from Section~\ref{ver_2} that for
QMTTL $k_{1,n} = 35k_{2,n}$ roughly minimizes bias in the Pareto case $%
P(|\epsilon_{t}| \geq a) = (1 + a)^{-2.5}$ when $n = 100$.
We show here that using $k_{1,n} = 35k_{2,n}$ even when $n = 800$
still promotes a sharp estimator. In simulations not reported here, we find
that the bias minimizing relation $k_{1,n} = 100k_{2,n}$ when $%
P(|\epsilon_{t}| \geq a) = (1 + a)^{-2.5}$ and $n = 800$
logically leads to even smaller bias, but bias is still low when
$k_{1,n} %
= 35k_{2,n}$. Recall\vspace*{2pt} also that any combination $\{k_{1,n},k_{2,n}\}$
works in the Gaussian case provided $k_{i,n} = \mathrm{o}(n^{1/2})$. This is
violated here since we use $k_{i,n} \sim Kn/\ln(n)$, however this
matters only asymptotically, and we demonstrate that using $k_{2,n}
\sim Kn/\ln(n)$ and $k_{1,n} = 35k_{2,n}$ for $n = 100$ and $800$ in
the thin tail case still leads to a competitive estimator in small samples.
Indeed, if we use $k_{i,n} \sim Kn^{1/2}/\ln(n)$ then the small sample
performance is essentially identical to what we see here.

Peng and Yao's \cite{pengyao03} Log-LAD criterion is $\sum_{t=2}^{n}|\ln y_{t}^{2} - \ln h_{t}(\theta)|$. The WLQML criterion is
$\sum_{t=2}^{n}\{\ln
h_{t}^{1/2}(\theta) + |y_{t}/h_{t}^{1/2}(\theta)|\}w_{t}$ where we
choose the weights $\{w_{t}\}$ as in Zhu and Ling \cite{ZhuLing},
equation (2.4): $w_{t} = (\max\{1,C^{-1}\sum_{i=1}^{\infty
}i^{-9}|y_{t-i}I(|y_{t-i}| > C)|\})^{-4}$ where $C = %
y_{(0.10n)}^{(a)} $ and $y_{t-i} = 0\ \forall i \geq t$.

\begin{sidewaystable*}
\tablewidth=\textwidth
\caption{Simulation estimation results for $%
\theta _{3}^{0}$}\label{tab1}
\begin{tabular*}{\textwidth}{@{\extracolsep{\fill}}lllllllllllll@{}}
\hline
& \multicolumn{6}{l}{$\epsilon _{t}\sim  \bar{P}_{2.5}$} &
\multicolumn{6}{l@{}}{$\epsilon _{t}\sim  N(0,1)$} \\
& \multicolumn{6}{l}{\hrulefill} &
\multicolumn{6}{l@{}}{\hrulefill} \\
 & \multicolumn{3}{l}{$n=100$}    & \multicolumn{3}{l}{$n=800$}   &
\multicolumn{3}{l}{$n=100$} &  \multicolumn{3}{l@{}}{$n=800$}
\\
& \multicolumn{3}{l}{\hrulefill}    & \multicolumn{3}{l}{\hrulefill}   &
\multicolumn{3}{l}{\hrulefill} &  \multicolumn{3}{l@{}}{\hrulefill}
\\
 & Bias & RMS\tabnoteref{a}  & KS\tabnoteref{b} &
  Bias & RMS & KS &    Bias & RMS & KS &    Bias & RMS
& KS \\
\hline
 QMTTL-SA\tabnoteref{c}  & $-$0.010 &0.092 & 1.75 &
  \hphantom{$-$}0.008 &0.045 & 1.45 &   $-$0.063 &0.095 & 3.87 &
\hphantom{$-$}0.001 &0.030 & 1.07 \\
 QMTTL-WA  & $-$0.031 &0.102 & 3.01 &  \hphantom{$-$}0.024 &0.038 & 2.76   & $-$0.060
 &0.089 & 4.76   &\hphantom{$-$}0.003 &0.030 & 1.34 \\
 QMTTL-S  & $-$0.041 &0.114 & 4.69
& \hphantom{$-$}0.016 &0.044 & 4.21 &   $-$0.069 &0.075 & 6.30
 &\hphantom{$-$}0.005 &0.032 & 1.89 \\
 MNWM-I\tabnoteref{d}  & $-$0.023 &0.111 & 2.31 &
  $-$0.010 &0.064 & 1.56 &   $-$0.029 &0.108 & 2.86 &
$-$0.008 &0.036 & 1.17 \\
 MNWM-T  & $-$0.019 &0.103 & 2.87   &
$-$0.012 &0.069 & 1.61 &  \hphantom{$-$}0.021 &0.113 & 3.16 &   $-$0.010 &0.039 & 1.20 \\
 MNWM-E  & $-$0.025 &0.117 & 3.13 &
$-$0.016 &0.058 & 1.50 &   $-$0.026 &0.097 & 3.02 &   $-$0.013 &0.037 & 1.31 \\
 WLQML\tabnoteref{e}  & $-$0.063 &0.124 & 5.92 &
  $-$0.135 &0.107 & 7.64 &   $-$0.092 &0.082 & 8.12 &
$-$0.088 &0.084 & 6.05 \\
 WLQML$_{E|\epsilon _{t}|=1}$  & $-$0.082 &0.219 & 8.48 &
  $-$0.072 &0.089 & 5.64 &   $-$0.075 &0.078 & 9.36 &
$-$0.065 &0.067 & 3.97 \\
 PQML$_{3}$  & $-$0.048 &0.085 & 6.17   & $-$0.039 &0.059 & 3.00 &   $-$0.065 &0.067 & 9.07 &
  \hphantom{$-$}0.005 &0.032 & 1.30 \\
 PQML$_{3.5}$  & $-$0.034 &0.083 & 4.74 &
  $-$0.018 &0.056 & 3.17 &   $-$0.064 &0.062 & 9.54 &
\hphantom{$-$}0.009 &0.029 & 2.75 \\
 PQMTTL$_{3.5}^{\mathrm{SA}}$  & $-$0.054 &0.116 & 6.23 &
  $-$0.017 &0.056 & 2.23 &  $-$0.061 &0.074 & 6.38 &
\hphantom{$-$}0.011 &0.028 & 2.65 \\[3pt]
PQMTTL$_{3.5}^{\mathrm{WA}}$ & $-$0.031 &0.074 & 4.28 &
 $-$0.012 &0.046 & 1.35 &   $-$0.051 &0.074 & 8.21 &
\hphantom{$-$}0.008 &0.028 & 2.12 \\[3pt]
PQMTTL$_{3.5}^{\mathrm{S}}$ & $-$0.027 &0.077 & 4.05 &
 $-$0.019 &0.057 & 2.43 &   $-$0.055 &0.069 & 8.30 &
\hphantom{$-$}0.011 &0.027 & 2.76 \\
Log-LAD & $-$0.217 &0.165 & 9.88
& $-$0.253 &0.149 & 9.12 &   $-$0.082 &0.100 & 7.01   & $-$0.019 &0.046 & 3.61 \\
QML & $-$0.073 &0.099 & 6.23  &
$-$0.054 &0.078 & 4.65 &   $-$0.112 &0.089 & 8.71   & $-$0.013 &0.034 & 1.64 \\
\hline
\end{tabular*}
\tabnotetext[a]{a}{The square root of the empirical mean squared error.}

\tabnotetext[b]{b}{The Kolmogorov--Smirnov statistic divided by the 5\%
critical value: KS $>1$ indicates rejection of normality at the 5\% level.}
\tabnotetext[c]{c}{\textit{Benchmark} QMTTL-SA (strong asymmetric) uses
fractiles $k_{1,n} = 35k_{2,n}$; QMTTL-WA (weak asymmetric)
uses $k_{1,n} = 10k_{2,n}$; QMTTL-S
(symmetric) uses $k_{1,n} = k_{2,n}$.}
\tabnotetext[d]{d}{\textit{Benchmark} MNWM-I uses the simple trimming
function $\psi (u,c) = uI(|u| \leq  c)$; MNWM-T and MNWM-E use
Tukey's bisquare and exponential transforms.}
\tabnotetext[e]{e}{WLQML is Weighted Laplace QML. WLQML$_{E|\epsilon
_{t}|=1}$ is WLQML for processes with $E|\epsilon _{t}|=1$. PQML$_{\vartheta
}$ is power-law QML with criterion index $\vartheta $.
PQMTTL$_{\vartheta }^{\mathrm{WA}}$ and PQMTTL$_{\vartheta }^{\mathrm{SA}}$ are
tail-trimmed PQML with \textit{weak} \textit{asymmetric} ($k_{1,T} = 5k_{2,T}$)
or \textit{strong asymmetric} ($k_{1,T} = 9k_{2,T}$) trimming.}
\end{sidewaystable*}

\begin{sidewaystable*}
\tablewidth=\textwidth
\caption{Test rejection frequencies\protect\tabnoteref{a1} at 5\% level
for $\theta _{3}^{0}$}\label{tab2}
\begin{tabular*}{\textwidth}{@{\extracolsep{\fill}}lllllllllllll@{}}
 \hline
& \multicolumn{6}{l}{$\epsilon _{t}\sim  \bar{P}_{2.5}$} &
\multicolumn{6}{l@{}}{$\epsilon _{t}\sim  N(0,1)$} \\
& \multicolumn{6}{l}{\hrulefill} &
\multicolumn{6}{l@{}}{\hrulefill} \\
 & \multicolumn{3}{l}{$n=100$}    & \multicolumn{3}{l}{$n=800$}   &
\multicolumn{3}{l}{$n=100$} &  \multicolumn{3}{l@{}}{$n=800$}
\\
& \multicolumn{3}{l}{\hrulefill}    & \multicolumn{3}{l}{\hrulefill}   &
\multicolumn{3}{l}{\hrulefill} &  \multicolumn{3}{l@{}}{\hrulefill}
\\
& $H_{0}$ & $H_{1}^{1}$ & $H_{1}^{2}$ &   $H_{0}$ & $H_{1}^{1}$ & $H_{1}^{2}
$ &   $H_{0}$ & $H_{1}^{1}$ & $H_{1}^{2}$ &   $H_{0}$ & $H_{1}^{1}$ & $%
H_{1}^{2}$ \\ \hline
QMTTL-SA\tabnoteref{b1} &0.054 &0.694 &0.995   &0.046 &0.951 &
0.999 &  0.059 &0.664 &0.789   &0.048 & 1.00 & 1.00 \\
QMTTL-WA &0.068 &0.431 &0.924 &0.041 & 1.00 & 1.00
&0.065 &0.067 &0.868 &0.045 & 1.00 & 1.00 \\
QMTTL-S &0.074 &0.256 &0.880 &0.036 & 1.00 & 1.00
&0.058 &0.166 &0.942 &0.040 & 1.00 & 1.00 \\
MNWM-I\tabnoteref{c1} &0.055 &0.521 &0.840 &0.054 &0.899 &
0.997 &0.058 &0.716 &0.927 &0.047 &0.998 & 1.00 \\
MNWM-T &0.043 &0.791 &0.981 &0.055 &0.878 &0.991 &0.031 &0.963 &0.998 &0.054 &0.992 & 1.00 \\
MNWM-E &0.050 &0.236 &0.907 &0.058 &0.867 &0.982
&0.062 &0.573 &0.981 &0.053 &0.988 & 1.00 \\
WLQML\tabnoteref{d1} &0.058 &0.045 &0.838 &0.038 &0.006 &0.260
&0.047 &0.031 &0.809 &0.043 &0.052 &0.793 \\
WLQML$_{E|\epsilon _{t}|=1}$ &0.041 &0.012 &0.493 &
0.041 &0.036 &0.436 &0.038 &0.021 &0.771 &0.044 &0.104 &0.965 \\
PQML$_{3}$ &0.058 &0.367 &0.955 &0.060 &0.688 & 1.00
&0.049 &0.243 &0.980 &0.034 & 1.00 & 1.00 \\
PQML$_{3.5}$ &0.051 &0.578 &0.980 &0.058 &0.891 &
0.998 &0.045 &0.246 &0.983 &0.039 & 1.00 & 1.00 \\
PQMTTL$_{3.5}^{\mathrm{SA}}$ &0.078 &0.122 &0.845 &
0.043 &0.863 & 1.00 &0.058 &0.260 &0.952 &0.043 & 1.00 & 1.00 \\[3pt]
PQMTTL$_{3.5}^{\mathrm{WA}}$ &0.053 &0.600 &0.978 &
0.056 &0.938 & 1.00 &0.060 &0.281 &0.955 &0.051 & 1.00 & 1.00 \\[3pt]
PQMTTL$_{3.5}^{\mathrm{S}}$ &0.053 &0.662 &0.978 &0.055
&0.888 & 1.00 &0.057 &0.385 &0.972 &0.040 & 1.00 & 1.00 \\
Log-LAD &0.061 &0.009 &0.000 &0.025 &0.000 &0.000 &0.058 &0.046 &0.785 &0.053 &0.951 & 1.00 \\
QML &0.065 &0.109 &0.789 &0.061 &0.342 &0.733 &
0.061 &0.000 &0.575 &0.051 &0.997 & 1.00\\
  \hline
\end{tabular*}
\tabnotetext[a]{a1}{The hypotheses are $H_{ 0 }$: $\theta
_{3}=\theta _{3}^{0}$, $H_{ 1} ^{ 1} $: $\theta _{3}=\theta
_{3}^{0}-0.2$, and $H_{ 1 }^{ 2} $: $\theta _{3}=\theta _{3}^{0}
-0.4$, where $\theta _{3}^{0}=0.9$.}
\tabnotetext[b]{b1}{\textit{Benchmark} QMTTL-SA (strong asymmetric) uses
fractiles $k_{1,n} = 35k_{2,n}$;
QMTTL-WA (weak asymmetric) uses $k_{1,n} =
10k_{2,n}$; QMTTL-S (symmetric) uses $k_{1,n} = k_{2,n}$.}
\tabnotetext[c]{c1}{\textit{Benchmark} MNWM-I uses the simple trimming
function $\psi (u,c) = uI(|u| \leq  c)$; MNWM-T and
MNWM-E use Tukey's bisquare and exponential
transforms.}
\tabnotetext[d]{d1}{WLQML is Weighted Laplace QML. WLQML$_{E|\epsilon
_{t}|=1}$ is WLQML for processes with $E|\epsilon _{t}|=1$.
PQML$_{\vartheta }$ is power-law QML with
criterion index $\vartheta $. PQMTTL$_{\vartheta }^{\mathrm{WA}}$ and
PQMTTL$
_{\vartheta }^{\mathrm{SA}}$ are tail-trimmed
PQML with \textit{weak asymmetric} ($k_{1,T} =
 5k_{2,T}$) or \textit{strong asymmetric} ($k_{1,T} = 9k_{2,T}$)
trimming.}%
\end{sidewaystable*}

The PQML estimator detailed in Berkes and Horvath \cite
{BerkesHorvath04}, Example~2.3, is based
on the criterion $-\sum_{t=2}^{n}\ln(h_{t}^{-1/2}(\theta
)f(y_{t}/h_{t}^{1/2}(\theta)))$ where $f(u) = K(1 + %
|u|)^{-\vartheta}$ with tail index $\vartheta > 1$. The value $K > %
0$ ensures $\int_{-\infty}^{\infty}f(u)\,\mathrm{d}u = 1$ and of course is
irrelevant for estimation, hence we simply set $K = 1$. Identification
of $\theta^{0}$ requires $E[|\epsilon_{t}|/(1 + |\epsilon_{t}|)] =
1/\vartheta$, while in the Pareto case $P(|\epsilon_{t}| \geq a) =
(1 + a)^{-\kappa}$ it is easily verified that $E[|\epsilon_{t}|/(1 %
+ |\epsilon_{t}|)] = 1/(\kappa + 1)$ hence we set $\vartheta %
= \kappa + 1 = 3.5$ in both Paretian and Gaussian cases.\footnote{%
Simply note $P(|\epsilon_{t}| \geq a) = (1 + a)^{-\kappa}$
implies $P(|\epsilon_{t}|/(1+|\epsilon_{t}|) > a)=P(|\epsilon_{t}| %
> a/(1 - a)) = (1 - a)^{\kappa}$ hence $E[|\epsilon_{t}|/(1 +
|\epsilon_{t}|)] = \int_{0}^{1}P(|\epsilon_{t}|/(1 +| %
\epsilon_{t}|) > a)\,\mathrm{d}a = \int_{0}^{1}(1 - a)^{\kappa}\,\mathrm{d}a = %
1/(1 + \kappa)$.} We also set $\vartheta = 3$ as a control case to
see if small sample bias increases when $\epsilon_{t}$ is Pareto, as it should.

%s5.2 #&#
\subsection{Simulation results}

Table~\ref{tab1} contains estimator bias, root mse [rmse], and the Kolmogorov--Smirnov
statistic scaled by its 5\% critical value. We only report results for $
\theta_{3}^{0}$ in order to conserve space, while the omitted results are
qualitatively similar. In Table~\ref{tab2}, we report $t$-test rejection
frequencies for
tests of the hypotheses $\theta_{3}^{0} = 0.9$, $\theta_{3}^{0} = %
0.70$ and $\theta_{3}^{0} = 0.50$,\vadjust{\goodbreak} where the first is true. If $\{\hat{%
\theta}_{3,n}^{(r)}\}_{r=1}^{R}$ is the sequence of $R = 10,000$
independent estimates of $\theta_{3}^{0}$, we use the empirical
variance $%
1/R\sum_{r=1}^{R}(\hat{\theta}_{3,n}^{(r)} - 1/R\sum_{r=1}^{R}\hat
{\theta%
}_{3,n}^{(r)})^{2}$ to standardize $\hat{\theta}_{3,n}^{(r)}$ for KS test
and $t$-test computation.

Log-LAD and WLQML perform poorly when $E[\epsilon_{t}^{4}] = \infty$:
in small samples they are sensitive to large error observations,
contrary to
their theoretical robustness properties asymptotically. Indeed, Log-LAD
leads to exceptionally poor inference when $E[\epsilon_{t}^{4}] = %
\infty$ due to a high degree of bias, and is worst overall. Further, WLQML
is sensitive to large errors even in the Gaussian case. It is not surprising
that Log-LAD and WLQML are similar since Laplace QML merely generalizes LAD
to a likelihood framework (Zhu and Ling \cite{ZhuLing}). QML performs
better than Log-LAD
and worse than WLQML when $E[\epsilon_{t}^{4}] = \infty$, and is
better than both when $\epsilon_{t}$ is normal.

PQML is more promising than QML, Log-LAD and WLQML. It performs better on
all measures and in nearly every case: Log-LAD and WLQML are closer to
normally distributed for Gaussian $\epsilon_{t}$ with small $n = 100$.
In particular, PQML has the smallest rmse of all estimators in this study,
suggesting that it exhibits very low empirical variance since it has higher
bias than QMTTL and MNWM. Identification is assured in the Pareto case $
\kappa = 2.5$ when $\vartheta = 3.5$, so it is not surprising that
bias in the Pareto case is higher when $\vartheta = 3$. Further, there
should be noticeable bias in the Gaussian case since identification fails,
yet bias is actually \textit{smaller} than for Paretian errors when $n
= %
800$. It is important to stress that PQML with index $\vartheta = 3.5$
is perfectly suited for our Paretian case $P(|\epsilon_{t}| \geq a) %
= (1 + a)^{-2.5}$ since this non-Gaussian QML leads to identification
and therefore Fischer consistency. However, even this estimator exhibits
more bias than QMTTL and MNWM evidently due to the adverse effects of sample
error extremes (see Section~\ref{sec:add_exp}).

The best estimators in this study are QMTTL (with strong asymmetric
trimming) and MNWM in terms of bias, approximate normality and test
performance, while only PQML has a smaller rmse. QMTTL with strong
asymmetric trimming ($k_{1,n} = 35k_{2,n}$), as required in the Paretian
case when $n = 100$, is superb when $\epsilon_{t}$ is Paretian for
either $n \in \{100,800\}$, and works very well in the Gaussian case
with a rmse close to PQML. Overall, \textit{QMTTL with strong asymmetric
trimming is the best estimator} since it beats MNWM in terms of bias and
approximate normality in nearly every case and has a small rmse in all cases.

QMTTL with weak asymmetric ($k_{1,n} = 10k_{2,n}$) or symmetric ($%
k_{1,n} = k_{2,n}$) trimming lead to greater bias when $\epsilon_{t}$
is Paretian, and to negligible bias when $\epsilon_{t}$ is Gaussian, in
each case as this estimator should. Nevertheless, QMTTL with weak asymmetric
or symmetric trimming is superior to QML, Log-LAD, and WLQML by all measure;
QMTTL with weak asymmetric trimming beats PQML by all measures except rmse;
and QMTTL with symmetric trimming beats PQML when $n = 800$. Our QMTTL
simulations strongly point to the use of strong asymmetric trimming in
general since it is valid for thin tailed errors, and necessary for heavy
tailed errors. They also reveal that using weak asymmetric of symmetric
trimming still leads to a competitive estimator.

Further, re-centering after trimming in the MNWM estimator in general leads
to higher mean-squared-error than QMTTL. Recall this estimator may be less
efficient than QMTTL, and QMTTL with strong asymmetric trimming results in
the lowest bias of all estimators in this study. Nevertheless, MNWM works
well, with the second smallest bias, and overall is closer to normal than
all estimators save QMTTL with strong asymmetric trimming. As discussed in
Section~\ref{sec:bias_remarks}, the preferred estimator depends on the
analyst's agenda: MNWM is always asymptotically unbiased with symmetric
trimming which is easy to implement, while QMTTL performs better in small
samples.\vspace*{-1pt}

%s5.3 #&#
\subsection{Addtional experiments for WLQML and PQML}
\label{sec:add_exp}

We now perform two additional experiments. First, recall WLQML requires
$%
E|\epsilon_{t}| = 1$ which does not hold for either Paretian or
Gaussian errors in this study. We now standardize $\epsilon_{t}$ such
that $%
E|\epsilon_{t}| = 1$ to see if ensuring identification helps in small
samples. The results are nevertheless qualitatively similar whether $%
E[\epsilon_{t}^{2}] = 1$ and $E|\epsilon_{t}| \neq 1$, or $%
E|\epsilon_{t}|= 1$, is true. See Tables~\ref{tab1} and \ref{tab2}. In fact, for heavy
tailed errors WLQML actually performs \textit{worse} in terms of bias and
approximate normality when identification is assured. Further,
inference is
still quite poor in many cases. This suggests the previous poor performance
of WLQML is not due to the identification condition failing to hold.

Second, recall that QMTTL has lower bias and is closer to normally
distributed that other estimators whether trimming is needed or not. We
therefore tail trim the PQML criterion to see if the benefits of trimming
carry over to non-Gaussian QML. Recall PQML with index $\vartheta > 1$
has the identification condition $E[u_{t}] = 0$ where $u_{t} := %
|\epsilon_{t}|/(1 + |\epsilon_{t}|)-1/\vartheta$. Define $%
u_{t}^{(-)}(\theta) := u_{t}(\theta)I(u_{t}(\theta) < 0)$ and $%
u_{t}^{(+)}(\theta) := u_{t}(\theta)I(u_{t}(\theta) \geq 0)$ and
their order statistics $u_{(1)}^{(-)}(\theta) \leq \cdots \leq %
u_{(n)}^{(-)}(\theta) \leq 0$ and $u_{(1)}^{(+)}(\theta) \geq %
\cdots \geq u_{(n)}^{(+)}(\theta) \geq 0$. Let $%
\{k_{1,n}^{(u)},k_{2,n}^{(u)}\}$ be intermediate order sequences and
let $%
\{c_{1,n}^{(u)},c_{2,n}^{(u)}\}$ be positive sequences satisfying $%
P(u_{t}(\theta) \leq -c_{1,n}^{(u)}) = k_{1,n}/n$ and $%
P(u_{t}(\theta) \geq c_{2,n}^{(u)}) = k_{2,n}/n$. The tail-trimmed
PQML (PQMTTL) criterion is $-\sum_{t=2}^{n}\ln(h_{t}^{-1/2}(\theta)\{
1 + |y_{t}/h_{t}^{1/2}(\theta)|\}^{-(\kappa
+1)})I(u_{(k_{1,n}^{(u)})}^{(-)}(\theta) \leq u_{t}(\theta) \leq %
u_{(k_{2,n}^{(u)})}^{(+)}(\theta))$.

If $\epsilon_{t}$ is Paretian $P(|\epsilon_{t}| \geq a) = (1 +
a)^{-2.5}$ it is straightforward to show $k_{1,n}^{(u)} = %
5k_{2,n}^{(u)}$ when $n = 100$ and $k_{1,n}^{(u)} = 9k_{2,n}^{(u)}$
when $n = 800$ renders roughly $E[u_{t}I(-c_{1,n}^{(u)} \leq %
u_{t}(\theta) \leq c_{2,n}^{(u)})] = 0$. We therefore set \textit{%
symmetric} ($k_{1,n}^{(u)} = k_{2,n}^{(u)}$), \textit{weak asymmetric}
($%
k_{1,n}^{(u)} = 5k_{2,n}^{(u)}$) or \textit{strong asymmetric} ($%
k_{1,n}^{(u)} = 9k_{2,n}^{(u)}$) trimming with $k_{2,n}^{(u)} = %
\max\{1,[0.025n/\ln(n)]\}$. Tables~\ref{tab1} and \ref{tab2} show PQMTTL with weak asymmetric
trimming performs better than PQML in all cases. If we use strong asymmetric
trimming then the over-trimming for $n = 100$ leads to greater bias, but
when $n = 800$ the estimator works well as it should, in particular it
is closer to normal and therefore has better inference than PQML.
Conversely, symmetric trimming leads to greater bias when $n = 800$ as
it should. QMTTL with strong asymmetric trimming and MNWM with simple or
exponential trimming are better than PQMTTL in terms of bias and approximate
normality in most cases. Consider when $n = 800$ then in the Pareto
case PQMTTL with weak asymmetric trimming is marginally closer to
normal and
slightly more biased than QMTTL, and in the Gaussian case PQMTTL is slightly
less biased and farther from normally distributed than QMTTL. Overall
tail-trimming seems to matter even for an inherently heavy tail robust
non-Gaussian QML estimator.

%s6 #&#
\section{Empirical application}\label{s6}

Finally, we apply our estimators to asset returns series generated from the
London Stock Exchange (FTSE-100),\vadjust{\goodbreak} the NASDAQ composite index (IXIC),
and the
Hang Seng Index. The period is Jan. 1, 2008--Dec. 31, 2010, representing
757, 757 and 756 daily observations respectively, net of market
closures. We
use log-returns $y_{t} = \ln(x_{t}/x_{t-1})$ where $x_{t}$ is the daily
open/close average of each index.\footnote{%
The data were obtained from \url{http://finance.yahoo.com}, and the open/close
average is
computed using the reported adjusted close values.}

As in Section~\ref{sim}, we compute MNWM using simple trimming denoted ``I'',
Tukey's bisquare and exponential transforms, with fractiles $k_{n} = %
\max\{1,[0.025n/\ln(n)]\}$ and $\tilde{k}_{n} = \max\{1,[0.1\ln(n)]\}$
for trimming by $\epsilon_{t}$ and $y_{t-1}$, respectively.
Similarly, QMTTL
is computed using \textit{strong asymmetric} ($k_{1,n} = 35k_{2,n}$),
\textit{weak asymmetric} ($k_{1,n} = 10k_{2,n}$), and \textit{symmetric} ($k_{1,n} = k_{2,n}$) error fractiles denoted ``SA'', ``WA'' and
``S'', with $k_{2,n} = \max\{1,[0.025n/\ln(n)]\}$, and $\tilde{k}_{n}$
for $y_{t-1}$. The parameter space is $\Theta = [\iota,2] \times %
[\iota,1-\iota] \times [\iota,1-\iota]$ where $\iota = %
10^{-10} $.

\begin{sidewaystable*}
\tablewidth=\textwidth
\caption{QMTTL and MNWM estimation results for
financial returns}\label{tab3}
\begin{tabular*}{\textwidth}{@{\extracolsep{\fill}}lllllll@{}}
\hline
 & $\omega $ & $\alpha $ & $\beta $ &   $\omega $ & $\alpha $ & $\beta $ \\
&\multicolumn{3}{l}{\hrulefill} &
\multicolumn{3}{l@{}}{\hrulefill} \\
 & \multicolumn{3}{l}{QMTTL-SA\tabnoteref{a2}: $k_{1,n}=35k_{2,n}$} &
\multicolumn{3}{l@{}}{MNWM-I\tabnoteref{b2}} \\
\hline
 NASDAQ\tabnoteref{c2}  & 0.029 (0.031)\tabnoteref{c2} &
0.113 (0.082) & 0.893 (0.069) &    0.016 (0.008) & 0.117 (0.017) & 0.884 (0.018) \\
 HSI\tabnoteref{d2}  & 0.058 (0.064) & 0.106 (0.151) & 0.878
(0.252) &    0.020 (0.009) & 0.078 (0.013) & 0.915 (0.015) \\
 LSE  & 0.066 (0.083) & 0.213 (0.156) & 0.743 (0.224) &    0.025
(0.006) & 0.119 (0.015) & 0.822 (0.020) \\[5pt]
& \multicolumn{3}{l}{QMTTL-WA: $k_{1,n}=10k_{2,n}$} &
\multicolumn{3}{l@{}}{MNWM-E} \\
 NASDAQ  & 0.017 (0.038) & 0.138 (0.113) & 0.849 (0.135)    &
0.032 (0.010) & 0.102 (0.021) & 0.886 (0.019) \\
 HSI  & 0.046 (0.086) & 0.082 (0.142) & 0.910 (0.211) &   0.021
(0.011) & 0.078 (0.016) & 0.915 (0.028) \\
 LSE  & 0.022 (0.065) & 0.179 (0.134) & 0.805 (0.194) &    0.030
(0.011) & 0.125 (0.019) & 0.824 (0.031) \\[5pt]
& \multicolumn{3}{l}{QMTTL-S: $k_{1,n}=k_{2,n}$} &    \multicolumn{3}{l@{}}{
MNWM-T} \\
  NASDAQ  & 0.012 (0.033) & 0.171 (0.145) & 0.839 (0.125)    &
0.033 (0.011) & 0.095 (0.031) & 0.901 (0.034) \\
 HSI  & 0.065 (0.092) & 0.092 (0.123) & 0.887 (0.163)    & 0.039
(0.011) & 0.076 (0.017) & 0.920 (0.027) \\
 LSE  & 0.034 (0.076) & 0.214 (0.189) & 0.752 (0.203) &    0.058
(0.016) & 0.122 (0.021) & 0.839 (0.024) \\ \hline
\end{tabular*}
\tabnotetext[a]{a2}{SA~$=$ strong asymmetry; WA~$=$ weak asymmetry; S~$=$
symmetric.}
\tabnotetext[b]{b2}{I~$=$ simple trimming, T~$=$ Tukey's bisquare, E~$=$
exponential.}
\tabnotetext[c]{c2}{Standard errors are in parentheses ($\cdot $).}
\tabnotetext[d]{d2}{HSI~$=$ Hang Seng; LSE~$=$ London Stock Exchange.}%
\end{sidewaystable*}

See Table~\ref{tab3} for estimation details where standard errors are computed
using (%
\ref{V_hat}) for QMTTL and its logical extension for MNWM. In each case a
GARCH model fits well, while QMTTL and MNWM produce qualitatively similar
estimates. The various MNWM estimates are similar across transform type,
especially exponential and simple trimming versions. The QMTTL
estimates are
somewhat similar across asymmetric and symmetric trimming. For example,
evidence for IGARCH or explosive GARCH $\hat{\alpha}_{n} + \hat{\beta}%
_{n} \geq 1$ exists only for the NASDAQ based on QMTTL-SA and MNWM-I,
while QMTTL-WA and QMTTL-S lead to smaller values. However, in all
cases $%
\hat{\beta}_{n}$ is near $0.9$ and $\hat{\alpha}_{n}$ is near $0.05$, in many
cases $\hat{\alpha}_{n} + \hat{\beta}_{n} \approx 1$, and for each
series the various estimates are quite similar. The latter suggests the
various asymmetric and symmetric trimming strategies for QMTTL work as well
as inherently asymptotically unbiased MNWM. This is matched by our
simulations where $n = 800$ aligns with the sample sizes in the present
empirical study: strong asymmetric trimming leads to the best QMTTL results
when $\epsilon_{t}$ has power law tails with a small index $\kappa$, but
each trimming strategy leads to similar results, especially when $n = %
800 $.

%s7 #&#
\section{Conclusion}

We develop tail-trimmed QML and Method of Moments estimators for GARCH
models with possibly heavy tailed errors $\epsilon_{t}$ that satisfy $%
E[\epsilon_{t}^{2}] = 1$. In the Method of Moments case, the model
errors are first negligibly transformed with a redescending function, and
then re-centered to control for small sample bias induced by the transform.
We show by Monte Carlo experiment that tail-trimming within a QML framework
dominates QML, Log-LAD and Weighted Laplace QML based on bias,
mean-squared-error, approximate normality, and inference, and trumps
Power-Law QML in all aspects except variance (Power-Law QML has higher bias
yet lower mean-squared-error). Only QMTTL and MNWM directly counter the
negative influence of large errors in small \textit{and} large samples.
Indeed, we show trimming leads to a better infeasible Power-Law QML
estimator in small samples. The next stage must involve a theoretical
development of data-dependent or automatic fractile selection, including
possibly bootstrap and covariance determinant methods. This is left for
future research.

%\section*{Acknowledgements}

%\setcounter{equation}{0} \renewcommand{\theequation}{{\thesection}.%
%\arabic{equation}} \appendix
\begin{appendix}

%s8 #&#
\section*{Appendix: Proofs of main theorems}\label{app}
%\label{app:proofs}
\setcounter{equation}{0}

Recall $\mathcal{E}_{t}(\theta) := \epsilon_{t}^{2}(\theta) - 1$,
and score and Jacobian equations are $m_{t}(\theta) = \mathcal{E}%
_{t}(\theta)\mathfrak{s}_{t}(\theta)$ and $G_{t}(\theta) = \mathcal
{E}%
_{t}(\theta)\mathfrak{d}_{t}(\theta) - \epsilon_{t}^{2}(\theta)%
\mathfrak{s}_{t}(\theta)\mathfrak{s}_{t}(\theta)^{\prime}$ where
\[
\mathfrak{s}_{t}(\theta):=\frac{1}{\sigma_{t}^{2}(\theta)}\frac
{\partial
}{\partial\theta}
\sigma_{t}^{2} ( \theta ) \quad\mbox{and}\quad%
\mathfrak{d}_{t}(\theta):=\frac{\partial}{\partial\theta}\mathfrak {s}%
_{t}(\theta).
\]
Define indicators, trimmed score equations and corresponding covariance and
Jacobian matrices:
\begin{eqnarray*}
 \hat{m}_{n,t}(\theta)&:=&m_{t}(\theta)\hat{I}_{n,t}^{(\mathcal
{E})}(
\theta )\quad\mbox{and}\quad m_{n,t}(\theta):=m_{t}(
\theta)I_{n,t}^{(\mathcal{E}%
)}(\theta),
\\
 \Sigma_{n}(\theta)&:=&E \bigl[ m_{n,t}(\theta)m_{n,t}(
\theta )^{\prime}%
\bigr] \quad\mbox{and}
\\
\mathcal{G}(\theta)&:=&-E
\bigl[ \mathfrak{s}%
_{t}(\theta)\mathfrak{s}_{t}(
\theta)^{\prime} \bigr] \quad\mbox{and} \quad
\mathcal{V}_{n}(
\theta)=n\mathcal{G}(\theta)^{\prime}\Sigma _{n}^{-1}(
\theta)\mathcal{G}(\theta),
\\
 \mathcal{S}_{n}(\theta)&:=&\frac{1}{n}E \Biggl[ \Biggl( \sum
_{t=1}^{n}m_{n,t}(\theta) \Biggr)
\Biggl( \sum_{t=1}^{n}m_{n,t}(
\theta )^{\prime} \Biggr) \Biggr],
\\
  \widehat{\mathcal{G}}_{n}(\theta)&:=&\frac{1}{n}%
\sum_{t=1}^{n}G_{t}(\theta)
\hat{I}_{n,t}^{(\mathcal{E})}(\theta)\quad\mbox {and}\quad\mathcal{
\check{G}}_{n}(\theta):=\frac{1}{n}\sum
_{t=1}^{n}G_{t}( %%
\theta)I_{n,t}^{(\mathcal{E})}(\theta).
\end{eqnarray*}
By independence and identification Assumption~\ref{ass2} $\mathcal{S}_{n} = %
\Sigma_{n} \times (1 + \mathrm{o} ( 1 ) )$.

We implicitly assume all functions in this paper satisfy Pollard's
(\cite{pollard84},
Appendix C) permissibility criteria, the measure space
that governs all random variables in this paper is complete, and therefore
all majorants are measurable. Cf. Dudley \cite{Dudley78}. Probability statements
are therefore with respect to outer probability, and expectations over
majorants are outer expectations.

%s8.1 #&#
\subsection{Theorems \texorpdfstring{\protect\ref{th:consist}}{2.1} and \texorpdfstring{\protect\ref{th:norm}}{2.2}}
\label{app:th21_22}

The proofs of QMTTL consistency and asymptotic normality Theorems \ref%
{th:consist} and \ref{th:norm} require supporting lemmas. We state them when
required and provide proofs in Appendix~\ref{app:lemmas}. Consistency
requires bounding $\sum_{t=1}^{n}\{\hat{m}_{n,t}(\theta) - %
m_{n,t}(\theta)\}$, variance bounds, and laws of large numbers. Unless
otherwise noted, Assumptions \ref{ass1} and \ref{ass2} hold.

%le8.1 #&#
%
\begin{lemma}[(Asymptotic approximation)]\label{lm:approx}
\textup{(a)}  $n^{-1/2}\Sigma_{n}^{-1/2}\sum_{t=1}^{n}\{\hat{m}%
_{n,t} - m_{n,t}\} =
\mathrm{o}_{p}(1) $;  \textup{(b)}~$\sup_{\theta\in\Theta}\|1/n\sum_{t=1}^{n}\{%
\hat{m}_{n,t}(\theta) - m_{n,t}(\theta)\}\| = \mathrm{o}_{p}(\sup_{\theta
\in\Theta}E\|m_{n,t}(\theta)\|)$.
\end{lemma}

%le8.2 #&#
%
\begin{lemma}[(Variance bounds)]\label{lm_var_bnds}
Under Assumption~\ref{ass1}\textup{(a)}   $\Sigma_{n} = \mathrm{o}(n/\ln(n))$;
\textup{(b)}  $\mathcal{S}_{n} = \mathrm{o}(n/\ln(n))$.
\end{lemma}

%re11 #&#
%
\begin{remark}
Under Assumption~\ref{ass2} $\mathcal{S}_{n} = \Sigma_{n}(1
+ \mathrm{o}(1))$
hence then (b) follows from (a).
\end{remark}

%le8.3 #&#
%
\begin{lemma}[(LLN and ULLN)]\label{lm:ulln}
\textup{(a)}  $1/n\sum_{t=1}^{n}m_{n,t} = \mathrm{o}_{p}(1)$; \textup{(b)}  $
\sup_{\theta\in\Theta}\{\|1/n\times\break \sum_{t=1}^{n}m_{n,t}(\theta) -
E[m_{n,t}(\theta)]\|\} =
\mathrm{o}_{p}(\sup_{\theta\in\Theta}E\|m_{n,t}(\theta)\|)$.\vadjust{\goodbreak}
\end{lemma}

Asymptotic normality requires an expansion, central limit theorem, and
Jacobian consistency.

%le8.4 #&#
%
\begin{lemma}[(Asymptotic expansion)]\label{lm:expan}
Let $\{\theta_{n}\}$ and $\{\tilde{\theta}%
_{n}\}$ be any sequences of random variables in $\Theta$
with probability limit $\theta^{0}$. Let $\theta_{n,\ast} \in
\Theta$ satisfy $\|\theta_{n,\ast} - \theta_{n}\| %
\leq \|\theta_{n} - \tilde{\theta}_{n}\|$ which may be
different in difference places.   \textup{(a)}  $1/n\sum_{t=1}^{n}\{m_{n,t}(\theta
_{n}) - m_{n,t}(\tilde{\theta}_{n})\} = \mathcal{\check{G}}%
_{n}(\theta_{n,\ast}) \times (\theta_{n} - \tilde{\theta}_{n})
\times (1 + %
\mathrm{o}_{p}(1))4$;   \textup{(b)} $1/n\sum_{t=1}^{n}\{\hat{m}_{n,t}(\theta_{n}) - \hat
{m}_{n,t}(\tilde{\theta}_{n})\} = \widehat{\mathcal{G}}_{n}(\theta
_{n,\ast}) \times (\theta_{n} - \tilde{\theta}_{n}) \times (1 + \mathrm{o}_{p}(1))$.
\end{lemma}

%le8.5 #&#
%
\begin{lemma}[(CLT)]\label{lm:clt}
$n^{-1/2}\Sigma_{n}^{-1/2}\sum_{t=1}^{n}m_{n,t} \stackrel
{d}{%
\rightarrow} N(0,I_{q})$.
\end{lemma}

%le8.6 #&#
%
\begin{lemma}[(Jacobian)]\label{lm:jac}
\textup{(a)}  $\widehat{\mathcal{G}}_{n}(\hat{\theta}_{n}^{\ast}) =
\mathcal{G} \times (1 + \mathrm{o}_{p}(1))$ and $\mathcal{\check{G}}_{n}%
\hat{\theta}_{n}^{\ast}) = \mathcal{G} \times (1 + \mathrm{o}_{p}(1))$
for any $\hat{\theta}_{n}^{\ast} \stackrel{p}{\rightarrow} %
\theta^{0}$;  \textup{(b)} $1/n\sum_{t=1}^{n}\mathfrak{s}_{t}(\hat{\theta}_{n})%
\mathfrak{s}_{t}^{\prime}(\hat{\theta}_{n}) = \mathcal{G} \times %
(1 + \mathrm{o}_{p}(1))$; \textup{(c)}  $(\partial/\partial\theta)E[m_{n,t} (
\theta ) ]|_{\theta^{0}} = \mathcal{G} \times (1 + \mathrm{o}(1))$; \textup{(d)}
$\lim\sup_{n\rightarrow
\infty}\sup_{\theta\in\Theta}E\|m_{n,t}(\theta)\| \leq K\|\mathcal
{%
G}\|\times(1 + \mathrm{o}(1))$.
\end{lemma}

We are now ready to prove Theorems \ref{th:consist} and \ref{th:norm}.

\begin{pf*}{Proof of Theorem~\ref{th:consist}} Define $\hat
{m}_{n}(\theta
) := 1/n\sum_{t=1}^{n}\hat{m}_{n,t}(\theta)$, $m_{n}(\theta) := %
1/n\sum_{t=1}^{n}m_{n,t}(\theta)$, $\mathcal{M}_{n}(\theta) := %
E[m_{n,t}(\theta)]$ and $\mathfrak{e}_{n} := \sup_{\theta\in\Theta
}E\|m_{n,t}(\theta)\|$. We use an argument in Pakes and Pollard \cite{pakespollard89}, pages~1038--1039.

\textit{Step 1.} We first prove a required inequality:%
%e1 #&#
\begin{equation}\label{ed}
\epsilon(\delta):=\liminf_{n\rightarrow\infty}\inf_{\theta\in
\Theta
:\|\theta-\theta^{0}\|>\delta}
\bigl\{ \bigl\llVert \mathcal {M}_{n}(\theta )\bigr\rrVert /
\mathfrak{e}_{n} \bigr\} >0\quad\quad\mbox{for any small }\delta>0.
\end{equation}
Note $E[m_{n,t}] \rightarrow E[\mathcal{E}_{t}\mathfrak{s}_{t}] = 0$
by dominated convergence and independence. By the definition of a derivative
and Lemma~\ref{lm:jac}(c) we have $E[m_{n,t} ( \theta ) ] = %
\mathcal{G} \times (\theta - \theta^{0}) \times (1 + %
\mathrm{o}(1))$ where $\mathcal{G} = -E[\mathfrak{s}_{t}\mathfrak{s}_{t}^{\prime
}]$, and bound Lemma~\ref{lm:jac}(d) states $\mathfrak{e}_{n} := %
\sup_{\theta\in\Theta}E\|m_{n,t} ( \theta ) ] \leq K\|%
\mathcal{G}\|\times(1 + \mathrm{o}(1))$. It~therefore follows for every $n %
\geq N$ and $\delta > 0$
\[
\inf_{\|\theta-\theta^{0}\|>\delta} \bigl\{ \mathfrak{e}%
_{n}^{-1}\bigl\llVert E \bigl[ m_{n,t}(\theta)
\bigr] \bigr\rrVert \bigr\} \geq K\inf_{\|\theta-\theta^{0}\|>\delta} \biggl\{ \biggl
\llVert \frac{%
\mathcal{G}}{\llVert \mathcal{G}\rrVert }\times \bigl( \theta -\theta ^{0} \bigr)
\biggr\rrVert \biggr\} \times \bigl( 1+\mathrm{o} ( 1 ) \bigr) >0.
\]
\textit{Step 2.} In view of (\ref{ed}) we have $P(\|\hat{\theta}_{n} %
- \theta^{0}\| > \delta) \leq P(\|\mathcal{M}_{n}(\hat{\theta}%
_{n})\|/\mathfrak{e}_{n} > \epsilon(\delta))$, hence it suffices to
show $\|\mathcal{M}_{n}(\hat{\theta}_{n})\|/\mathfrak{e}_{n}= \mathrm{o}_{p}(1)$ in
order to prove $\hat{\theta}_{n} \stackrel{p}{\rightarrow} \theta^{0}$.
By Minkowski's inequality
\[
\bigl\llVert \mathcal{M}_{n}(\hat{\theta}_{n})\bigr\rrVert /
\mathfrak{e}%
_{n}\leq\bigl\llVert \hat{m}_{n}(\hat{
\theta}_{n})\bigr\rrVert /\mathfrak {e}%
_{n}+\bigl
\llVert \hat{m}_{n}(\hat{\theta}_{n})-\mathcal{M}_{n}(
\hat {\theta}%
_{n})\bigr\rrVert /\mathfrak{e}_{n}=
\mathcal{A}_{1,n}(\hat{\theta}_{n})+%
\mathcal{A}_{2,n}(\hat{\theta}_{n}),
\]
say. The proof is complete if we show $\mathcal{A}_{1,n}(\hat{\theta}_{n})$
and $\mathcal{A}_{2,n}(\hat{\theta}_{n})$ are $\mathrm{o}_{p}(1)$.

Consider $\mathcal{A}_{1,n}(\hat{\theta}_{n})$. We exploit theory developed
in Cizek \cite{cizek08}, Lemma~2.1, page 29. By distribution
continuity and linearity of the volatility process $\{\sigma_{t}^{2}\}
$, $%
\hat{Q}_{n}(\theta) := 1/n\sum_{t=1}^{n}(\ln\sigma_{t}^{2}(\theta) %
+ y_{t}^{2}/\sigma_{t}^{2}(\theta))\hat{I}_{n,t}^{(\mathcal
{E})}(\theta
) $ is almost surely twice differentiable at $\hat{\theta}_{n}$. In
particular, up to a scalar constant $(\partial/\partial\theta)\hat{Q}
_{n}(\theta)|_{\hat{\theta}_{n}} = \hat{m}_{n}(\hat{\theta}_{n})$ a.s.
By $\hat{\theta}_{n}$ a minimum $\hat{Q}_{n}(\hat{\theta}_{n}) \leq %
\hat{Q}_{n}(\theta) \forall\theta \in \Theta$ it follows $\|\hat{m%
}_{n}(\hat{\theta}_{n})\| = 0$ a.s., while $\lim\inf_{n\rightarrow
\infty}\mathfrak{e}_{n} > 0$ by distribution non-degeneracy and
trimming negligibility, hence $\mathcal{A}_{1,n}(\hat{\theta}_{n}) = 0$
a.s.

Next $\mathcal{A}_{2,n}(\hat{\theta}_{n})$. By Lemma~\ref{lm:approx}(b) $
\sup_{\theta\in\Theta}\|\hat{m}_{n}(\theta) - m_{n}(\theta)\|/%
\mathfrak{e}_{n} = \mathrm{o}_{p}(1)$, and\break  $\sup_{\theta\in\Theta
}\|m_{n}(\theta) - \mathcal{M}_{n}(\theta)\|/\mathfrak{e}_{n} =
\mathrm{o}_{p}(1)$ by ULLN Lemma~\ref{lm:ulln}(b). Hence\vspace*{1pt}
\[
\sup_{\theta\in\Theta} \bigl\{ \mathcal{A}_{2,n}(\theta) \bigr
\} \leq \sup_{\theta\in\Theta}\frac{\llVert \hat{m}_{n}(\theta
)-m_{n}(\theta
)\rrVert }{\mathfrak{e}_{n}}+\sup_{\theta\in\Theta}
\frac{\llVert
m_{n}(\theta)-\mathcal{M}_{n}(\theta)\rrVert }{\mathfrak{e}_{n}}=%
\mathrm{o}_{p}(1).
\]
\upqed\end{pf*}

\begin{pf*}{Proof of Theorem~\ref{th:norm}}
Use $1/n\sum_{t=1}^{n}\hat{m}%
_{n,t}(\hat{\theta}_{n}) = 0$  a.s.  by the proof of Theorem~2.1,
and expansion Lemma~\ref{lm:expan}(b) to deduce for some $\hat{\theta}%
_{n}^{\ast}$, $\|\hat{\theta}_{n}^{\ast} - \theta^{0}\| \leq \|%
\hat{\theta}_{n} - \theta^{0}\|$:\vspace*{-1pt}
%e2 #&#
\begin{equation}
\widehat{\mathcal{G}}_{n}\bigl(\hat{\theta}_{n}^{\ast}
\bigr) \bigl( \hat{\theta}%
_{n}-\theta^{0} \bigr)
\bigl( 1+\mathrm{o}_{p} ( 1 ) \bigr) +\frac
{1}{n}%
\sum
_{t=1}^{n}\hat{m}_{n,t}=0 \quad\quad \mbox{a.s.}
\label{foc_expand}
\end{equation}
Consistency $\|\hat{\theta}_{n}^{\ast} - \theta^{0}\| \leq \|\hat{%
\theta}_{n} - \theta^{0}\| \stackrel{p}{\rightarrow} 0$ by Theorem~2.1 ensures $\widehat{\mathcal{G}}_{n}(\hat{\theta}_{n}^{\ast}) = %
\mathcal{G}(1 + \mathrm{o}_{p}(1))$ by Lemma~\ref{lm:jac}(a). Multiply both sides
of (\ref{foc_expand}) by $n^{1/2}\Sigma_{n}^{-1/2}$, rearrange terms and
use $\mathcal{V}_{n} = n\mathcal{G}^{\prime}\Sigma_{n}^{-1}\mathcal{G}$
to deduce $\mathcal{V}_{n}^{1/2}(\hat{\theta}_{n}-\theta^{0}) = %
-n^{-1/2}\Sigma_{n}^{-1/2}\sum_{t=1}^{n}\hat{m}_{n,t} \times (1 + %
\mathrm{o}_{p}(1))$. In view of $n^{-1/2}\Sigma_{n}^{-1/2}\sum_{t=1}^{n}\{\hat
{m}%
_{n,t}-m_{n,t}\} = \mathrm{o}_{p}(1)$ by Lemma~\ref{lm:approx}(a), we have\vspace*{-1pt}
\[
\mathcal{V}_{n}^{1/2} \bigl( \hat{\theta}_{n}-
\theta^{0} \bigr) =-\Sigma _{n}^{-1/2}
\frac{1}{n^{1/2}}\sum_{t=1}^{n}m_{n,t}
\times \bigl( 1+\mathrm{o}_{p} ( 1 ) \bigr) ,
\]
hence $\mathcal{V}_{n}^{1/2}(\hat{\theta}_{n} - \theta^{0}) \stackrel{d
}{\rightarrow} N(0,I_{q})$ by Lemma~\ref{lm:clt}. Finally $\mathcal{V}%
_{i,i,n} \rightarrow \infty$ follows from the fact that $\|\mathcal{G}
\| > 0$, and $\|n\Sigma_{n}^{-1}\| \rightarrow \infty$ by
Lemma~\ref{lm_var_bnds}(a).\vspace*{-2pt}
\end{pf*}

%s8.2 #&#
\subsection{Remaining theorems}\label{app:remain}

Define $h_{t}^{\theta}(\theta) := (\partial/\partial\theta
)h_{t}(\theta)$ and $h_{t}^{\theta,\theta}(\theta) := (\partial
/\partial\theta)h_{t}^{\theta}(\theta)$. We require stationary solutions
$\{h_{t}^{\ast}(\theta), h_{i,t}^{\ast\theta}(\theta), %
h_{i,t}^{\ast\theta,\theta}(\theta)\}$ of the volatility process $%
\{h_{t}(\theta), h_{i,t}^{\theta}(\theta), h_{i,t}^{\theta,\theta
}(\theta)\}$\vspace*{1pt} in order to prove the asymptotic equivalence of the infeasible
and feasible QMTTL estimators.

Let $\{\mathfrak{s}_{t}^{\ast}(\theta),\mathfrak{d}_{t}^{\ast
}(\theta)\}$
denote $\{\mathfrak{s}_{t}(\theta),\mathfrak{d}_{t}(\theta)\}$ evaluated
with $\{h_{t}^{\ast}(\theta),h_{i,t}^{\ast\theta}(\theta
),h_{i,t}^{\ast
\theta,\theta}(\theta)\}$. Define error and volatility derivatives
evaluated at $\{h_{t}(\theta),h_{t}^{\theta}(\theta),h_{t}^{\theta
,\theta}(\theta)\}$\vspace*{-1pt}
\begin{eqnarray*}
 \tilde{\epsilon}_{t}(\theta)&:=&\frac{y_{t}}{\sqrt{h_{t}(\theta
)}},\quad\quad\mathcal{
\tilde{E}}_{t}(\theta):=\tilde{\epsilon}_{t}^{2}(
\theta )-1,
\\
\mathfrak{\tilde{s}}_{t}(\theta)&:=&
\frac{1}{h_{t}(\theta)}\frac{%
\partial}{\partial\theta}h_{t} ( \theta ) \quad\mbox{and}\quad%
\mathfrak{\tilde{d}}_{t}(\theta)=\frac{\partial}{\partial\theta}%
\mathfrak{\tilde{s}}_{t}(\theta),
\\
 \tilde{m}_{t}(\theta)&:=&\mathcal{\tilde{E}}_{t}(\theta)
\mathfrak {\tilde{s}%
}_{t}(\theta),\quad\quad \tilde{G}_{t}(
\theta):=\frac{\partial}{\partial
\theta}\tilde{m}_{t}(\theta)\quad\mbox{and}\quad\mathcal{
\tilde{G}}:=-E \bigl[ \mathfrak{\tilde{s}}_{t}(\theta)\mathfrak{
\tilde{s}}_{t}^{\prime
}(\theta)%
\bigr] .
\end{eqnarray*}
Define $\widehat{\tilde{I}}_{n,t}^{(\mathcal{E})}(\theta) := I(\mathcal
{%
\tilde{E}}_{(k_{1,n})}^{(-)}(\theta) \leq \mathcal{\tilde{E}}%
_{t}(\theta) \leq \mathcal{\tilde{E}}_{(k_{2,n})}^{(+)}(\theta))$ and
let $\{\mathcal{\tilde{L}}_{n}(\theta),\widetilde{\mathcal
{U}}_{n}(\theta
)\}$ satisfy $P(\mathcal{\tilde{E}}_{t}(\theta) \leq -\mathcal{\tilde
{L%
}}_{n}(\theta)) = k_{1,n}/n$ and $P(\mathcal{\tilde{E}}_{t}(\theta) %
\geq \widetilde{\mathcal{U}}_{n}(\theta)) = k_{2,n}/n$. Similarly $%
\tilde{I}_{n,t}^{(\mathcal{E})}(\theta) := I(-\mathcal{\tilde{L}}%
_{n}(\theta) \leq \mathcal{\tilde{E}}_{t}(\theta) \leq %
\widetilde{\mathcal{U}}_{n}(\theta))$. Define trimmed variants
$\widehat{%
\tilde{m}}_{n,t}(\theta) := \tilde{m}_{t}(\theta)\widehat{\tilde{I}}%
_{n,t}^{(\epsilon)}(\theta)$ and $\tilde{m}_{n,t}(\theta) := \tilde
{m}%
_{t}(\theta)\times\break \tilde{I}_{n,t}^{(\epsilon)}(\theta)$.\vadjust{\goodbreak}

%le8.7 #&#
%
\begin{lemma}[(Stationary solution)]
\label{lm:sol} Let $a_{t}(\theta) \in
\{h_{t}(\theta),h_{i,t}^{\theta}(\theta),h_{i,j,t}^{\theta,\theta
}(\theta)\}$ and $a_{t}^{\ast}(\theta) \in \{h_{t}^{\ast}(\theta
), h_{i,t}^{\ast\theta}(\theta), %
h_{i,j,t}^{\ast\theta,\theta}(\theta)\}$.
\begin{enumerate}[(d)]
 \item[(a)] A stationary and ergodic solution $a_{t}^{\ast}(\theta)$
exists for each $\theta \in \Theta$,
it is $\sigma(y_{\tau} \dvtx  \tau \leq t - 1)$-measurable, and $\inf_{\theta\in\Theta}a_{t}^{\ast}(\theta)
> 0$ a.s. Further, $h_{t}^{\ast}(\theta
^{0}) = \sigma_{t}^{2}$ a.s., $%
h_{t}^{\ast\theta}(\theta) = (\partial/\partial
\theta)h_{t}^{\ast}(\theta)$ and $h_{t}^{\ast\theta,\theta
}(\theta) = (\partial/\partial\theta)h_{t}^{\ast
\theta}(\theta)$ a.s.

\item[(b)] $E[\sup_{\theta\in\Theta}|a_{t}^{\ast}(\theta)|^{\iota
}] < \infty$ for some tiny $\iota > 0$.

 \item[(c)] If $a_{t}(\theta)$ is any other stationary solution
then $E[(\sup_{\theta\in\Theta}|a_{t}^{\ast}(\theta) -
a_{t}(\theta)|)^{\iota}] = \mathrm{o}(\rho
^{t}) $ for some $\rho \in (0,1)$.

\item[(d)] $E[\sup_{\theta\in\Theta}|w_{t}^{\ast}(\theta) - w_{t}(\theta
)|] = \mathrm{o}(\rho^{t})$
for each $w_{t}(\theta) \in \{\mathfrak{s}_{i,t}(\theta),%
\mathfrak{d}_{i,j,t}(\theta)\}$.

\item[(e)] $1/n\sum_{t=1}^{n}E[\sup_{\theta\in\Theta}|\tilde{I}%
_{n,t}^{(\mathcal{E})}(\theta) - I_{n,t}^{(\mathcal{%
E})}(\theta)|]$ and $1/n\sum_{t=1}^{n}E[\sup_{\theta\in\Theta}|%
\widehat{\tilde{I}}_{n,t}^{(\mathcal{E})}(\theta) -
\hat{I}_{n,t}^{(\mathcal{E})}(\theta)|]$ are~$\mathrm{o}(1)$.
\end{enumerate}
\end{lemma}

\begin{pf*}{Proof of Theorem~\ref{th:feas}}
We first
characterize properties of random variables based on $h_{t}(\theta)$. We~then prove consistency of the feasible QMTTL estimator $\tilde{\theta
}_{n} \stackrel{p}{\rightarrow} \theta^{0}$. Lastly, we
prove the claim $\mathcal{V}_{n}^{1/2}(\tilde{\theta}_{n} - \hat{\theta}_{n}) \stackrel{p}{\rightarrow} 0$.\vspace*{2pt}

Define $\widehat{\tilde{m}}_{n}(\theta) := 1/n\sum_{t=1}^{n}\widehat{%
\tilde{m}}_{n,t}(\theta)$, $\tilde{m}_{n}(\theta) := 1/n\sum_{t=1}^{n}
\tilde{m}_{n,t}(\theta)$, $\mathcal{\tilde{M}}_{n}(\theta) := %
1/n\times\break \sum_{t=1}^{n}\hspace*{-0.5pt}E[\tilde{m}_{n,t}(\theta)]$, and $\mathfrak{\tilde
{e}}_{n} := 1/n\sum_{t=1}^{n}\hspace*{-0.5pt}\sup_{\theta\in\Theta}E\|\tilde
{m}_{n,t}(\theta
)\| $, and recall $\mathfrak{e}_{n} :=\break \sup_{\theta\in\Theta
}E\|m_{n,t}(\theta)\|$.

\textit{Step 1:}\vspace*{1.5pt}
Use Lemma~\ref{lm:sol} to obtain $|\mathfrak{\tilde{e}%
}_{n} - \mathfrak{e}_{n}| \leq \sup_{\theta\in\Theta
}1/n\sum_{t=1}^{n}\|\tilde{m}_{n,t}(\theta) - m_{n,t}(\theta)]\| = %
\mathrm{O}_{p}(1/n) = \mathrm{o}_{p}(1)$. Similarly, $\|1/n\sum_{t=1}^{n}\|\tilde{G}%
_{t}(\theta)\tilde{I}_{n,t}^{(\mathcal{E})}(\theta) - G_{t}(\theta
)I_{n,t}^{(\mathcal{E})}(\theta)\|$, $\|1/n\sum_{t=1}^{n}\|\tilde{G}%
_{t}(\theta)\*\widehat{\tilde{I}}_{n,t}^{(\mathcal{E})}(\theta) - %
G_{t}(\theta)\hat{I}_{n,t}^{(\mathcal{E})}(\theta)]\|$, and
$\|\mathcal{%
\tilde{G}}(\theta) - \mathcal{G}(\theta)\|$ are uniformly $\mathrm{o}_{p}(1)$,
and for any sequence of positive numbers $\{g_{n}\}$, $g_{n}\rightarrow
\infty$, $\sup_{\theta\in\Theta}\{1/g_{n}\sum_{t=1}^{n}\mid |\tilde
{%
\epsilon}_{t}(\theta)| - |\epsilon_{t}(\theta)|\mid = \mathrm{o}_{p}(1)$.
Use the latter to deduce $\sup_{\theta\in\Theta}k_{n}^{1/2}|\mathcal{
\tilde{E}}_{(k_{n})}^{(a)}(\theta) - \mathcal{E}_{(k_{n})}^{(a)}(\theta
)| \stackrel{p}{\rightarrow} 0$, hence by Lemma B.2 $\sup_{\theta\in
\Theta}|\mathcal{\tilde{E}}_{(k_{1,n})}^{(-)}(\theta)/\mathcal{L}%
_{n}(\theta) + 1| = \mathrm{O}_{p}(1/k_{1,n}^{1/2})$ and $\sup_{\theta\in
\Theta}|\mathcal{\tilde{E}}_{(k_{2,n})}^{(+)}(\theta)/\mathcal{U}%
_{n}(\theta) - 1| = \mathrm{O}_{p}(1/k_{1,n}^{1/2})$. By similar arguments
and Lemma~\ref{lm:sol} it is straightforward to verify Lemmas \ref
{lm:approx}%
, \ref{lm:ulln} and \ref{lm:expan} extend to $\widehat{\tilde
{m}}_{n}(\theta
)$ and $\tilde{m}_{n,t}(\theta)$.

\textit{Step 2 ($\tilde{\theta}_{n} \stackrel{p}{\rightarrow} \theta
^{0}$):}
We follow the proof of Theorem~\ref{th:consist}. By the Lemma~\ref{lm:jac}(c), (d) arguments and $\|\mathcal{\tilde{G}}
- \mathcal{G}\| = \mathrm{o}(1)$ it follows $1/n\sum_{t=1}^{n}E[\tilde{m}%
_{n,t} ( \theta ) ] = \mathcal{G} \times (\theta - %
\theta^{0}) \times (1 + \mathrm{o}(1))$ and $\mathfrak{\tilde{e}}_{n} %
\leq K\|\mathcal{G}\|$. Since $\|\mathcal{G}\| > 0$ it follows $%
\tilde{\epsilon}(\delta) := \liminf_{n\rightarrow\infty
}\inf_{\|\theta-\theta^{0}\|>\delta}\{\mathfrak{\tilde{e}}%
_{n}^{-1}\|1/n\times\break \sum_{t=1}^{n}E[\tilde{m}_{n,t}(\theta)]\|\} > 0$ for
every $n \geq N$ and $\delta > 0$. Therefore $P(\|\tilde{\theta}%
_{n} - \theta^{0}\| > \delta) \leq P(\|\mathcal{\tilde{M}}%
_{n}(\tilde{\theta}_{n})\|/\mathfrak{\tilde{e}}_{n} > \tilde{\epsilon}%
(\delta))$. It remains to show$\|\mathcal{\tilde{M}}_{n}(\tilde{\theta}
_{n})\|/\mathfrak{\tilde{e}}_{n}= \mathrm{o}_{p}(1)$.\vspace*{3pt}

Note $\|\mathcal{\tilde{M}}_{n}(\tilde{\theta}_{n})\| \leq \|\widehat{%
\tilde{m}}_{n}(\tilde{\theta}_{n})\| + \|\widehat{\tilde{m}}_{n}(\tilde{
\theta}_{n}) - \mathcal{\tilde{M}}_{n}(\tilde{\theta}_{n})\|$, where $%
\widehat{\tilde{m}}_{n}(\tilde{\theta}_{n}) = 0$ a.s. by $\tilde{\theta%
}_{n}$ a minimizer. It remains to show $\|\widehat{\tilde{m}}_{n}(\tilde
{%
\theta}_{n}) - \mathcal{\tilde{M}}_{n}(\tilde{\theta}_{n})\|/\mathfrak{%
\tilde{e}}_{n} = \mathrm{o}_{p}(1)$. Note
\begin{eqnarray*}
&&\sup_{\theta\in\Theta}\bigl\llVert \widehat{\tilde{m}}_{n}(
\theta)-%
\mathcal{\tilde{M}}_{n}(\theta)\bigr\rrVert
\\
&&\quad\leq\sup
_{\theta\in\Theta
}\bigl\llVert \widehat{\tilde{m}}_{n}(\theta)-
\hat{m}_{n}(\theta)\bigr\rrVert +\sup_{\theta\in\Theta}\bigl
\llVert \hat{m}_{n}(\theta)-\mathcal{M}%
_{n}(
\theta)\bigr\rrVert +\sup_{\theta\in\Theta}\bigl\llVert
\mathcal{%
\tilde{M}}_{n}(\theta)-\mathcal{M}_{n}(
\theta)\bigr\rrVert .
\end{eqnarray*}
The first and third terms on the right-hand side are $\mathrm{o}_{p}(1)$ by Step 1.
The second is $\mathrm{o}_{p}(\mathfrak{e}_{n})$ by the proof of Theorem~2.1.
Since $|%
\mathfrak{\tilde{e}}_{n} - \mathfrak{e}_{n}| = \mathrm{o}_{p}(1)$ we have
shown $\sup_{\theta\in\Theta}\|\widehat{\tilde{m}}_{n}(\theta) - %
\mathcal{\tilde{M}}_{n}(\theta)\| = \mathrm{o}_{p}(\mathfrak{\tilde{e}}_{n})$
hence $\|\widehat{\tilde{m}}_{n}(\tilde{\theta}_{n}) - \mathcal{\tilde
{M}%
}_{n}(\tilde{\theta}_{n})\|/\mathfrak{\tilde{e}}_{n} = \mathrm{o}_{p}(1)$ as
required.

\textit{Step 3 ($\mathcal{V}_{n}^{1/2}(\tilde{\theta}_{n} - \hat{\theta
}_{n}) \stackrel{p}{\rightarrow} 0$):}
The first order conditions are $\sum_{t=1}^{n}\hat{m}_{n,t}(\hat{%
\theta}_{n}) = 0$ a.s. and $\sum_{t=1}^{n}\widehat{\tilde{m}}_{n,t}(%
\tilde{\theta}_{n}) = 0$ a.s. Combine $\tilde{\theta}_{n} \stackrel{p}%
{\rightarrow} \theta^{0}$, $\sup_{\theta\in\Theta
}\|1/n\sum_{t=1}^{n}\|\tilde{G}_{t}(\theta)\widehat{\tilde{I}}_{n,t}^{(
\mathcal{E})}(\theta) - G_{t}(\theta)\hat{I}_{n,t}^{(\mathcal{E}%
)}(\theta)]\| = \mathrm{o}_{p}(1)$, and $\|\mathcal{\tilde{G}} - \mathcal{G}%
\| = \mathrm{o}(1)$ to deduce by Lemma~\ref{lm:jac} $1/n\sum_{t=1}^{n}\tilde{G}%
_{t}(\tilde{\theta}_{n})\widehat{\tilde{I}}_{n,t}^{(\mathcal{E})} = %
\mathcal{G} \times (1 + \mathrm{o}_{p}(1))$. Therefore, in view of
consistency of the infeasible estimator $\hat{\theta}_{n} \stackrel{p}{%
\rightarrow} \theta^{0}$, expansion Lemma~\ref{lm:expan}, and the
construction $\mathcal{V}_{n} = n\mathcal{G}\Sigma_{n}^{-1}\mathcal{G}$,
it follows
%e3 #&#
\begin{eqnarray}\label{nSmm1}
\frac{1}{n^{1/2}}\Sigma_{n}^{-1/2}\sum
_{t=1}^{n} \bigl\{ \widehat{\tilde {m}}%
_{n,t}(\tilde{\theta}_{n})-\hat{m}_{n,t}(\tilde{
\theta}_{n}) \bigr\} &=& %%
\frac{1}{n^{1/2}}
\Sigma_{n}^{-1/2}\sum_{t=1}^{n}
\bigl\{ \hat {m}_{n,t}(\hat{%
\theta}_{n})-
\hat{m}_{n,t}(\tilde{\theta}_{n}) \bigr\} \nonumber
\\[-8pt]\\[-8pt]
&=&\mathcal{V}_{n}^{1/2} ( \hat{\theta}_{n}-\tilde{
\theta}_{n} ) \bigl( 1+\mathrm{o}_{p} ( 1 ) \bigr) .
\nonumber
\end{eqnarray}
Further, by two applications of Lemmas \ref{lm:approx}(a), \ref
{lm:expan} and %
\ref{lm:jac}, and cancelling the terms $\mathcal{V}_{n}^{1/2}(\tilde
{\theta}%
_{n} - \theta^{0}) = n^{1/2}\Sigma_{n}^{-1/2}\mathcal{G}(\tilde{%
\theta}_{n} - \theta^{0})$, we have
%e4 #&#
\begin{eqnarray}\label{nSmm2}
&&\frac{1}{n^{1/2}}\Sigma_{n}^{-1/2}\sum
_{t=1}^{n} \bigl\{ \widehat {\tilde{m}%
}_{n,t}(\tilde{\theta}_{n})-\hat{m}_{n,t}(\tilde{
\theta}_{n}) \bigr\}\nonumber
\\
&&\quad =\frac{1}{n^{1/2}}\Sigma _{n}^{-1/2}\sum
_{t=1}^{n} \bigl\{ \widehat{\tilde{m}}_{n,t}(
\tilde{\theta }%
_{n})-\tilde{m}_{n,t} \bigr\}
\nonumber
\\
&&\quad\quad{} -\frac{1}{n^{1/2}}\Sigma _{n}^{-1/2}\sum
_{t=1}^{n} \bigl\{ \hat{m}_{n,t}(\tilde{
\theta}%
_{n})-m_{n,t} \bigr\} +\frac{1}{n^{1/2}}
\Sigma _{n}^{-1/2}\sum_{t=1}^{n}
\{ \tilde{m}_{n,t}-m_{n,t} \}
\nonumber
\\[-8pt]\\[-8pt]
&&\quad =\frac{1}{n^{1/2}}\Sigma _{n}^{-1/2}\sum
_{t=1}^{n} \{ \tilde{m}_{n,t}-m_{n,t}
\}
\nonumber
\\
&&\quad\quad{} +n^{1/2}\Sigma_{n}^{-1/2}%
\mathcal{G}
\bigl( \tilde{\theta}_{n}-\theta^{0} \bigr) \bigl(
1+\mathrm{o}_{p} ( 1 ) \bigr) -n^{1/2}\Sigma_{n}^{-1/2}
\mathcal{G} \bigl( \tilde {\theta}%
_{n}-\theta^{0}
\bigr) \bigl( 1+\mathrm{o}_{p} ( 1 ) \bigr)
\nonumber
\\
&&\quad =\frac{1}{n^{1/2}}\Sigma _{n}^{-1/2}\sum
_{t=1}^{n} \{ \tilde{m}_{n,t}-m_{n,t}
\} +\mathrm{o}_{p} \bigl( \bigl\llVert \mathcal{V}_{n}^{1/2}
\bigl( \tilde{\theta}%
_{n}-\theta^{0} \bigr) \bigr
\rrVert \bigr) .
\nonumber
\end{eqnarray}
Combine (\ref{nSmm1}), (\ref{nSmm2}) and Theorem~\ref{th:norm} to obtain $\mathcal
{V}_{n}^{1/2}(%
\tilde{\theta}_{n}-\hat{\theta}_{n}) = n^{-1/2}\Sigma
_{n}^{-1/2}\sum_{t=1}^{n}\{\tilde{m}_{n,t} - m_{n,t}\}(1 + \mathrm{o}_{p}(1))$%
. By Lo\`{e}ve's inequality, $\lim\inf_{n\rightarrow\infty}\|\Sigma
_{n}\| > 0$ in view of non-degeneracy and trimming negligibility, and
Lemma~\ref{lm:sol}(d), it follows for tiny $\iota > 0$, $\rho \in %
(0,1)$, and\vadjust{\goodbreak} sufficiently large $n$ and $K$
\[
E\Biggl\llvert \frac{1}{n^{1/2}}\Sigma_{n}^{-1/2}\sum
_{t=1}^{n} \{ \tilde{m%
}_{n,t}-m_{n,t} \} \Biggr\rrvert ^{\iota}\leq K
\frac{1}{n^{\iota
/2}}%
\sum_{t=1}^{n}E
\llvert \tilde{m}_{n,t}-m_{n,t}\rrvert ^{\iota
}\leq
K%
\frac{1}{n^{\iota/2}}\sum_{t=1}^{n}
\rho^{t}=\mathrm{o} ( 1 ) .
\]
Therefore, $\mathcal{V}_{n}^{1/2}(\tilde{\theta}_{n} - \hat{\theta
}_{n}) = \mathrm{o}_{p}(1)$ by Markov's inequality.
\end{pf*}

\begin{pf*}{Proof of Theorem~\ref{th:mnwm}}
By Assumption~\ref{ass4} $\psi(u,c) %
= u\varpi(u,c)I(|u| \leq c)$ behaves like $uI(|u| \leq c)$ as $%
c \rightarrow \infty$. See (\ref{psi_I}). In the following, we
therefore only treat the simple trimming transform $\psi(u,c) = uI(|u|
\leq c)$. The general case with properties (\ref{psi_I}) and (\ref{dh})
has a similar proof.

Lemmas \ref{lm:approx}--\ref{lm:jac} extend to cover the equations
\begin{eqnarray*}
\widehat{\check{m}}_{n,t}(\theta)&=& \Biggl( \epsilon_{t}^{2}(
\theta )\hat{I}%
_{n,t}^{(\epsilon)}(\theta)-
\frac{1}{n}\sum_{t=1}^{n}\epsilon
_{t}^{2}(\theta)\hat{I}_{n,t}^{(\epsilon)}(
\theta) \Biggr) \times \mathfrak{s}_{t}(\theta),
\\
\check{m}_{n,t}(\theta)&=& \bigl( \epsilon_{t}^{2}(
\theta )I_{n,t}^{(\epsilon)}(\theta)-E \bigl[ \epsilon_{t}^{2}(
\theta )I_{n,t}^{(\epsilon)}(\theta) \bigr] \bigr) \times \bigl(
\mathfrak {s}%
_{t}(\theta)-E \bigl[ \mathfrak{s}_{t}(
\theta) \bigr] \bigr) .
\end{eqnarray*}
Consider Lemma~\ref{lm:approx}(a). By Lemma~\ref{lm:approx}, it follows
\begin{eqnarray*}
&&\Sigma_{n}^{-1/2}\frac{1}{n^{1/2}}\sum
_{t=1}^{n} \Biggl( \epsilon _{t}^{2}%
\hat{I}_{n,t}^{(\epsilon)}-\frac{1}{n}\sum
_{t=1}^{n}\epsilon _{t}^{2}
\hat{I}%
_{n,t}^{(\epsilon)} \Biggr) \times
\mathfrak{s}_{t}
\\
&&\quad=\Sigma _{n}^{-1/2}\frac{1%
}{n^{1/2}}
\sum_{t=1}^{n} \Biggl( \epsilon_{t}^{2}I_{n,t}^{(\epsilon
)}-
\frac{1%
}{n}\sum_{t=1}^{n}
\epsilon_{t}^{2}I_{n,t}^{(\epsilon)} \Biggr)
\times \mathfrak{s}_{t}+\mathrm{o}_{p} ( 1 ) ,
\end{eqnarray*}
where by independence and dominated convergence $\Sigma_{n} \sim %
E[(\epsilon_{t}^{2}I_{n,t}^{(\epsilon)} - E[\epsilon
_{t}^{2}I_{n,t}^{(\epsilon)}])^{2}] \times E[\mathfrak{s}_{t}\mathfrak
{%
s}_{t}^{\prime}] =: \sigma_{n}^{2}\mathfrak{S}$. Now add and subtract
$%
E[\epsilon_{t}^{2}I_{n,t}^{(\epsilon)}]$ and $E[\mathfrak{s}_{t}]$ to
deduce
%e5 #&#
\begin{eqnarray}\label{Ses}
&&\Sigma_{n}^{-1/2}\frac{1}{n^{1/2}}\sum
_{t=1}^{n} \Biggl( \epsilon _{t}^{2}%
\hat{I}_{n,t}^{(\epsilon)}-\frac{1}{n}\sum
_{t=1}^{n}\epsilon _{t}^{2}
\hat{I}%
_{n,t}^{(\epsilon)} \Biggr) \times
\mathfrak{s}_{t}\nonumber
\\
&& \quad=\Sigma_{n}^{-1/2}\frac{1}{n^{1/2}}%
\sum
_{t=1}^{n} \bigl( \epsilon_{t}^{2}I_{n,t}^{(\epsilon)}-E
\bigl[ \epsilon _{t}^{2}I_{n,t}^{(\epsilon)}
\bigr] \bigr) \times \bigl( \mathfrak {s}_{t}-E%
[
\mathfrak{s}_{t} ] \bigr) +\mathrm{o}_{p} ( 1 )
\\
&& \quad\quad{}-\frac{1}{\sigma_{n}n^{1/2}}%
\sum_{t=1}^{n}
\bigl( \epsilon_{t}^{2}I_{n,t}^{(\epsilon)}-E
\bigl[ \epsilon _{t}^{2}I_{n,t}^{(\epsilon)}
\bigr] \bigr) \times\mathfrak{S}\times \Biggl( \frac{1}{n}\sum
_{t=1}^{n}\mathfrak{s}_{t}-E [ \mathfrak
{s}_{t}%
] \Biggr) \times \bigl( 1+\mathrm{o}_{p} ( 1 )
\bigr) .
\nonumber
\end{eqnarray}
Under Assumption~\ref{ass1} $\mathfrak{s}_{t}$ is stationary, ergodic and integrable,
hence $1/n\sum_{t=1}^{n}\mathfrak{s}_{t} - E[\mathfrak{s}_{t}] = %
\mathrm{o}_{p}(1)$, and by a generalization of central limit theorem
Lemma~\ref%
{lm:clt} $\sigma_{n}^{-1}n^{-1/2}\sum_{t=1}^{n}(\epsilon
_{t}^{2}I_{n,t}^{(\epsilon)} - E[\epsilon_{t}^{2}I_{n,t}^{(\epsilon
)}]) = \mathrm{O}_{p}(1)$. The\vspace*{2pt} second term in (\ref{Ses}) is therefore
$\mathrm{o}_{p}(1)$, hence $\Sigma_{n}^{-1/2}n^{-1/2}\*\sum_{t=1}^{n}(\widehat{\check
{m}}_{n,t} - \check{m}_{n,t}) = \mathrm{o}_{p}(1)$ which extends Lemma~\ref{lm:approx}(a)
to $\{\widehat{\check{m}}_{n,t},\check{m}_{n,t}\}$. In view of
$L_{2+\iota}$%
-boundedness of $\sup_{\theta\in\mathcal{N}_{0}}\|\mathfrak
{s}_{t}(\theta
)\|$ for some compact $\mathcal{N}_{0} \subset \Theta$ with positive
Lebesgue measure and containing $\theta^{0}$, and independence of $%
\epsilon_{t}$, the arguments used to prove\vspace*{2pt} Lemmas \ref{lm:approx}(b),
\mbox{\ref{lm_var_bnds}--\ref{lm:jac}} carry over with simple modifications to
cover $\{%
\widehat{\check{m}}_{n,t},\check{m}_{n,t}\}$. The claims therefore
follow by
imitating the proofs of Theorems 2.1 and 2.2, and by the constructions
of $%
\mathcal{\accentset{\sdbullet}{V}}_{n}$ and $\mathcal{V}_{n}$.
\end{pf*}

%le8.8 #&#
%
\begin{lemma}
\label{lm:mom_est}$1/n\sum_{t=1}^{n}\mathcal{E}_{t}^{2}(\hat{\theta
}_{n})%
\hat{I}_{n,t}^{(\mathcal{E})}(\hat{\theta}_{n})/E[\mathcal{E}%
_{t}^{2}I_{n,t}^{(\mathcal{E})}] \stackrel{p}{\rightarrow} 1$.
\end{lemma}

\begin{pf*}{Proof of Theorem~\ref{th:scale}}
The claim follows
from Jacobian consistency Lemma~\ref{lm:jac}(b) and Lemma~\ref{lm:mom_est}.
\end{pf*}

%s8.3 #&#
\subsection{Proofs of supporting lemmas}
\label{app:lemmas}

In order to decrease the number of cases we augment Assumption~\ref{ass1}(b) and
impose power law tails on $\epsilon_{t}$ in general:%
%e6 #&#
\begin{equation}
P \bigl( \llvert \epsilon_{t}\rrvert >a \bigr) =da^{-\kappa
}
\bigl( 1+\mathrm{o} ( 1 ) \bigr) \quad\quad\mbox{where }d\in ( 0,\infty ) \mbox{ and }\kappa\in ( 2,
\infty ) . \label{pow_e}
\end{equation}
Notice $\epsilon_{t}(\theta)$ is stationary and ergodic on $\Theta$
by (%
\ref{theta}), and also has a power law tail. The latter follows by
noting
\[
\epsilon_{t}(\theta)=\frac{y_{t}}{\sigma_{t}(\theta)}=\epsilon _{t}
\frac{%
\sigma_{t}}{\sigma_{t}(\theta)}\quad\mbox{and}\quad\mathcal{E}_{t}(\theta )=
\epsilon_{t}^{2}(\theta)-1,
\]
where $E(\sup_{\theta\in\Theta}|\sigma_{t}^{2}/\sigma
_{t}^{2}(\theta
)|)^{p} < \infty$ for any $p > 0$ under Assumption~\ref{ass1}. Since $%
\epsilon_{t}$ is independent of $\sigma_{t}/\sigma_{t}(\theta)$ the
product convolution $\epsilon_{t} \times (\sigma_{t}/\sigma
_{t}(\theta))$ has tail (\ref{pow_e}) with the same index $\kappa >
2$ %
(Breiman \cite{Breiman65}). In general $\lim_{a\rightarrow\infty}\sup_{\theta\in
\Theta}\{|c^{\kappa}P(|\epsilon_{t}(\theta)| > a) - d(\theta
)|\} = 0$ and $\inf_{\theta\in\Theta}\{d(\theta)\} > 0$ and $%
\sup_{\theta\in\Theta}\{d(\theta)\} < \infty$. Hence, in view of (%
\ref{pow_E}), $\mathcal{E}_{t}(\theta) := \epsilon_{t}^{2}(\theta) %
- 1$ also satisfies
%e7 #&#
\begin{eqnarray}\label{pow_E_T}
&&\lim_{a\rightarrow\infty}\sup_{\theta\in\Theta} \bigl\{ \bigl\llvert
a^{\kappa/2}P \bigl( \bigl\llvert \mathcal{E}_{t}(\theta)\bigr
\rrvert >a \bigr) -d(\theta)\bigr\rrvert \bigr\} =0\nonumber
\\[-8pt]\\[-8pt]
&&\quad \mbox{where } \inf
_{\theta
\in\Theta} \bigl\{ d(\theta) \bigr\} >0\mbox{ and }\sup
_{\theta\in
\Theta} \bigl\{ d(\theta) \bigr\} <\infty. \nonumber
\end{eqnarray}
Recall $P(|\mathcal{E}_{t}(\theta)| > \mathcal{C}_{n}(\theta)) = %
k_{n}/n$ holds for $\mathcal{C}_{n}(\theta) = \mathcal{U}_{n}(\theta)$
and $k_{n} = k_{2,n}$. Then by (\ref{pow_E_T})%
%e8 #&#
\begin{equation}
\mathcal{C}_{n}(\theta)=d(\theta)^{2/\kappa} ( n/k_{n} )
^{2/\kappa}. \label{FT}
\end{equation}
Further, by (\ref{pow_E_T}) and an application of Karamata's theorem:%
%e9 #&#
\begin{eqnarray}\label{Karam}
&&\mbox{if }\kappa=4\mbox{:}\quad E \bigl[ \mathcal{E}_{t}^{2}(\theta)I
\bigl( \bigl\llvert \mathcal{E}_{t}(\theta)\bigr\rrvert \leq\mathcal
{C}_{n}(\theta ) \bigr) \bigr] \sim d(\theta)\ln(n),\nonumber
\\
&&\mbox{if }\kappa<4\mbox{:}\quad E \bigl[ \mathcal{E}_{t}^{2}(\theta)I
\bigl( \bigl\llvert \mathcal{E}_{t}(\theta)\bigr\rrvert \leq\mathcal
{C}_{n}(\theta ) \bigr) \bigr] \sim\frac{\kappa}{4-\kappa}\mathcal
{C}_{n}^{2}(\theta )P \bigl( \bigl\llvert
\mathcal{E}_{t}(\theta)\bigr\rrvert >\mathcal{C}%
_{n}(\theta) \bigr)
\\
&&\hphantom{\mbox{if }\kappa<4\mbox{:}\quad E \bigl[ \mathcal{E}_{t}^{2}(\theta)I
\bigl( \bigl\llvert \mathcal{E}_{t}(\theta)\bigr\rrvert \leq\mathcal
{C}_{n}(\theta ) \bigr) \bigr]}=\frac{\kappa}{4-\kappa}d(\theta)^{4/\kappa
}(n/k_{n})^{4/\kappa-1}.
\nonumber
\end{eqnarray}
Uniform bounds are similar given (\ref{pow_E_T})--(\ref{Karam}). For example,
when $\kappa < 4$:
%e10 #&#
\begin{equation}
\sup_{\theta\in\Theta} \biggl\{ \frac{n}{k_{n}}\frac{\mathcal{C}%
_{n}^{2}(\theta)}{E [ \mathcal{E}_{t}^{2}(\theta)I ( \llvert
\mathcal{E}_{t}(\theta)\rrvert \leq\mathcal{C}_{n}(\theta)
) %
] }
\biggr\} \rightarrow ( 0,\infty ) . \label{TM}
\end{equation}

Unless otherwise noted, and in view of (\ref{pow_E}) and (\ref
{pow_E_T}), we
assume two-tailed trimming to reduce notation, hence thresholds and
fractiles are simply $\mathcal{C}_{n}(\theta)$ and $k_{n}$, and order
statistics are $\mathcal{E}_{(k_{n})}^{(a)}(\theta)$ where $\mathcal{E}
_{t}^{(a)}(\theta) := |\mathcal{E}_{t}(\theta)|$.

The proofs of Lemmas \ref{lm:approx}--\ref{lm:mom_est} require two supporting
results. See the supplementary material Hill \cite{Hillsupp} for proofs. First,
trimming indicators satisfy a uniform CLT.\def\theLemma{B.\arabic{Lemma}}

\begin{Lemma}[(Uniform indicator CLT)]\label{lB.1}
Define $%
\mathcal{I}_{n,t}(\theta) := ((n/k_{n})^{1/2})\{I(|\mathcal{E}%
_{t}(\theta)| \leq \mathcal{C}_{n}(\theta)) -
E[I(|\mathcal{E}_{t}(\theta)| \leq \mathcal{C}_{n}(\theta))]\}$.
Then $\{n^{-1/2}\sum_{t=1}^{n}\mathcal{I}_{n,t}(\theta) \dvtx  \theta
\in \Theta\}
\Longrightarrow^{\ast} \{\mathcal{I}(\theta) \dvtx  \theta \in \Theta\}
$, where
$\mathcal{I}(\theta)$ is a Gaussian process with uniformly bounded
and uniformly continuous sample paths with respect to $L_{2}$-norm,
and $\Longrightarrow^{\ast}$ denotes weak convergence on a
Polish space (Hoffman-J{\o}rgensen \cite{HoffJorg1984}).
\end{Lemma}

Second, intermediate order statistics are uniformly bounded in probability.

\begin{Lemma}[(Uniform order statistic bound)]\label{lB.2} $%
\sup_{\theta\in\Theta}|\mathcal{E}_{(k_{n})}^{(a)}(\theta)/\mathcal
{C}%
_{n}(\theta) - 1| =
\mathrm{O}_{p}(1/k_{n}^{1/2})$.
\end{Lemma}

Lemmas \ref{lm:approx}, \ref{lm:ulln}, \ref{lm:expan} and \ref{lm:jac} are
similar to results proven in Hill \cite{Hillltts}, Appendix A,
hence their proofs are relegated to the supplementary material Hill \cite
{Hillsupp}.

\begin{pf*}{Proof of Lemma~\ref{lm_var_bnds}}

\textit{Claim (a):} $\mathfrak{s}_{t}$ is $L_{2+\iota}$-bounded by
Assumption~\ref{ass1}, hence by error independence $\Sigma_{i,i,n} \sim E[%
\mathcal{E}_{t}^{2}I_{n,t}^{(\mathcal{E})}] \times E[\mathfrak{s}%
_{i,t}^{2}] \sim KE[\mathcal{E}_{t}^{2}I_{n,t}^{(\mathcal{E})}]$. The
claim now follows from arguments leading to Theorem~\ref{th:converg}%
.

\textit{Claim (b):} We prove the claim for $\mathcal{S}_{i,i,n}$, so
let $m_{t}$ denote $m_{i,t}$, hence $\mathfrak{s}_{t}$ denotes
$\mathfrak{s}%
_{i,t}$, and express $\mathcal{S}_{i,i,n}$ as $\mathcal{S}_{n}$. Note $%
\mathcal{S}_{n} \sim E[m_{n,t}^{2}]+2\sum_{i=1}^{n-1}(1 - %
i/n)E[m_{n,1}m_{n,i+1}]$. If $E[m_{t}^{2}] < \infty$ then $\mathcal{S}%
_{n}\sim K = \mathrm{o}(n/\ln(n))$ in view of geometric $\beta$-mixing %
(cf. Ibragimov \cite{Ibrag1962}).

Now assume $E[m_{t}^{2}] = \infty$. We first characterize the tails of
$%
m_{t} = (\epsilon_{t}^{2} - 1)\mathfrak{s}_{t}$, and then bound $%
\sum_{i=1}^{n-1}|E[m_{n,1}m_{n,i+1}]|$.

\textit{Step 1:} By Assumption~\ref{ass1}(b) and (\ref{pow_e}) independent $%
\epsilon_{t}$ has a power law tail with index $\kappa \in (2,4]$, and
since $\alpha^{0} + \beta^{0} > 0$ it follows $E[\mathfrak{s}%
_{t}^{2}] < \infty$. Therefore $m_{t}$ has a power law tail with index
$%
\kappa_{m} := \kappa/2 \in (1,2]$, cf. Breiman \cite{Breiman65}.

\textit{Step 2:} Define quantile functions $Q_{n}(u) = \inf\{m %
\geq 0 : P(|m_{n,t}| > m) \leq u\}$ and $Q(u) = \inf\{m \geq 0 \dvtx
P(|m_{t}| > m) \leq u\}$ for $u \in [0,1]$,
recall geometric $\beta$-mixing implies $\alpha$-mixing with
coefficients $%
\alpha_{h} \leq K\rho^{h}$ for $\rho \in (0,1)$. By Theorem~1.1
of Rio \cite{Rio93}
\[
\sum_{i=1}^{n-1}\bigl\llvert E [
m_{n,1}m_{n,i+1} ] \bigr\rrvert \leq 2\sum
_{i=1}^{n-1}\int_{0}^{2\alpha_{i}}Q_{n}^{2}(u)\,\mathrm{d}u
\leq 2\sum_{i=1}^{n-1}\int
_{0}^{K\rho^{i}}Q_{n}^{2}(u)\,\mathrm{d}u.
\]
Tail-trimming $m_{n,t} = m_{t}I_{n,t}^{(\mathcal{E})}$ coupled with
distribution continuity imply $P(m_{n,t} = 0) = k_{n}/n$. Thus $%
Q_{n}(u) = 0$ for $u \in [0,k_{n}/n]$ and $Q_{n}(u) = Q(u)$
for $u \in (k_{n}/n,1]$. Further, under the Step~1 power law properties
$Q(u) = \mathrm{O}(u^{-2/\kappa})$. Therefore
\begin{eqnarray*}
\sum_{i=1}^{n-1}\bigl\llvert E [
m_{n,1}m_{n,i+1} ] \bigr\rrvert &\leq &K\sum
_{i=1}^{n-1}\int_{k_{n}/n}^{K\rho^{i}}u^{-4/\kappa}\,\mathrm{d}u
\leq K\sum_{i=1}^{n-1}\max \bigl\{ 0, (
n/k_{n} ) ^{ ( 4/\kappa
-1 ) }-K\rho^{-i ( 4/\kappa-1 ) } \bigr\}
\\
&=&K\sum_{i=1}^{K\ln(n/k_{n})} \bigl\{ (
n/k_{n} ) ^{ (
4/\kappa-1 ) }-K\rho^{-i ( 4/\kappa-1 ) } \bigr\} .
\end{eqnarray*}
Moreover $\sum_{i=1}^{K\ln(n/k_{n})}\{(n/k_{n})^{4/\kappa-1} - K\rho
^{-i(4/\kappa-1)}\} = K\ln(n/k_{n})\times(n/k_{n})^{4/\kappa-1}(1 %
+ \mathrm{O}(1))$ and $k_{n} = \mathrm{o}(n)$ hence
\begin{eqnarray*}
\sum_{i=1}^{n-1}\bigl\llvert E [
m_{n,1}m_{n,i+1} ] \bigr\rrvert &\leq &K\ln ( n/k_{n}
) \times ( n/k_{n} ) ^{4/\kappa
-1} \bigl( 1+\mathrm{O}(1) \bigr)
\\
&\leq&K\ln ( n/k_{n} ) \times ( n/k_{n} ) ^{4/\kappa
-1}
\leq K\ln ( n ) ( n/k_{n} ) ^{4/\kappa-1}.
\end{eqnarray*}
Further, $\ln(n)(n/k_{n})^{4/\kappa-1} = \mathrm{o}(n/\ln(n))$ since $k_{n} %
\rightarrow \infty$ and $\kappa \in (2,4]$. Finally, by Step 1 and
(\ref{Karam}) $E[m_{n,t}^{2}] \sim K(n/k_{n})^{4/\kappa-1}$ if $\kappa
< 4$ and $E[m_{n,t}^{2}] \sim K\ln(n)$ if $\kappa = 4$.
Therefore $\mathcal{S}_{n} \leq K\ln(n)(n/k_{n})^{4/\kappa-1} = %
\mathrm{o}(n/\ln(n))$ which completes the proof.
\end{pf*}

\begin{pf*}{Proof of Lemma~\ref{lm:clt}} By identification Assumption~\ref{ass2}
$%
n^{-1/2}\Sigma_{n}^{-1/2}\sum_{t=1}^{n}m_{n,t} =
n^{-1/2}\Sigma_{n}^{-1/2}\*\sum_{t=1}^{n}\{m_{n,t} - E[m_{n,t}]\} + %
\mathrm{o}(1)$. Define $z_{n,t} := r^{\prime}\Sigma_{n}^{-1/2}\{m_{n,t} - %
E[m_{n,t}]\}$ for any $r \in \mathbb{R}^{q}$, \mbox{$r^{\prime}r = 1$}.
Note by error independence, dominated convergence and (\ref{sigma}): $%
E(\sum_{t=1}^{n}z_{n,t})^{2} \sim n$. We will prove $%
1/n^{1/2}\sum_{t=1}^{n}z_{n,t} \stackrel{d}{\rightarrow} N(0,1)$, hence
the claim will follow from the Cram\'{e}r--Wold Theorem. Define $\Im
_{t} %
:= \sigma(y_{\tau} \dvtx  \tau \leq t)$.

In view of geometric $\beta$-mixing and stationarity under Assumption~\ref{ass1},
and $E[z_{n,t}^{2}] = 1$ it suffices to show the three conditions of
Theorem~2.1 in Peligrad \cite{Peligrad96} hold.\footnote{We require a result like Theorem~2.1 in Peligrad \cite{Peligrad96} since
asymptotically $z_{n,t}$ need not have finite moments higher than two.} The
first two are $\sup_{n\geq1}1/n\sum_{t=1}^{n}E[z_{n,t}^{2}] < \infty$
and the Lindeberg condition $1/n\sum_{t=1}^{n}E[z_{n,t}^{2}I(|z_{n,t}|
> %
\varepsilon n^{1/2})] \rightarrow 0\ \forall\varepsilon > 0$. By
construction $E[z_{n,t}^{2}] = 1$ hence $1/n\sum_{t=1}^{n}E[z_{n,t}^{2}] = 1 + \mathrm{o}(1)$ which verifies the first.

The Lindeberg condition holds if $\kappa > 4$ since $E|\epsilon
_{t}|^{4+\iota} < \infty$ and $E|\mathfrak{s}_{i,t}|^{2+\iota} < %
\infty$ for some $\iota > 0$, hence $\lim\sup_{n\rightarrow\infty
}E|z_{n,t}|^{2+\iota} < \infty$. Now suppose $\kappa \leq 4$,
assume $E[\mathcal{E}_{t}I_{n,t}^{(\mathcal{E})}] = 0$ to simplify
notation, and note $z_{n,t} = \mathcal{E}_{t}I_{n,t}^{(\mathcal{E}%
)}r^{\prime}\Sigma_{n}^{-1/2}\mathfrak{s}_{t}$. By independence and
$L_{2}$%
-boundedness of $\mathfrak{s}_{t}$ it follows $\Sigma_{n} = E[\mathcal
{E%
}_{t}^{2}I_{n,t}^{(\mathcal{E})}] \times \mathfrak{S}$ where $\mathfrak
{%
S} = E[\mathfrak{s}_{t}\mathfrak{s}_{t}^{\prime}]$ is finite and
positive definite. By construction $\lim\inf_{n\rightarrow\infty
}\inf_{r^{\prime}r=1}\|\Sigma_{n}\|(r^{\prime}\Sigma
_{n}^{-1/2}\mathfrak{%
s}_{t})^{2} > 0$ a.s., by independence $\mathcal{E}_{t}^{2} \times
r^{\prime}\Sigma_{n}^{-1/2}\mathfrak{s}_{t}\mathfrak{s}_{t}^{\prime
}\Sigma_{n}^{-1/2}r\times\|\Sigma_{n}\|$ has Paretian tails with
index $%
\kappa/4 \leq 1$, and by trimming $|\mathcal{E}_{t}^{2}I_{n,t}^{(%
\mathcal{E})}| \leq K\mathcal{C}_{n}^{2}$. Therefore, for finite $K >
0$ that may be different in different places,
\begin{eqnarray*}
E \bigl[ z_{n,t}^{2}I \bigl( z_{n,t}^{2}>
\varepsilon^{2}n \bigr) \bigr] &\leq&KE \biggl( \bigl( r^{\prime}
\Sigma_{n}^{-1/2}\mathfrak {s}_{t} \bigr)
^{2}\times E \biggl[ \mathcal{E}_{t}^{2}I_{n,t}^{(\mathcal{E})}I
\biggl( \mathcal{E}_{t}^{2}I_{n,t}^{(\mathcal{E})}>
\frac{\varepsilon
^{2}n}{ (
r^{\prime}\Sigma_{n}^{-1/2}\mathfrak{s}_{t} ) ^{2}} \biggr) \Big|\Im _{t-1}%
\biggr] \biggr)
\\
&\leq&KE \bigl( \bigl( r^{\prime}\Sigma_{n}^{-1/2}
\mathfrak {s}_{t} \bigr) ^{2}\times E \bigl[
\mathcal{E}_{t}^{2}I_{n,t}^{(\mathcal{E})}I \bigl(
\mathcal{E}_{t}^{2}I_{n,t}^{(\mathcal{E})}>K
\varepsilon^{2}nE \bigl[ \mathcal{E}_{t}^{2}I_{n,t}^{(\mathcal{E})}
\bigr] \bigr) |\Im _{t-1} \bigr] \bigr)
\\
&\leq&KE \biggl( \bigl( r^{\prime}\Sigma_{n}^{-1/2}
\mathfrak {s}_{t} \bigr) ^{2}\times E \biggl[ \int
_{K\varepsilon^{2}nE[\mathcal
{E}_{t}^{2}I_{n,t}^{(%
\mathcal{E})}]}^{K\mathcal{C}_{n}^{2}}u^{-\kappa/4}\,\mathrm{d}u \biggr] \biggr) .
\end{eqnarray*}
In general $\mathcal{C}_{n}^{2} = K(n/k_{n})^{4/\kappa}$. If $\kappa %
= 4$ then $E[\mathcal{E}_{t}^{2}I_{n,t}^{(\mathcal{E})}] \sim K\ln
(n))$ hence $\mathcal{C}_{n}^{2}=K(n/k_{n}) < K\varepsilon^{2}nE[%
\mathcal{E}_{t}^{2}I_{n,t}^{(\mathcal{E})}]$ as $n \rightarrow \infty$. 
This implies for some $N \in \mathbb{N}$ and all $n \geq N$ that\break $%
\int_{K\varepsilon^{2}nE[\mathcal{E}_{t}^{2}I_{n,t}^{(\mathcal
{E})}]}^{K%
\mathcal{C}_{n}^{2}}u^{-\kappa/4}\,\mathrm{d}u = 0$. If $\kappa < 4$
then $E[%
\mathcal{E}_{t}^{2}I_{n,t}^{(\mathcal{E})}] \sim K\mathcal{C}%
_{n}^{2}(k_{n}/n) = K(n/k_{n})^{4/\kappa-1}$, hence again $\mathcal{C}%
_{n}^{2}=K(n/k_{n}) <Kn(n/k_{n})^{4/\kappa-1} = K\varepsilon^{2}nE[%
\mathcal{E}_{t}^{2}I_{n,t}^{(\mathcal{E})}]$ as $n \rightarrow \infty$.
Therefore, $E[z_{n,t}^{2}I(z_{n,t}^{2} > \varepsilon^{2}n)] = 0$
for some $N \in \mathbb{N}$ and all $n \geq N$. This proves $%
1/n\sum_{t=1}^{n}E[z_{n,t}^{2}\times\break I(|z_{n,t}| > \varepsilon n^{1/2})] %
\rightarrow 0$  $\forall\varepsilon > 0$.

The third condition concerns the maximum correlation coefficient $\rho(
\mathcal{A},\mathcal{B}) :=\break  \sup_{f\in L_{2}(\mathcal{A}%
),g\in L_{2}(\mathcal{B})}|\operatorname{corr}(f,g)|$ defined on
$L_{2}(\mathfrak{F})$ the\vspace*{2pt}
space of $L_{2}$-bounded $\mathfrak{F}$-measurable random variables. We
require the interlaced coefficient $\rho_{k}^{\ast} := \sup_{n\geq
1}\sup_{S_{k},T_{k}}\rho(\sigma(z_{n,i} \dvtx  i \in T_{k}), \sigma
(z_{n,j} \dvtx  j \in S_{k}))$ to satisfy $\lim_{k\rightarrow\infty
}\rho_{k}^{\ast} < 1$, where $T_{k} , S_{k} \subset \{1,\ldots,n\}$
are non-empty subsets with $\inf_{s\in S_{k},t\in T_{k}}\{|s - t|\} %
\geq k$, and $\sup_{S_{k},T_{k}}$ is taken over all sets $
\{S_{k},T_{k}\} $ for a given distance $k$. See equations (1.2), (1.7) and
(1.8) in Peligrad \cite{Peligrad96}. In view of the GARCH process and
Assumption~\ref{ass1}, $%
\{z_{n,t} \dvtx  1 \leq t \leq n\}_{n\geq1}$ is a first order
Markov chain that is stationary over $1 \leq t \leq n$, and by
geometric $\beta$-mixing it is also geometric $\alpha$-mixing. Since
$\rho
_{1}^{\ast} < 1$ as a consequence of independence of $\epsilon_{t}$,
it therefore follows $\rho_{k}^{\ast} \rightarrow 0$ by an extension
of Theorem~3.3 in Bradley \cite{Bradley05} to triangular arrays.
\end{pf*}

\begin{pf*}{Proof of Lemma~\ref{lm:sol}} Claims (a)--(c) follow from the
Assumption~\ref{ass3} response Lipschitz properties. See Francq and Zako\"{\i}an
\cite{FZ04,FZ10} and Meitz and Saikkonen \cite{meitzsaikkonen08}. Claim
(d) follows from stationarity, independence of $%
\epsilon_{t}$, and (b) and (c).

Consider (e). We will prove $1/n\sum_{t=1}^{n}E[\sup_{\theta\in\Theta
}|%
\tilde{I}_{n,t}^{(\mathcal{E})}(\theta) - I_{n,t}^{(%
\mathcal{E})}(\theta)|] = \mathrm{o}(1)$, the second claim being similar. We can
approximate $I(u) := I(u \leq 0)$ with the \textit{regular}
sequence $\{\mathfrak{I}_{n}(u)\}_{n\geq1}$, defined by $\mathfrak{I}%
_{n}(u) := \int_{-\infty}^{\infty}I(\varpi)\mathrm{S}(\mathcal{N}%
_{n}(\varpi - u))\mathcal{N}_{n}\mathrm{e}^{-\varpi^{2}/\mathcal{N}%
_{n}^{2}}\,\mathrm{d}\varpi$ where $\mathrm{S}(\xi) = \mathrm{e}^{-1/(1-\xi
^{2})}/\int_{-1}^{1}\mathrm{e}^{-1/(1-w^{2})}\,\mathrm{d}w$ if $|\xi| < 1$ and $\mathrm{S}%
(\xi) = 0$ if $|\xi| \geq 1$. Here $\{\mathcal{N}_{n}\}$ is a
sequence of finite positive numbers, $\mathcal{N}_{n} \rightarrow %
\infty$, the rate to be chosen below. See Lighthill \cite{Lighthill}.
$\mathfrak{I}%
_{n}(u)$ is uniformly bounded in $u$, continuous and differentiable.
Also, $%
(\partial/\partial u)I(u)$ has a regular sequence $\mathfrak{D}_{n}(u)
:= (\mathcal{N}_{n}/\uppi)^{1/2}\exp\{-\mathcal{N}_{n}u^{2}\}$.

Define $\mathfrak{e}_{t}(a) := |\mathcal{E}_{t}| - a$ and $\mathfrak{%
\tilde{e}}_{t}(a) := |\mathcal{\tilde{E}}_{t}| - a$, and let $%
\mathcal{\tilde{C}}_{n}(\theta)$ satisfy $P(|\mathcal{\tilde
{E}}_{t}(\theta
)| \geq \mathcal{\tilde{C}}_{n}(\theta)) = k_{n}/n$. Hence $\tilde{%
I}_{n,t}^{(\mathcal{E})}(\theta) = I(\mathfrak{\tilde{e}}_{t}(\mathcal{
\tilde{C}}_{n}(\theta))$ and $I_{n,t}^{(\mathcal{E})} = I(\mathfrak{e}%
_{t}(\mathcal{C}_{n}(\theta))$. Note $\mathcal{N}_{n} \rightarrow %
\infty$ can be made as fast as we choose such that $\sup_{\theta\in
\Theta
}|I(\mathfrak{e}_{t}(\mathcal{C}_{n}(\theta)) - I(\mathfrak{\tilde{e}}%
_{t}(\mathcal{\tilde{C}}_{n}(\theta))| \leq K\sup_{\theta\in\Theta
}|%
\mathfrak{I}_{n}(\mathfrak{e}_{t}(\mathcal{C}_{n}(\theta)) - \mathfrak
{I%
}_{n}(\mathfrak{\tilde{e}}_{t}(\mathcal{\tilde{C}}_{n}(\theta))| + %
\mathrm{o}_{p}(1)$, and $\mathfrak{D}_{n}(u) \rightarrow 0$ as fast as we
choose. Hence, by the mean-value-theorem and boundedness of $\mathfrak
{D}%
_{n}(u)$ it follows $\sup_{\theta\in\Theta}|I(\mathfrak
{e}_{t}(\mathcal{C}%
_{n}(\theta)) - I(\mathfrak{\tilde{e}}_{t}(\mathcal{\tilde{C}}%
_{n}(\theta))| \leq K\sup_{\theta\in\Theta}|\mathcal{\tilde{E}}%
_{t}(\theta) - \mathcal{E}_{t}(\theta)| + K\sup_{\theta\in\Theta
}|\mathcal{\tilde{C}}_{n}(\theta) - \mathcal{C}_{n}(\theta)|$. By
(iii) $\sup_{\theta\in\Theta}|\mathcal{\tilde{E}}_{t}(\theta) - %
\mathcal{E}_{t}(\theta)| = \mathrm{o}_{p}(\rho^{t})$. Similarly $\sup_{\theta
\in\Theta}\sum_{t=1}^{n}| |\mathcal{\tilde{E}}_{t}(\theta)| - |%
\mathcal{E}_{t}(\theta)| | = \mathrm{O}_{p}(1)$ hence $\sup_{\theta\in\Theta
}|\mathcal{\tilde{E}}_{(k_{n})}^{(a)}(\theta) -\break  \mathcal{E}%
_{(k_{n})}^{(a)}(\theta)| \stackrel{p}{\rightarrow} 0$, hence by Lemma
B.2 $\sup_{\theta\in\Theta}|\mathcal{\tilde{C}}_{n}(\theta)-\mathcal
{C}%
_{n}(\theta)| \rightarrow 0$. Therefore by dominated convergence $%
1/n\sum_{t=1}^{n}E[\sup_{\theta\in\Theta}|\tilde{I}_{n,t}^{(\mathcal
{E}%
)}(\theta) - I_{n,t}^{(\mathcal{E})}(\theta)|] %
\leq Kn^{-1}\sum_{t=1}^{n}\rho^{t} + \mathrm{o}(1) = \mathrm{o}(1)$.
\end{pf*}

\begin{pf*}{Proof of Lemma~\ref{lm:mom_est}} Define $\mathfrak{Z}%
_{n,t}(\theta) := \mathcal{E}_{t}^{2}(\theta)I_{n,t}^{(\mathcal{E}%
)}(\theta)/E[\mathcal{E}_{t}^{2}(\theta)I_{n,t}^{(\mathcal
{E})}(\theta)]$%
. By the same arguments used to prove approximation Lemma~\ref
{lm:approx}: $%
1/n\sum_{t=1}^{n}\mathcal{E}_{t}^{2}(\hat{\theta}_{n})\hat{I}_{n,t}^{(%
\mathcal{E})}(\hat{\theta}_{n}) = 1/n\times\break \sum_{t=1}^{n}\mathcal{E}_{t}^{2}(%
\hat{\theta}_{n})\*I_{n,t}^{(\mathcal{E})}(\hat{\theta}_{n})(1 + %
\mathrm{o}_{p} ( 1 ) )$. Since $\mathfrak{Z}_{n,t}(\theta)$ is uniformly
integrable and geometrically $\beta$-mixing by Assumption~\ref{ass1}(d), it
follows $%
1/n\sum_{t=1}^{n}\mathfrak{Z}_{n,t}(\theta) \stackrel{p}{\rightarrow} 1$
by Theorem~2 and Example~4 in Andrews \cite{Andrews88}. Moreover, since
$\mathfrak{Z}%
_{n,t}(\theta)$ is trivially $L_{1}$-bounded uniformly in $\Theta$, $%
\mathfrak{Z}_{n,t}(\theta)$ belongs to a separable Banach space, hence
$%
L_{1}$-bracketing numbers satisfy $N_{[\ ]}(\varepsilon,\Theta
,\|\cdot\|_{1}) < \infty$ (Dudley \cite{Dudley99}, Proposition~7.1.7). Combine the pointwise law and $N_{[\ ]}(\varepsilon
,\Theta,\|\cdot\|_{1}) < \infty$ to deduce $\sup_{\theta\in\Theta
}|1/n\sum_{t=1}^{n}\mathfrak{Z}_{n,t} ( \theta ) |
\stackrel{p}{\rightarrow} 0$ by Theorem~7.1.5 of Dudley \cite{Dudley99}.
Therefore $1/n\sum_{t=1}^{n}\mathcal{E}_{t}^{2}(\hat{\theta
}_{n})I_{n,t}^{(%
\mathcal{E})}(\hat{\theta}_{n})/E[\mathcal{E}_{t}^{2}(\hat{\theta}%
_{n})I_{n,t}^{(\mathcal{E})}(\hat{\theta}_{n})] \stackrel{p}{%
\rightarrow} 1$. %%
Further, by the definition of a derivative: $|E[\mathcal
{E}_{t}^{2}(\theta
)I_{n,t}^{(\mathcal{E})}(\theta)] - E[\mathcal{E}_{t}^{2}I_{n,t}^{(%
\mathcal{E})}]| \leq \|(\partial/\partial\theta)E[\mathcal{E}%
_{t}^{2}(\theta)I_{n,t}^{(\mathcal{E})}(\theta)]|_{\theta
^{0}}\|\times
\|\theta - \theta^{0}\| \times (1 + \mathrm{o}(1))$. By the same
argument as Lemma~\ref{lm:jac}(c) we can write
\begin{eqnarray*}
\frac{\partial}{\partial\theta}E \bigl[ \mathcal{E}_{t}^{2}(\theta
)I_{n,t}^{(\mathcal{E})}(\theta) \bigr] |_{\theta^{0}}&=&E \biggl[
\frac{%
\partial}{\partial\theta}\mathcal{E}_{t}^{2}(\theta)\bigg|_{\theta
^{0}}
\times I_{n,t}^{(\mathcal{E})} \biggr] \times \bigl( 1+\mathrm{o}(1) \bigr)
\\
&=&-2E%
\bigl[ \mathcal{E}_{t}\epsilon_{t}^{2}I_{n,t}^{(\mathcal{E})}
\mathfrak {s}%
_{t} \bigr] \times \bigl( 1+\mathrm{o}(1) \bigr) ,
\end{eqnarray*}
and trivially $E[\mathcal{E}_{t}\epsilon_{t}^{2}I_{n,t}^{(\mathcal{E})}
\mathfrak{s}_{t}] = E[\mathcal{E}_{t}^{2}\mathfrak{s}_{t}I_{n,t}^{(%
\mathcal{E})}] - E[\mathcal{E}_{t}\mathfrak{s}_{t}I_{n,t}^{(\mathcal{E}%
)}] = E[\mathcal{E}_{t}^{2}\mathfrak{s}_{t}I_{n,t}^{(\mathcal{E})}] =
E[\mathcal{E}_{t}^{2}I_{n,t}^{(\mathcal{E})}] \times E[\mathfrak{s}%
_{t}]$. Therefore $|E[\mathcal{E}_{t}^{2}(\theta)I_{n,t}^{(\mathcal{E}%
)}(\theta)] - E[\mathcal{E}_{t}^{2}I_{n,t}^{(\mathcal{E})}]| \leq %
K|E[\mathcal{E}_{t}^{2}I_{n,t}^{(\mathcal{E})}]|\times\|\theta - %
\theta^{0}\| \times (1 + \mathrm{o}(1))$. Now use $\hat{\theta}_{n} %
\stackrel{p}{\rightarrow} \theta^{0}$ by Theorem~\ref{th:consist} and
$%
\inf_{n\geq N}E[\mathcal{E}_{t}^{2}I_{n,t}^{(\mathcal{E})}] > 0$ for
some $N \geq 1$ to deduce $E[\mathcal{E}_{t}^{2}(\hat{\theta}%
_{n})\*I_{n,t}^{(\mathcal{E})}(\hat{\theta}_{n})]/E[\mathcal{E}%
_{t}^{2}I_{n,t}^{(\mathcal{E})}] \rightarrow 1$. This proves $%
1/n\sum_{t=1}^{n}\mathcal{E}_{t}^{2}(\hat{\theta}_{n})I_{n,t}^{(\mathcal
{E}%
)}(\hat{\theta}_{n})/E[\mathcal{E}_{t}^{2}I_{n,t}^{(\mathcal{E})}]
\stackrel{p}{\rightarrow} 1$.
\end{pf*}

%\section{}
\end{appendix}

\section*{Acknowledgements}
We thank two referees, an associate editor and editor Richard Davis for
helpful comments. We also thank Bruno Remillard for pointing out several
typos.

\begin{supplement}%[id=suppA]
\stitle{Supplement to ``Robust estimation and inference for heavy tailed GARCH''}
\slink[doi,text={10.3150/ 14-BEJ616SUPP}]{10.3150/14-BEJ616SUPP} %[doi,text={...}] - jei reikia suskaldyti doi
\sdatatype{.pdf}
\sfilename{BEJ616\_supp.pdf}
\sdescription{We prove Lemmas \ref{lm:approx},   \ref{lm:ulln},
\ref{lm:expan} and \ref{lm:jac}, and Lemmas
\ref{lB.1} and \ref{lB.2}.
Assume all functions satisfy Pollard's \cite{pollard84} permissibility
criteria, the measure space that governs all random variables
in this paper is complete, and therefore all majorants are
measurable. Cf. Dudley \cite{Dudley78}. Probability statements are
therefore with respect to outer probability, and expectations
over majorants are outer expectations.}
\end{supplement}

% imsref loaded by arune.pranskunaite, 2014-04-23 10:40:28
%

\printhistory


\begin{thebibliography}{61}
% pybtex-1.07. Style name=bej, version=1.4, label_style=nolabel,
%sorting_style=complex, cfg=None, language=None.

%b1 ###
\bibitem{Andrewsetal72}
%
\begin{bbook}[mr]
\bauthor{\bsnm{Andrews},~\bfnm{D.~F.}\binits{D.F.}},
\bauthor{\bsnm{Bickel},~\bfnm{P.~J.}\binits{P.J.}},
\bauthor{\bsnm{Hampel},~\bfnm{F.~R.}\binits{F.R.}},
\bauthor{\bsnm{Huber},~\bfnm{P.~J.}\binits{P.J.}},
\bauthor{\bsnm{Rogers},~\bfnm{W.~H.}\binits{W.H.}} \AND
\bauthor{\bsnm{Tukey},~\bfnm{J.~W.}\binits{J.W.}}
(\byear{1972}).
\btitle{Robust Estimates of Location: {S}urvey and Advances}.
\blocation{Princeton, NJ}:
\bpublisher{Princeton Univ. Press}.
\bid{mr={0331595}}
\end{bbook}
%
\bptok{imsref}%
% NOT OUTPUTED:
% fpage = ix+373
\endbibitem

%b2 ###
\bibitem{Andrews88}
%
\begin{barticle}[mr]
\bauthor{\bsnm{Andrews},~\bfnm{Donald~W.~K.}\binits{D.W.K.}}
(\byear{1988}).
\btitle{Laws of large numbers for dependent nonidentically distributed
random variables}.
\bjournal{Econometric Theory}
\bvolume{4}
\bpages{458--467}.
\bid{doi={10.1017/S0266466600013396}, issn={0266-4666}, mr={0985156}}
\end{barticle}
%
\bptok{imsref}%
% NOT OUTPUTED:
% issn = 0266-4666
% url = http://dx.doi.org/10.1017/S0266466600013396
% number = 3
% fjournal = Econometric Theory
\endbibitem

%b3 ###
\bibitem{Andrews99}
%
\begin{barticle}[mr]
\bauthor{\bsnm{Andrews},~\bfnm{Donald~W.~K.}\binits{D.W.K.}}
(\byear{1999}).
\btitle{Estimation when a parameter is on a boundary}.
\bjournal{Econometrica}
\bvolume{67}
\bpages{1341--1383}.
\bid{doi={10.1111/1468-0262.00082}, issn={0012-9682}, mr={1720781}}
\end{barticle}
%
\bptok{imsref}%
% NOT OUTPUTED:
% issn = 0012-9682
% url = http://dx.doi.org/10.1111/1468-0262.00082
% number = 6
% coden = ECMTA7
% fjournal = Econometrica. Journal of the Econometric Society
\endbibitem

%b4 ###
\bibitem{Basrak2002}
%
\begin{barticle}[mr]
\bauthor{\bsnm{Basrak},~\bfnm{Bojan}\binits{B.}},
\bauthor{\bsnm{Davis},~\bfnm{Richard~A.}\binits{R.A.}} \AND
\bauthor{\bsnm{Mikosch},~\bfnm{Thomas}\binits{T.}}
(\byear{2002}).
\btitle{Regular variation of {GARCH} processes}.
\bjournal{Stochastic Process. Appl.}
\bvolume{99}
\bpages{95--115}.
\bid{doi={10.1016/S0304-4149(01)00156-9}, issn={0304-4149}, mr={1894253}}
\end{barticle}
%
\bptok{imsref}%
% NOT OUTPUTED:
% issn = 0304-4149
% url = http://dx.doi.org/10.1016/S0304-4149(01)00156-9
% number = 1
% coden = STOPB7
% fjournal = Stochastic Processes and their Applications
\endbibitem

%b5 ###
\bibitem{BerkesHorvath04}
%
\begin{barticle}[mr]
\bauthor{\bsnm{Berkes},~\bfnm{Istv{\'a}n}\binits{I.}} \AND
\bauthor{\bsnm{Horv{\'a}th},~\bfnm{Lajos}\binits{L.}}
(\byear{2004}).
\btitle{The efficiency of the estimators of the parameters in {GARCH}
processes}.
\bjournal{Ann. Statist.}
\bvolume{32}
\bpages{633--655}.
\bid{doi={10.1214/009053604000000120}, issn={0090-5364}, mr={2060172}}
\end{barticle}
%
\bptok{imsref}%
% NOT OUTPUTED:
% issn = 0090-5364
% url = http://dx.doi.org/10.1214/009053604000000120
% number = 2
% coden = ASTSC7
% fjournal = The Annals of Statistics
\endbibitem

%b6 ###
\bibitem{Berkesal03}
%
\begin{barticle}[mr]
\bauthor{\bsnm{Berkes},~\bfnm{Istv{\'a}n}\binits{I.}},
\bauthor{\bsnm{Horv{\'a}th},~\bfnm{Lajos}\binits{L.}} \AND
\bauthor{\bsnm{Kokoszka},~\bfnm{Piotr}\binits{P.}}
(\byear{2003}).
\btitle{G{ARCH} processes: Structure and estimation}.
\bjournal{Bernoulli}
\bvolume{9}
\bpages{201--227}.
\bid{doi={10.3150/bj/1068128975}, issn={1350-7265}, mr={1997027}}
\end{barticle}
%
\bptok{imsref}%
% NOT OUTPUTED:
% issn = 1350-7265
% url = http://dx.doi.org/10.3150/bj/1068128975
% number = 2
% fjournal = Bernoulli. Official Journal of the Bernoulli Society for
%Mathematical Statistics and Probability
\endbibitem

%b7 ###
\bibitem{bollerslev86}
%
\begin{barticle}[mr]
\bauthor{\bsnm{Bollerslev},~\bfnm{Tim}\binits{T.}}
(\byear{1986}).
\btitle{Generalized autoregressive conditional heteroskedasticity}.
\bjournal{J. Econometrics}
\bvolume{31}
\bpages{307--327}.
\bid{doi={10.1016/0304-4076(86)90063-1}, issn={0304-4076}, mr={0853051}}
\end{barticle}
%
\bptok{imsref}%
% NOT OUTPUTED:
% issn = 0304-4076
% url = http://dx.doi.org/10.1016/0304-4076(86)90063-1
% number = 3
% coden = JECMB6
% fjournal = Journal of Econometrics
\endbibitem

%b8 ###
\bibitem{Boudtetal2011}
%
\begin{barticle}[auto:STB|2014/02/12|14:17:21]
\bauthor{\bsnm{Boudt},~\bfnm{K.}\binits{K.}},
\bauthor{\bsnm{Danielsson},~\bfnm{J.}\binits{J.}} \AND
\bauthor{\bsnm{Laurent},~\bfnm{S.}\binits{S.}}
(\byear{2011}).
\btitle{Robust forecasting of dynamic conditional correlation garch models}.
\bjournal{Internat. J. Forecasting}
\bvolume{29}
\bpages{244--257}.
\end{barticle}
%
\bptok{imsref}%
\endbibitem

%b9 ###
\bibitem{Bradley05}
%
\begin{barticle}[mr]
\bauthor{\bsnm{Bradley},~\bfnm{Richard~C.}\binits{R.C.}}
(\byear{2005}).
\btitle{Basic properties of strong mixing conditions. {A} survey and
some open questions}.
\bjournal{Probab. Surv.}
\bvolume{2}
\bpages{107--144}.
\bnote{Update of, and a supplement to, the 1986 original}.
\bid{doi={10.1214/154957805100000104}, issn={1549-5787}, mr={2178042}}
\end{barticle}
%
\bptok{imsref}%
% NOT OUTPUTED:
% issn = 1549-5787
% url = http://dx.doi.org/10.1214/154957805100000104
% fjournal = Probability Surveys
\endbibitem

%b10 ###
\bibitem{Breiman65}
%
\begin{barticle}[mr]
\bauthor{\bsnm{Breiman},~\bfnm{L.}\binits{L.}}
(\byear{1965}).
\btitle{On some limit theorems similar to the arc-sin law}.
\bjournal{Teor. Veroyatn. Primen.}
\bvolume{10}
\bpages{351--360}.
\bid{issn={0040-361X}, mr={0184274}}
\end{barticle}
%
\bptok{imsref}%
% NOT OUTPUTED:
% issn = 0040-361x
% fjournal = Akademija Nauk SSSR. Teorija Verojatnoste\u\i\ i ee
%Primenenija
\endbibitem

%b11 ###
\bibitem{CantoniRonchetti}
%
\begin{barticle}[mr]
\bauthor{\bsnm{Cantoni},~\bfnm{Eva}\binits{E.}} \AND
\bauthor{\bsnm{Ronchetti},~\bfnm{Elvezio}\binits{E.}}
(\byear{2001}).
\btitle{Robust inference for generalized linear models}.
\bjournal{J. Amer. Statist. Assoc.}
\bvolume{96}
\bpages{1022--1030}.
\bid{doi={10.1198/016214501753209004}, issn={0162-1459}, mr={1947250}}
\end{barticle}
%
\bptok{imsref}%
% NOT OUTPUTED:
% issn = 0162-1459
% url = http://dx.doi.org/10.1198/016214501753209004
% number = 455
% coden = JSTNAL
% fjournal = Journal of the American Statistical Association
\endbibitem

%b12 ###
\bibitem{carrascochen02}
%
\begin{barticle}[mr]
\bauthor{\bsnm{Carrasco},~\bfnm{Marine}\binits{M.}} \AND
\bauthor{\bsnm{Chen},~\bfnm{Xiaohong}\binits{X.}}
(\byear{2002}).
\btitle{Mixing and moment properties of various {GARCH} and stochastic
volatility models}.
\bjournal{Econometric Theory}
\bvolume{18}
\bpages{17--39}.
\bid{doi={10.1017/S0266466602181023}, issn={0266-4666}, mr={1885348}}
\end{barticle}
%
\bptok{imsref}%
% NOT OUTPUTED:
% issn = 0266-4666
% url = http://dx.doi.org/10.1017/S0266466602181023
% number = 1
% fjournal = Econometric Theory
\endbibitem

%b13 ###
\bibitem{CavaliereGeorgiev09}
%
\begin{barticle}[mr]
\bauthor{\bsnm{Cavaliere},~\bfnm{Giuseppe}\binits{G.}} \AND
\bauthor{\bsnm{Georgiev},~\bfnm{Iliyan}\binits{I.}}
(\byear{2009}).
\btitle{Robust inference in autoregressions with multiple outliers}.
\bjournal{Econometric Theory}
\bvolume{25}
\bpages{1625--1661}.
\bid{doi={10.1017/S0266466609990272}, issn={0266-4666}, mr={2557576}}
\end{barticle}
%
\bptok{imsref}%
% NOT OUTPUTED:
% issn = 0266-4666
% url = http://dx.doi.org/10.1017/S0266466609990272
% number = 6
% fjournal = Econometric Theory
\endbibitem

%b14 ###
\bibitem{CharlesDarne}
%
\begin{barticle}[mr]
\bauthor{\bsnm{Charles},~\bfnm{Am{\'e}lie}\binits{A.}} \AND
\bauthor{\bsnm{Darn{\'e}},~\bfnm{Olivier}\binits{O.}}
(\byear{2005}).
\btitle{Outliers and {GARCH} models in financial data}.
\bjournal{Econom. Lett.}
\bvolume{86}
\bpages{347--352}.
\bid{doi={10.1016/j.econlet.2004.07.019}, issn={0165-1765}, mr={2124418}}
\end{barticle}
%
\bptok{imsref}%
% NOT OUTPUTED:
% issn = 0165-1765
% url = http://dx.doi.org/10.1016/j.econlet.2004.07.019
% number = 3
% coden = ECLEDS
% fjournal = Economics Letters
\endbibitem

%b15 ###
\bibitem{ChenLiu93}
%
\begin{barticle}[auto:STB|2014/02/12|14:17:21]
\bauthor{\bsnm{Chen},~\bfnm{C.}\binits{C.}} \AND
\bauthor{\bsnm{Liu},~\bfnm{L.-M.}\binits{L.-M.}}
(\byear{1993}).
\btitle{Joint estimation of model parameters and outlier effects in
time series}.
\bjournal{J. Amer. Statist. Assoc.}
\bvolume{88}
\bpages{284--297}.
\end{barticle}
%
\bptok{imsref}%
\endbibitem

%b16 ###
\bibitem{cizek08}
%
\begin{barticle}[mr]
\bauthor{\bsnm{{\v{C}}{\'{\i}}{\v{z}}ek},~\bfnm{Pavel}\binits{P.}}
(\byear{2008}).
\btitle{General trimmed estimation: Robust approach to nonlinear and
limited dependent variable models}.
\bjournal{Econometric Theory}
\bvolume{24}
\bpages{1500--1529}.
\bid{doi={10.1017/S0266466608080596}, issn={0266-4666}, mr={2456536}}
\end{barticle}
%
\bptok{imsref}%
% NOT OUTPUTED:
% issn = 0266-4666
% url = http://dx.doi.org/10.1017/S0266466608080596
% number = 6
% fjournal = Econometric Theory
\endbibitem

%b17 ###
\bibitem{Davis2010}
%
\begin{bincollection}[auto:STB|2014/02/12|14:17:21]
\bauthor{\bsnm{Davis},~\bfnm{R.~A.}\binits{R.A.}}
(\byear{2010}).
\btitle{Heavy tails in financial time series}.
In \bbooktitle{Encyclopedia of Quantitative Finance}
(\beditor{\bfnm{R.}\binits{R.}~\bsnm{Cont}}, ed.).
\blocation{New York}:
\bpublisher{Wiley}.
\end{bincollection}
%
\bptok{imsref}%
\endbibitem

%b18 ###
\bibitem{Mendes00}
%
\begin{barticle}[mr]
\bauthor{\bsnm{De Melo Mendes},~\bfnm{Beatriz~Vaz}\binits{B.V.}}
(\byear{2000}).
\btitle{Assessing the bias of maximum likelihood estimates of
contaminated {G}arch models}.
\bjournal{J. Stat. Comput. Simul.}
\bvolume{67}
\bpages{359--376}.
\bid{doi={10.1080/00949650008812051}, issn={0094-9655}, mr={1806901}}
\end{barticle}
%
\bptok{imsref}%
% NOT OUTPUTED:
% issn = 0094-9655
% url = http://dx.doi.org/10.1080/00949650008812051
% number = 4
% coden = JSCSAT
% fjournal = Journal of Statistical Computation and Simulation
\endbibitem

%b19 ###
\bibitem{Dudley78}
%
\begin{barticle}[mr]
\bauthor{\bsnm{Dudley},~\bfnm{R.~M.}\binits{R.M.}}
(\byear{1978}).
\btitle{Central limit theorems for empirical measures}.
\bjournal{Ann. Probab.}
\bvolume{6}
\bpages{899--929 (1979)}.
\bid{issn={0091-1798}, mr={0512411}}
\end{barticle}
%
\bptok{imsref}%
% NOT OUTPUTED:
% issn = 0091-1798
% url =
%http://links.jstor.org/sici?sici=0091-1798(197812)6:6<899:CLTFEM>2.0.CO;2-O&origin=MSN
% number = 6
% coden = APBYAE
% fjournal = The Annals of Probability
\endbibitem

%b20 ###
\bibitem{Dudley99}
%
\begin{bbook}[mr]
\bauthor{\bsnm{Dudley},~\bfnm{R.~M.}\binits{R.M.}}
(\byear{1999}).
\btitle{Uniform Central Limit Theorems}.
\bseries{Cambridge Studies in Advanced Mathematics}
\bvolume{63}.
\blocation{Cambridge}:
\bpublisher{Cambridge Univ. Press}.
\bid{doi={10.1017/CBO9780511665622}, mr={1720712}}
\end{bbook}
%
\bptok{imsref}%
% NOT OUTPUTED:
% isbn = 0-521-46102-2
% url = http://dx.doi.org/10.1017/CBO9780511665622
% fpage = xiv+436
\endbibitem

%b21 ###
\bibitem{Embetal97}
%
\begin{bbook}[mr]
\bauthor{\bsnm{Embrechts},~\bfnm{Paul}\binits{P.}},
\bauthor{\bsnm{Kl{\"u}ppelberg},~\bfnm{Claudia}\binits{C.}} \AND
\bauthor{\bsnm{Mikosch},~\bfnm{Thomas}\binits{T.}}
(\byear{1997}).
\btitle{Modelling Extremal Events. For Insurance and Finance}.
\bseries{Applications of Mathematics (New York)}
\bvolume{33}.
\blocation{Berlin}:
\bpublisher{Springer}.
\bid{mr={1458613}}
\end{bbook}
%
\bptok{imsref}%
% NOT OUTPUTED:
% isbn = 3-540-60931-8
% fpage = xvi+645
\endbibitem

%b22 ###
\bibitem{engleng93}
%
\begin{barticle}[auto:STB|2014/02/12|14:17:21]
\bauthor{\bsnm{Engle},~\bfnm{R.~F.}\binits{R.F.}} \AND
\bauthor{\bsnm{Ng},~\bfnm{V.~K.}\binits{V.K.}}
(\byear{1993}).
\btitle{Measuring and testing the impact of news on volatility}.
\bjournal{J. Finance}
\bvolume{48}
\bpages{1749--1778}.
\end{barticle}
%
\bptok{imsref}%
\endbibitem

%b23 ###
\bibitem{Fanetal12}
%
\begin{bmisc}[auto:STB|2014/02/12|14:17:21]
\bauthor{\bsnm{Fan},~\bfnm{J.}\binits{J.}},
\bauthor{\bsnm{Qi},~\bfnm{L.}\binits{L.}} \AND
\bauthor{\bsnm{Xiu},~\bfnm{D.}\binits{D.}}
(\byear{2014}).
\bhowpublished{Quasi maximum likelihood estimation of GARCH models with
heavy-tailed likelihoods. \textit{J. Business and Economic Statistics}. To appear}.
\end{bmisc}
%
\bptok{imsref}%
% NOT OUTPUTED:
% sortkey = Fan(2012
\endbibitem

%b24 ###
\bibitem{FZ04}
%
\begin{barticle}[mr]
\bauthor{\bsnm{Francq},~\bfnm{Christian}\binits{C.}} \AND
\bauthor{\bsnm{Zako{\"{\i}}an},~\bfnm{Jean-Michel}\binits{J.-M.}}
(\byear{2004}).
\btitle{Maximum likelihood estimation of pure {GARCH} and
{ARMA}--{GARCH} processes}.
\bjournal{Bernoulli}
\bvolume{10}
\bpages{605--637}.
\bid{doi={10.3150/bj/1093265632}, issn={1350-7265}, mr={2076065}}
\end{barticle}
%
\bptok{imsref}%
% NOT OUTPUTED:
% issn = 1350-7265
% url = http://dx.doi.org/10.3150/bj/1093265632
% number = 4
% fjournal = Bernoulli. Official Journal of the Bernoulli Society for
%Mathematical Statistics and Probability
\endbibitem

%b25 ###
\bibitem{FZ10}
%
\begin{bbook}[auto:STB|2014/02/12|14:17:21]
\bauthor{\bsnm{Francq},~\bfnm{C.}\binits{C.}} \AND
\bauthor{\bsnm{Zako{\"{\i}}an},~\bfnm{J.-M.}\binits{J.-M.}}
(\byear{2010}).
\btitle{GARCH Models: Structure, Statistical Inference and Financial
Applications}.
\blocation{New York}:
\bpublisher{Wiley}.
\end{bbook}
%
\bptok{imsref}%
\endbibitem

%b26 ###
\bibitem{hallyao03}
%
\begin{barticle}[mr]
\bauthor{\bsnm{Hall},~\bfnm{Peter}\binits{P.}} \AND
\bauthor{\bsnm{Yao},~\bfnm{Qiwei}\binits{Q.}}
(\byear{2003}).
\btitle{Inference in {ARCH} and {GARCH} models with heavy-tailed errors}.
\bjournal{Econometrica}
\bvolume{71}
\bpages{285--317}.
\bid{doi={10.1111/1468-0262.00396}, issn={0012-9682}, mr={1956860}}
\end{barticle}
%
\bptok{imsref}%
% NOT OUTPUTED:
% issn = 0012-9682
% url = http://dx.doi.org/10.1111/1468-0262.00396
% number = 1
% coden = ECMTA7
% fjournal = Econometrica. Journal of the Econometric Society
\endbibitem

%b27 ###
\bibitem{Hampleetal86}
%
\begin{bbook}[mr]
\bauthor{\bsnm{Hampel},~\bfnm{Frank~R.}\binits{F.R.}},
\bauthor{\bsnm{Ronchetti},~\bfnm{Elvezio~M.}\binits{E.M.}},
\bauthor{\bsnm{Rousseeuw},~\bfnm{Peter~J.}\binits{P.J.}} \AND
\bauthor{\bsnm{Stahel},~\bfnm{Werner~A.}\binits{W.A.}}
(\byear{1986}).
\btitle{Robust Statistics. The Approach Based on Influence Functions}.
\bseries{Wiley Series in Probability and Mathematical Statistics:
Probability and Mathematical Statistics}.
\blocation{New York}:
\bpublisher{Wiley}.
%\bnote{}.
\bid{mr={0829458}}
\end{bbook}
%
\bptok{imsref}%
% NOT OUTPUTED:
% isbn = 0-471-82921-8
% fpage = xxiv+502
\endbibitem

%b28 ###
\bibitem{hansen82}
%
\begin{barticle}[mr]
\bauthor{\bsnm{Hansen},~\bfnm{Lars~Peter}\binits{L.P.}}
(\byear{1982}).
\btitle{Large sample properties of generalized method of moments estimators}.
\bjournal{Econometrica}
\bvolume{50}
\bpages{1029--1054}.
\bid{doi={10.2307/1912775}, issn={0012-9682}, mr={0666123}}
\end{barticle}
%
\bptok{imsref}%
% NOT OUTPUTED:
% issn = 0012-9682
% url = http://dx.doi.org/10.2307/1912775
% number = 4
% coden = ECMTA7
% fjournal = Econometrica. Journal of the Econometric Society
\endbibitem

%b29 ###
\bibitem{Hill10}
%
\begin{barticle}[mr]
\bauthor{\bsnm{Hill},~\bfnm{Jonathan~B.}\binits{J.B.}}
(\byear{2010}).
\btitle{On tail index estimation for dependent, heterogeneous data}.
\bjournal{Econometric Theory}
\bvolume{26}
\bpages{1398--1436}.
\bid{doi={10.1017/S0266466609990624}, issn={0266-4666}, mr={2684790}}
\end{barticle}
%
\bptok{imsref}%
% NOT OUTPUTED:
% issn = 0266-4666
% url = http://dx.doi.org/10.1017/S0266466609990624
% number = 5
% fjournal = Econometric Theory
\endbibitem

%b30 ###
\bibitem{Hilltest}
%
\begin{bincollection}[auto:STB|2014/02/12|14:17:21]
\bauthor{\bsnm{Hill},~\bfnm{J.~B.}\binits{J.B.}}
(\byear{2012}).
\btitle{Heavy-tail and plug-in robust consistent conditional moment
tests of functional form}.
In \bbooktitle{Festschrift in Honor of Hal White}
(\beditor{\bfnm{X.}\binits{X.}~\bsnm{Chen}} \AND
\beditor{\bfnm{N.}\binits{N.}~\bsnm{Swanson}}, eds.)
\bpages{241--274}.
\blocation{New York}:
\bpublisher{Springer}.
\end{bincollection}
%
\bptok{imsref}%
\endbibitem

%b31 ###
\bibitem{Hillltts}
%
\begin{barticle}[mr]
\bauthor{\bsnm{Hill},~\bfnm{Jonathan~B.}\binits{J.B.}}
(\byear{2013}).
\btitle{Least tail-trimmed squares for infinite variance autoregressions}.
\bjournal{J. Time Series Anal.}
\bvolume{34}
\bpages{168--186}.
\bid{doi={10.1111/jtsa.12005}, issn={0143-9782}, mr={3028364}}
\end{barticle}
%
\bptok{imsref}%
% NOT OUTPUTED:
% issn = 0143-9782
% url = http://dx.doi.org/10.1111/jtsa.12005
% number = 2
% fjournal = Journal of Time Series Analysis
\endbibitem

%b32 ###
\bibitem{Hillsupp}
%
\begin{bmisc}[auto:STB|2014/02/12|14:17:21]
\bauthor{\bsnm{Hill},~\bfnm{J.~B.}\binits{J.B.}}
(\byear{2014}).
\bhowpublished{Supplement to ``Robust estimation and
inference for heavy tailed GARCH.'' DOI:\doiurl{10.3150/14-BEJ616SUPP}}.
\end{bmisc}
%
\bptok{imsref}%
\endbibitem

%b33 ###
\bibitem{HillAguilar13}
%
\begin{barticle}[mr]
\bauthor{\bsnm{Hill},~\bfnm{Jonathan~B.}\binits{J.B.}} \AND
\bauthor{\bsnm{Aguilar},~\bfnm{Mike}\binits{M.}}
(\byear{2013}).
\btitle{Moment condition tests for heavy tailed time series}.
\bjournal{J. Econometrics}
\bvolume{172}
\bpages{255--274}.
\bid{doi={10.1016/j.jeconom.2012.08.013}, issn={0304-4076}, mr={3010616}}
\end{barticle}
%
\bptok{imsref}%
% NOT OUTPUTED:
% issn = 0304-4076
% url = http://dx.doi.org/10.1016/j.jeconom.2012.08.013
% number = 2
% coden = JECMB6
% fjournal = Journal of Econometrics
\endbibitem

%b34 ###
\bibitem{Huber1964}
%
\begin{barticle}[mr]
\bauthor{\bsnm{Huber},~\bfnm{Peter~J.}\binits{P.J.}}
(\byear{1964}).
\btitle{Robust estimation of a location parameter}.
\bjournal{Ann. Math. Statist.}
\bvolume{35}
\bpages{73--101}.
\bid{issn={0003-4851}, mr={0161415}}
\end{barticle}
%
\bptok{imsref}%
% NOT OUTPUTED:
% issn = 0003-4851
% fjournal = Annals of Mathematical Statistics
\endbibitem

%b61 ###
\bibitem{HoffJorg1984}
%
\begin{bmisc}[auto:STB|2014/02/12|14:17:21]
\bauthor{\bsnm{Hoffman-J{\o}rgensen},~\bfnm{J.}\binits{J.}}
(\byear{1984}).
\bhowpublished{Convergence of stochastic processes on Polish spaces.
Mimeo, Aarhus Univ., Denmark}.
\end{bmisc}
%
\bptok{imsref}%
% NOT OUTPUTED:
% sortkey = Hoffman(1984
\endbibitem

%b35 ###
\bibitem{Ibrag1962}
%
\begin{barticle}[mr]
\bauthor{\bsnm{Ibragimov},~\bfnm{I.~A.}\binits{I.A.}}
(\byear{1962}).
\btitle{Some limit theorems for stationary processes}.
\bjournal{Teor. Veroyatn. Primen.}
\bvolume{7}
\bpages{361--392}.
\bid{issn={0040-361X}, mr={0148125}}
\end{barticle}
%
\bptok{imsref}%
% NOT OUTPUTED:
% issn = 0040-361x
% fjournal = Akademija Nauk SSSR. Teorija Verojatnoste\u\i\ i ee
%Primenenija
\endbibitem

%b36 ###
\bibitem{JurSen96}
%
\begin{bbook}[mr]
\bauthor{\bsnm{Jure{\v{c}}kov{\'a}},~\bfnm{Jana}\binits{J.}} \AND
\bauthor{\bsnm{Sen},~\bfnm{Pranab~Kumar}\binits{P.K.}}
(\byear{1996}).
\btitle{Robust Statistical Procedures. Asymptotics and Interrelations}.
\bseries{Wiley Series in Probability and Statistics: Applied
Probability and Statistics}.
\blocation{New York}:
\bpublisher{Wiley}.
\bid{mr={1387346}}
\end{bbook}
%
\bptok{imsref}%
% NOT OUTPUTED:
% isbn = 0-471-82221-3
% fpage = xvi+466
\endbibitem

%b37 ###
\bibitem{leadbetteretal83}
%
\begin{bbook}[mr]
\bauthor{\bsnm{Leadbetter},~\bfnm{M.~R.}\binits{M.R.}},
\bauthor{\bsnm{Lindgren},~\bfnm{Georg}\binits{G.}} \AND
\bauthor{\bsnm{Rootz{\'e}n},~\bfnm{Holger}\binits{H.}}
(\byear{1983}).
\btitle{Extremes and Related Properties of Random Sequences and Processes}.
\bseries{Springer Series in Statistics}.
\blocation{New York}:
\bpublisher{Springer}.
\bid{mr={0691492}}
\end{bbook}
%
\bptok{imsref}%
% NOT OUTPUTED:
% isbn = 0-387-90731-9
% fpage = xii+336
\endbibitem

%b38 ###
\bibitem{leehansen94}
%
\begin{barticle}[mr]
\bauthor{\bsnm{Lee},~\bfnm{Sang-Won}\binits{S.-W.}} \AND
\bauthor{\bsnm{Hansen},~\bfnm{Bruce~E.}\binits{B.E.}}
(\byear{1994}).
\btitle{Asymptotic theory for the $\operatorname{GARCH}(1,1)$
quasi-maximum likelihood estimator}.
\bjournal{Econometric Theory}
\bvolume{10}
\bpages{29--52}.
\bid{doi={10.1017/S0266466600008215}, issn={0266-4666}, mr={1279689}}
\end{barticle}
%
\bptok{imsref}%
% NOT OUTPUTED:
% issn = 0266-4666
% url = http://dx.doi.org/10.1017/S0266466600008215
% number = 1
% fjournal = Econometric Theory
\endbibitem

%b39 ###
\bibitem{Lighthill}
%
\begin{bbook}[mr]
\bauthor{\bsnm{Lighthill},~\bfnm{M.~J.}\binits{M.J.}}
(\byear{1958}).
\btitle{Introduction to {F}ourier Analysis and Generalised Functions}.
\bseries{Cambridge Monographs on Mechanics and Applied Mathematics}.
\blocation{Cambridge}:
\bpublisher{Cambridge Univ. Press}.
\bid{mr={0092119}}
\end{bbook}
%
\bptok{imsref}%
% NOT OUTPUTED:
% fpage = viii+79
\endbibitem

%b40 ###
\bibitem{ling07}
%
\begin{barticle}[mr]
\bauthor{\bsnm{Ling},~\bfnm{Shiqing}\binits{S.}}
(\byear{2007}).
\btitle{Self-weighted and local quasi-maximum likelihood estimators for
{ARMA}--{GARCH}/{IGARCH} models}.
\bjournal{J. Econometrics}
\bvolume{140}
\bpages{849--873}.
\bid{doi={10.1016/j.jeconom.2006.07.016}, issn={0304-4076}, mr={2408929}}
\end{barticle}
%
\bptok{imsref}%
% NOT OUTPUTED:
% issn = 0304-4076
% url = http://dx.doi.org/10.1016/j.jeconom.2006.07.016
% number = 2
% coden = JECMB6
% fjournal = Journal of Econometrics
\endbibitem

%b41 ###
\bibitem{Liu2006}
%
\begin{barticle}[mr]
\bauthor{\bsnm{Liu},~\bfnm{Ji-Chun}\binits{J.-C.}}
(\byear{2006}).
\btitle{On the tail behaviors of a family of {GARCH} processes}.
\bjournal{Econometric Theory}
\bvolume{22}
\bpages{852--862}.
\bid{doi={10.1017/S0266466606060397}, issn={0266-4666}, mr={2291220}}
\end{barticle}
%
\bptok{imsref}%
% NOT OUTPUTED:
% issn = 0266-4666
% url = http://dx.doi.org/10.1017/S0266466606060397
% number = 5
% fjournal = Econometric Theory
\endbibitem

%b42 ###
\bibitem{Mancinietal05}
%
\begin{barticle}[mr]
\bauthor{\bsnm{Mancini},~\bfnm{Loriano}\binits{L.}},
\bauthor{\bsnm{Ronchetti},~\bfnm{Elvezio}\binits{E.}} \AND
\bauthor{\bsnm{Trojani},~\bfnm{Fabio}\binits{F.}}
(\byear{2005}).
\btitle{Optimal conditionally unbiased bounded-influence inference in
dynamic location and scale models}.
\bjournal{J. Amer. Statist. Assoc.}
\bvolume{100}
\bpages{628--641}.
\bid{doi={10.1198/016214504000001402}, issn={0162-1459}, mr={2160565}}
\end{barticle}
%
\bptok{imsref}%
% NOT OUTPUTED:
% issn = 0162-1459
% url = http://dx.doi.org/10.1198/016214504000001402
% number = 470
% coden = JSTNAL
% fjournal = Journal of the American Statistical Association
\endbibitem

%b43 ###
\bibitem{meitzsaikkonen08}
%
\begin{barticle}[mr]
\bauthor{\bsnm{Meitz},~\bfnm{Mika}\binits{M.}} \AND
\bauthor{\bsnm{Saikkonen},~\bfnm{Pentti}\binits{P.}}
(\byear{2008}).
\btitle{Stability of nonlinear {AR}-{GARCH} models}.
\bjournal{J. Time Ser. Anal.}
\bvolume{29}
\bpages{453--475}.
\bid{doi={10.1111/j.1467-9892.2007.00562.x}, issn={0143-9782}, mr={2410184}}
\end{barticle}
%
\bptok{imsref}%
% NOT OUTPUTED:
% issn = 0143-9782
% url = http://dx.doi.org/10.1111/j.1467-9892.2007.00562.x
% number = 3
% fjournal = Journal of Time Series Analysis
\endbibitem

%b44 ###
\bibitem{meitzsaik11}
%
\begin{barticle}[mr]
\bauthor{\bsnm{Meitz},~\bfnm{Mika}\binits{M.}} \AND
\bauthor{\bsnm{Saikkonen},~\bfnm{Pentti}\binits{P.}}
(\byear{2011}).
\btitle{Parameter estimation in nonlinear {AR}-{GARCH} models}.
\bjournal{Econometric Theory}
\bvolume{27}
\bpages{1236--1278}.
\bid{doi={10.1017/S0266466611000041}, issn={0266-4666}, mr={2868839}}
\end{barticle}
%
\bptok{imsref}%
% NOT OUTPUTED:
% issn = 0266-4666
% url = http://dx.doi.org/10.1017/S0266466611000041
% number = 6
% fjournal = Econometric Theory
\endbibitem

%b45 ###
\bibitem{Mulleretal09}
%
\begin{barticle}[mr]
\bauthor{\bsnm{Muler},~\bfnm{Nora}\binits{N.}},
\bauthor{\bsnm{Pe{\~n}a},~\bfnm{Daniel}\binits{D.}} \AND
\bauthor{\bsnm{Yohai},~\bfnm{V{\'{\i}}ctor~J.}\binits{V.J.}}
(\byear{2009}).
\btitle{Robust estimation for {ARMA} models}.
\bjournal{Ann. Statist.}
\bvolume{37}
\bpages{816--840}.
\bid{doi={10.1214/07-AOS570}, issn={0090-5364}, mr={2502652}}
\end{barticle}
%
\bptok{imsref}%
% NOT OUTPUTED:
% issn = 0090-5364
% url = http://dx.doi.org/10.1214/07-AOS570
% number = 2
% coden = ASTSC7
% fjournal = The Annals of Statistics
\endbibitem

%b46 ###
\bibitem{MulerYohai08}
%
\begin{barticle}[mr]
\bauthor{\bsnm{Muler},~\bfnm{Nora}\binits{N.}} \AND
\bauthor{\bsnm{Yohai},~\bfnm{Victor~J.}\binits{V.J.}}
(\byear{2008}).
\btitle{Robust estimates for {GARCH} models}.
\bjournal{J. Statist. Plann. Inference}
\bvolume{138}
\bpages{2918--2940}.
\bid{doi={10.1016/j.jspi.2007.11.003}, issn={0378-3758}, mr={2442223}}
\end{barticle}
%
\bptok{imsref}%
% NOT OUTPUTED:
% issn = 0378-3758
% url = http://dx.doi.org/10.1016/j.jspi.2007.11.003
% number = 10
% coden = JSPIDN
% fjournal = Journal of Statistical Planning and Inference
\endbibitem

%b47 ###
\bibitem{Nelson90}
%
\begin{barticle}[mr]
\bauthor{\bsnm{Nelson},~\bfnm{Daniel~B.}\binits{D.B.}}
(\byear{1990}).
\btitle{Stationarity and persistence in the {GARCH{$(1,1)$}} model}.
\bjournal{Econometric Theory}
\bvolume{6}
\bpages{318--334}.
\bid{doi={10.1017/S0266466600005296}, issn={0266-4666}, mr={1085577}}
\end{barticle}
%
\bptok{imsref}%
% NOT OUTPUTED:
% issn = 0266-4666
% url = http://dx.doi.org/10.1017/S0266466600005296
% number = 3
% fjournal = Econometric Theory
\endbibitem

%b48 ###
\bibitem{NeweySteigerwald97}
%
\begin{barticle}[mr]
\bauthor{\bsnm{Newey},~\bfnm{Whitney~K.}\binits{W.K.}} \AND
\bauthor{\bsnm{Steigerwald},~\bfnm{Douglas~G.}\binits{D.G.}}
(\byear{1997}).
\btitle{Asymptotic bias for quasi-maximum-likelihood estimators in
conditional heteroskedasticity models}.
\bjournal{Econometrica}
\bvolume{65}
\bpages{587--599}.
\bid{doi={10.2307/2171754}, issn={0012-9682}, mr={1445623}}
\end{barticle}
%
\bptok{imsref}%
% NOT OUTPUTED:
% issn = 0012-9682
% url = http://dx.doi.org/10.2307/2171754
% number = 3
% coden = ECMTA7
% fjournal = Econometrica. Journal of the Econometric Society
\endbibitem

%b49 ###
\bibitem{pakespollard89}
%
\begin{barticle}[mr]
\bauthor{\bsnm{Pakes},~\bfnm{Ari{\'e}l}\binits{A.}} \AND
\bauthor{\bsnm{Pollard},~\bfnm{David}\binits{D.}}
(\byear{1989}).
\btitle{Simulation and the asymptotics of optimization estimators}.
\bjournal{Econometrica}
\bvolume{57}
\bpages{1027--1057}.
\bid{doi={10.2307/1913622}, issn={0012-9682}, mr={1014540}}
\end{barticle}
%
\bptok{imsref}%
% NOT OUTPUTED:
% issn = 0012-9682
% url = http://dx.doi.org/10.2307/1913622
% number = 5
% coden = ECMTA7
% fjournal = Econometrica. Journal of the Econometric Society
\endbibitem

%b50 ###
\bibitem{Peligrad96}
%
\begin{barticle}[mr]
\bauthor{\bsnm{Peligrad},~\bfnm{Magda}\binits{M.}}
(\byear{1996}).
\btitle{On the asymptotic normality of sequences of weak dependent
random variables}.
\bjournal{J. Theoret. Probab.}
\bvolume{9}
\bpages{703--715}.
\bid{doi={10.1007/BF02214083}, issn={0894-9840}, mr={1400595}}
\end{barticle}
%
\bptok{imsref}%
% NOT OUTPUTED:
% issn = 0894-9840
% url = http://dx.doi.org/10.1007/BF02214083
% number = 3
% coden = JTPREO
% fjournal = Journal of Theoretical Probability
\endbibitem

%b51 ###
\bibitem{pengyao03}
%
\begin{barticle}[mr]
\bauthor{\bsnm{Peng},~\bfnm{Liang}\binits{L.}} \AND
\bauthor{\bsnm{Yao},~\bfnm{Qiwei}\binits{Q.}}
(\byear{2003}).
\btitle{Least absolute deviations estimation for {ARCH} and {GARCH} models}.
\bjournal{Biometrika}
\bvolume{90}
\bpages{967--975}.
\bid{doi={10.1093/biomet/90.4.967}, issn={0006-3444}, mr={2024770}}
\end{barticle}
%
\bptok{imsref}%
% NOT OUTPUTED:
% issn = 0006-3444
% url = http://dx.doi.org/10.1093/biomet/90.4.967
% number = 4
% coden = BIOKAX
% fjournal = Biometrika
\endbibitem

%b52 ###
\bibitem{Petersen83}
%
\begin{bbook}[mr]
\bauthor{\bsnm{Petersen},~\bfnm{Karl}\binits{K.}}
(\byear{1983}).
\btitle{Ergodic Theory}.
\bseries{Cambridge Studies in Advanced Mathematics}
\bvolume{2}.
\blocation{Cambridge}:
\bpublisher{Cambridge Univ. Press}.
\bid{mr={0833286}}
\end{bbook}
%
\bptok{imsref}%
% NOT OUTPUTED:
% isbn = 0-521-23632-0
% fpage = xii+329
\endbibitem

%b53 ###
\bibitem{pollard84}
%
\begin{bbook}[mr]
\bauthor{\bsnm{Pollard},~\bfnm{David}\binits{D.}}
(\byear{1984}).
\btitle{Convergence of Stochastic Processes}.
\bseries{Springer Series in Statistics}.
\blocation{New York}:
\bpublisher{Springer}.
\bid{doi={10.1007/978-1-4612-5254-2}, mr={0762984}}
\end{bbook}
%
\bptok{imsref}%
% NOT OUTPUTED:
% isbn = 0-387-90990-7
% url = http://dx.doi.org/10.1007/978-1-4612-5254-2
% fpage = xiv+215
\endbibitem

%b54 ###
\bibitem{Resnick87}
%
\begin{bbook}[mr]
\bauthor{\bsnm{Resnick},~\bfnm{Sidney~I.}\binits{S.I.}}
(\byear{1987}).
\btitle{Extreme Values, Regular Variation, and Point Processes}.
\bseries{Applied Probability. A~Series of the Applied Probability Trust}
\bvolume{4}.
\blocation{New York}:
\bpublisher{Springer}.
\bid{mr={0900810}}
\end{bbook}
%
\bptok{imsref}%
% NOT OUTPUTED:
% isbn = 0-387-96481-9
% fpage = xii+320
\endbibitem

%b55 ###
\bibitem{Rio93}
%
\begin{barticle}[mr]
\bauthor{\bsnm{Rio},~\bfnm{Emmanuel}\binits{E.}}
(\byear{1993}).
\btitle{Covariance inequalities for strongly mixing processes}.
\bjournal{Ann. Inst. Henri Poincar\'e Probab. Stat.}
\bvolume{29}
\bpages{587--597}.
\bid{issn={0246-0203}, mr={1251142}}
\end{barticle}
%
\bptok{imsref}%
% NOT OUTPUTED:
% issn = 0246-0203
% url = http://www.numdam.org/item?id=AIHPB_1993__29_4_587_0
% number = 4
% coden = AHPBAR
% fjournal = Annales de l'Institut Henri Poincar\'e. Probabilit\'es et
%Statistiques
\endbibitem

%b56 ###
\bibitem{ronchettitrojani01}
%
\begin{barticle}[mr]
\bauthor{\bsnm{Ronchetti},~\bfnm{Elvezio}\binits{E.}} \AND
\bauthor{\bsnm{Trojani},~\bfnm{Fabio}\binits{F.}}
(\byear{2001}).
\btitle{Robust inference with {GMM} estimators}.
\bjournal{J. Econometrics}
\bvolume{101}
\bpages{37--69}.
\bid{doi={10.1016/S0304-4076(00)00073-7}, issn={0304-4076}, mr={1805872}}
\end{barticle}
%
\bptok{imsref}%
% NOT OUTPUTED:
% issn = 0304-4076
% url = http://dx.doi.org/10.1016/S0304-4076(00)00073-7
% number = 1
% coden = JECMB6
% fjournal = Journal of Econometrics
\endbibitem

%b57 ###
\bibitem{SakataWhite98}
%
\begin{barticle}[auto:STB|2014/02/12|14:17:21]
\bauthor{\bsnm{Sakata},~\bfnm{S.}\binits{S.}} \AND
\bauthor{\bsnm{White},~\bfnm{H.}\binits{H.}}
(\byear{1998}).
\btitle{High breakdown point conditional dispersion estimation with
application to s\&p 500 daily returns volatility}.
\bjournal{Econometrica}
\bvolume{66}
\bpages{529--567}.
\end{barticle}
%
\bptok{imsref}%
\endbibitem

%b58 ###
\bibitem{Shevlyakovetal08}
%
\begin{barticle}[mr]
\bauthor{\bsnm{Shevlyakov},~\bfnm{Georgy}\binits{G.}},
\bauthor{\bsnm{Morgenthaler},~\bfnm{Stephan}\binits{S.}} \AND
\bauthor{\bsnm{Shurygin},~\bfnm{Alexander}\binits{A.}}
(\byear{2008}).
\btitle{Redescending $M$-estimators}.
\bjournal{J. Statist. Plann. Inference}
\bvolume{138}
\bpages{2906--2917}.
\bid{doi={10.1016/j.jspi.2007.11.008}, issn={0378-3758}, mr={2526216}}
\end{barticle}
%
\bptok{imsref}%
% NOT OUTPUTED:
% issn = 0378-3758
% url = http://dx.doi.org/10.1016/j.jspi.2007.11.008
% number = 10
% coden = JSPIDN
% fjournal = Journal of Statistical Planning and Inference
\endbibitem

%b59 ###
\bibitem{StraumannMikosch}
%
\begin{barticle}[mr]
\bauthor{\bsnm{Straumann},~\bfnm{Daniel}\binits{D.}} \AND
\bauthor{\bsnm{Mikosch},~\bfnm{Thomas}\binits{T.}}
(\byear{2006}).
\btitle{Quasi-maximum-likelihood estimation in conditionally
heteroscedastic time series: A stochastic recurrence equations approach}.
\bjournal{Ann. Statist.}
\bvolume{34}
\bpages{2449--2495}.
\bid{doi={10.1214/009053606000000803}, issn={0090-5364}, mr={2291507}}
\end{barticle}
%
\bptok{imsref}%
% NOT OUTPUTED:
% issn = 0090-5364
% url = http://dx.doi.org/10.1214/009053606000000803
% number = 5
% coden = ASTSC7
% fjournal = The Annals of Statistics
\endbibitem

%b60 ###
\bibitem{ZhuLing}
%
\begin{barticle}[mr]
\bauthor{\bsnm{Zhu},~\bfnm{Ke}\binits{K.}} \AND
\bauthor{\bsnm{Ling},~\bfnm{Shiqing}\binits{S.}}
(\byear{2011}).
\btitle{Global self-weighted and local quasi-maximum exponential
likelihood estimators for {ARMA}--{GARCH}/{IGARCH} models}.
\bjournal{Ann. Statist.}
\bvolume{39}
\bpages{2131--2163}.
\bid{doi={10.1214/11-AOS895}, issn={0090-5364}, mr={2893864}}
\end{barticle}
%
\bptok{imsref}%
% NOT OUTPUTED:
% issn = 0090-5364
% url = http://dx.doi.org/10.1214/11-AOS895
% number = 4
% coden = ASTSC7
% fjournal = The Annals of Statistics
\endbibitem

\end{thebibliography}
\end{document}